\documentclass[leqno,10pt]{article}

\usepackage[usenames, dvipsnames]{xcolor}
\usepackage{tikz}

\usepackage{amssymb,hyperref,amsthm}
\usepackage{euscript}
\usepackage{flafter}
\usepackage{pstricks}
\usepackage[latin1]{inputenc}
\usepackage{amsmath}
\usepackage{amsfonts}
\usepackage{pstricks-add}
\usepackage{dsfont}
\usepackage{bm}
\usepackage{variations}

\hypersetup{
   linktoc=page,
   linkcolor=red,          % color of internal links
  citecolor=blue,        % color of links to bibliography
    filecolor=blue,      % color of file links
   urlcolor=cyan,
    colorlinks=true           % color of external links
}

\usepackage{mathrsfs} %%% Nice caligraphic fonts
\evensidemargin\oddsidemargin
\newcommand{\eqnsection}{
\renewcommand{\theequation}{\thesection.\arabic{equation}}
   \makeatletter
   \csname  @addtoreset\endcsname{equation}{section}
   \makeatother}
\eqnsection
\def\r{{\mathbb R}}
\def\P{{\bf P}}
\def\E{{\bf E}}

\def\N{{\mathbb N}}
\def\1{{\mathds{1}}}
\def\d{\, \mathrm{d}}
\def\d{\mathtt{d}}
\def\o{\mathfrak{O}}

\newtheorem{theorem}{Theorem}[section]

\newtheorem{Definition}[theorem]{Definition}
\newtheorem{Lemma}[theorem]{Lemma}
\newtheorem{Proposition}[theorem]{Proposition}
\newtheorem{Remark}[theorem]{Remark}

\newtheorem{corollary}[theorem]{Corollary}
\newtheorem*{assa}{Assumption (A)}

\title{Glassy phase and freezing of log-correlated Gaussian potentials}
\date{}
\begin{document}

\maketitle
\begin{center}

{ Thomas Madaule \footnotemark[1],\ R\'emi Rhodes \footnotemark[2],\ 
Vincent Vargas \footnotemark[3]}

\bigskip

 \footnotetext[1]{Universit{\'e} Paris-13.} 
\footnotetext[2]{Universit{\'e} Paris-Dauphine, Ceremade, F-75016 Paris, France. Partially supported by grant ANR-11-JCJC  CHAMU}
\footnotetext[3]{Ecole normale sup\'erieure, Paris, France. Partially supported by grant ANR-11-JCJC  CHAMU}

\end{center}

\begin{abstract}
In this paper, we consider the Gibbs measure associated to a logarithmically correlated random potential (including two dimensional free fields) at low temperature. We prove that the energy landscape freezes and enters in the so-called glassy phase. The limiting Gibbs weights are integrated atomic random measures with random intensity expressed in terms of the critical Gaussian multiplicative chaos constructed in \cite{Rnew7,Rnew12}. This could be seen as a first rigorous step in the renormalization theory of super-critical Gaussian multiplicative chaos.   
\end{abstract}
\vspace{0.3cm}
\footnotesize

%\noindent{\bf Short Title.} 

\noindent{\bf Key words or phrases:} Gaussian multiplicative chaos, supercritical, renormalization, freezing, glassy phase, derivative multiplicative chaos.

\medskip
%\noindent{\bf MSC 2000 subject classifications: 60G57, 60G15, 60G25, 28A80}

\normalsize
% \newpage 

\tableofcontents

%%%%%%%%%%%%%%%%%%%%%%%%%%%%%%%%%%%%%%%%%%%%%%%%%%%%%%%%%%%%%%%%%%%%%%%%%%%%%%%%%%%%%%%%%%%%%%%%%%%%%%%%%%%%%%%%%%%%%%%%%%%%%%%%%%%%%%%%%%%%%%%%%%%%%%%%%%%%%%

\section{Introduction}
%%%%%%%%%%%%%%%%%%%%%%%%%%%%%%%%%%%%%%%%%%%%%%%%%%%%%%%%%%%%%%%%%%%%%%%%%%%%%%%%%%%%%%%%%%%%%%%%%%%%%%%%%%%%%%%%%%%%%%%%%%%%%%%%%%%%%%%%%%%%%%%%%%%%%%%%%%%%%%
Consider a log-correlated random distribution $(X(x))_{x\in\r^\d}$ on (a subdomain of) $\r^\d$ and apply a cut-off regularization procedure to get a field $(X_t(x))_{x\in\r^\d}$ with variance of order $t$, i.e. $\E[X_t(x)^2]\simeq t$ as $t\to\infty$. One may for instance think of the convolution of $X$ with a mollifying sequence, the projection of $X$ onto a finite dimensional set of functions or a white noise decomposition of $X$. We will be interested in the study of the behaviour of the random measure on the Borel sets of $\r^\d$:
$$M_t(dx)=e^{\gamma X_t(x)}\,dx,$$
where $\gamma>0$ is a parameter that stands for the inverse temperature. The high temperature phase is well known since the original work of Kahane \cite{cf:Kah} where it is proved that for $\gamma^2<2\d$ the renormalized measure
$$e^{-\frac{\gamma^2}{2}t}M_t(dx)$$ almost surely weakly converges towards a non-trivial measure $M_\gamma(dx)$, which is diffuse. At the critical temperature $\gamma^2=2\d$, the renormalized measure
$$\sqrt{t}e^{-\d t}M_t(dx)$$ weakly converges in probability towards a non-trivial diffuse measure $M'(dx)$, which is  called derivative multiplicative chaos \cite{Rnew7,Rnew12}. The purpose of this paper is to study the supercritical/low temperature phase $\gamma^2>2\d$ and to prove that 
the renormalized measure
\begin{equation}\label{intro:surcrit}
t^{\frac{3\gamma}{2\sqrt{2\d}}}e^{(\gamma\sqrt{2\d}-\d)t}M_t(dx)
\end{equation}
 weakly converges in law towards a purely atomic stable random measure $S_\gamma$ with (random intensity) $M'$, up to a deterministic multiplicative constant, call it $C(\gamma)$ (see section \ref{setup} for a rigorous statement). 
  
This is a longstanding problem, which has received much attention by physicists. It was first raised in \cite{mandelbrotstar,derrida} on dyadic trees, and then followed by \cite{CarDou,Fyo,rosso} for log-correlated Gaussian random fields. Following our notations, these papers essentially derived the statistics of the size ordered atoms of the measure $\frac{M_t(dx)}{M_t([0,1]^\d)}$, the so-called Poisson-Dirichlet statistics characteristic of stable Levy processes. However, these papers did not investigate the problem of the localization of these atoms. 

A few years later, the mathematical community caught up on this problem. In the context of Branching random walks, convergence of the measures \eqref{intro:surcrit} is investigated in \cite{madaule,Webb}. Built on these works,  the limit is identified in \cite{Rnew3} and is expressed as a stable transform of the so-called derivative martingale. In the context of log-correlated Gaussian potentials, the authors  in  \cite{Rnew4}  conjecture that results similar to Branching Random Walks should hold. The first  rigorous and important result for log-correlated Gaussian fields appeared in \cite{arguin} where the authors established the Poisson-Dirichlet statistics of the limiting measure in dimension $1$ (renormalized by its total mass) via spin glass techniques, hence confirming the prediction of \cite{CarDou} (these results were recently extended by the same authors in \cite{arguin2} to cover the case of the discrete GFF in a bounded domain).    

%in a way, this is a study of the masses of the atoms of $S_\gamma$, without focusing on their spatial localization.   

Roughly speaking, the terminology freezing comes from the linearization of the free energy of the measure $M_t$   beyond the value $\gamma^2=2\d$ (see \cite{derrida,CarDou,Fyo,rosso} for further comments). The terminology {\it glassy phase} comes from the fact that for $\gamma^2>2\d$, the measure $M_n$ is essentially dominated by a few points, the local extreme values of the field $X_t$ (along with the neighborhood of these extreme values). Therefore, this paper possesses strong connections with the study of the extreme values of the field $X_t$. This was conjectured in \cite{Rnew7} and important advances on this topic have recently appeared in \cite{Louisdor,DZ} in the context of the discrete GFF and in \cite{Mad13} for a large class of log-correlated fields. However, the description of the local maxima obtained in \cite{Louisdor} is not sufficient to obtain the so-called freezing theorems that will be established in this paper.   
%and these connections are further discussed in section \ref{setup}.

Finally, we would like to stress that we will only deal with the case of white noise cut-off of the Gaussian distribution $X$, building on techniques developed in \cite{Mad13}. We will then extend our results to two dimensional free fields. It is natural to wonder whether the nature of the cut-off may affect the structure  of the limiting measure. We will prove that the freezing theorem does not depend on the chosen cutoff family provided the cutoff is not too far from a white noise decomposition. From a more general angle, we believe that the glassy phase does not depend on the chosen cut-off, except  at the level of the multiplicative constant $C(\gamma)$. For instance, given a smooth mollifier $\theta$ and setting $\theta_\epsilon= \frac{1}{\epsilon^\d}\theta (\frac{.}{\epsilon})$, similar theorems should hold for measures built on approximations of the form $\theta_\epsilon \ast X$: in this setting, one would obtain an analog of theorem \ref{main} where the constant $C(\gamma)$ is replaced by a constant $C(\theta)$ depending on $\theta$.

%It would in fact be interesting to generalize the theory of Gaussian multiplicative chaos to the supercritical case. More precisely, if $\theta$ is smooth mo  

%We explain heuristically in section \ref{setup} why this should be true in the context of the two dimensional discrete Gaussian Free Field.
%%%%%%%%%%%%%%%%%%%%%%%%%%%%%%%%%%%%%%%%%%%%%%%%%%%%%%%%%%%%%%%%%%%%%%%%%%%%%%%%%%%%%%%%%%%%%%%%%%%%%%%%%%%%%%%%%%%%%%%%%%%%%%%%%%%%%%%%%%%%%%%%%%%%%%%%%%%%%%
\section{Setup and main results}\label{setup}
%%%%%%%%%%%%%%%%%%%%%%%%%%%%%%%%%%%%%%%%%%%%%%%%%%%%%%%%%%%%%%%%%%%%%%%%%%%%%%%%%%%%%%%%%%%%%%%%%%%%%%%%%%%%%%%%%%%%%%%%%%%%%%%%%%%%%%%%%%%%%%%%%%%%%%%%%%%%%%
%\subsection{Notations}
%%%%%%%%%%%%%%%%%%%%%%%%%%%%%%%%%%%%%%%%%%%%%%%%%

\subsection{Star scale invariant fields}\label{setupssik}
%%%%%%%%%%%%%%%%%%%%%%%%%%%%%%%%%%%%%%%%%%%%%%%%%

We denote by $\mathcal{B}_b(\r^\d)$ the Borel subsets of $\r^\d$. Let us introduce a canonical family of log-correlated Gaussian distributions, called star scale invariant, and their cut-off approximations, which we will work with in the first part of the paper.   Let us consider a continuous covariance kernel $k$ on $\r^{\d}$ such that:
\begin{assa}  The kernel $k$  satisfies  the following assumptions, for some constant  $C$ independent of   $x\in\r^\d$:
\begin{description}
\item[A1.] $k$ is continuous, nonnegative  and normalized by the condition $k(0)=1$,
\item[A2.] $k$ has compact support.
\item[A3.] $|k(x)-k(0)|\leq C|x|  $ for some constant $C:=C_k$ independent of   $x\in\r^{\d}$.
\end{description}
\end{assa}

We set for $t\geq 0$ and $x\in\r^{\d}$
 \begin{equation}\label{def:sskernel}
K_t(x) =  \int_1^{e^t}  \frac{k(xu)}{u}du.
 \end{equation}
We consider  a family of centered Gaussian processes $(X_t(x))_{x\in\r^{\d},t\geq 0}$  with covariance kernel given by:
 \begin{equation}
\forall t,s\geq 0,\quad  \E[ X_t(x) X_s(y)   ] 
= K_{t \wedge s}  (y-x ),
 \end{equation}
where $t \wedge s:= \min (t,s)$. The construction of such fields is possible via a white noise decomposition as explained in \cite{Rnew1}. We  set: 
\begin{align*}%
\mathcal{F}_{t} = \sigma \lbrace X_u(x);x\in\r^{\d},u \leq t\rbrace.\quad% & \text{and}\quad  \mathcal{F}^X   = \sigma \lbrace X_u(x);x\in\r^d,u \in ]0,1]\rbrace   .
\end{align*} 
 We stress that, for $s>t $, the field $(X_s(x)-X_t(x))_{x \in \r^{\d}}$ is independent from $\mathcal{F}_{t}$. 

We introduce for $t>0$ and $\gamma>0$, the random measures $M_t'(dx)$ and $M_t^\gamma(dx)$
\begin{equation}
\label{deuxdef}
M_t'(A):= \int_A( \sqrt{2\d}t-X_t(x))e^{\sqrt{2\d}X_t(x)-\d t}dx,\quad M_t^\gamma(A):= \int_Ae^{\gamma X_t(x)-\frac{\gamma^2}{2}t}dx,\quad \forall A\in \mathcal{B}_b(\r^\d).
\end{equation}
 Recall that (see \cite{Rnew7})   %,Rnew12})
\begin{theorem}\label{convderiv}
For each bounded open set $A\subset \r^\d$, the martingale $(M_t'(A))_{t\geq 0}$ converges almost surely towards a positive random variable denoted by $M'(A)$. 

Furthermore, the family of random measures $(M_t'(dx))_{t\geq 0}$ almost surely weakly converges towards a random measure $M'(dx)$, which is atom free and has full support.
%The variable 
%$(\sqrt{t} M_t^{\sqrt{2d}}(A))_{t\geq 0}$ also converges in probability towards $\sqrt{\frac{2}{\pi}} M'(A)$.
\end{theorem}
 
 In fact, one can reinforce the above statement to a convergence in the space of Radon measures.

\subsection{Results for star scale invariant fields}
%%%%%%%%%%%%%%%%%%%%%%%%%%%%%%%%%%%%%%%%%%%%%%%%%
The main purpose of this paper is to establish the following result which was conjectured in \cite{Rnew7}:
\begin{theorem}\label{main}{\bf (Freezing theorem)}
\label{Lapl}
For any $\gamma>\sqrt{2\d}$, there exists a constant $C(\gamma)>0$ such that for any smooth nonnegative function $f$ on $[0,1]^\d$ 
\begin{equation}
\label{eqLapl}\underset{t\to\infty}{\lim}  \E\left( \exp(- t^{\frac{3\gamma}{2\sqrt{2\d}}}e^{t(\frac{\gamma}{\sqrt{2}}-\sqrt{\d})^2}\int_{[0,1]^\d} f(x)M_t^\gamma(dx))\right)=\E\left(\exp(- C (\gamma) \int_{[0,1]^\d} f(x)^\frac{\sqrt{2\d}}{\gamma}M'(dx))\right).
\end{equation}
\end{theorem}
 
As a consequence, we deduce
\begin{corollary} \label{coro:meas}
For any $\gamma>\sqrt{2\d}$, the family of random measures $(  t^{\frac{3\gamma}{2\sqrt{2\d}}}e^{t(\frac{\gamma}{\sqrt{2}}-\sqrt{\d})^2}M_t^\gamma(dx))_{t\geq 0}$ weakly converges in law towards a purely atomic random measure denoted by $S_\gamma (dx)$. The law of $S_\gamma$ can be described as follows: conditionally on $M'$, $S_\gamma$ is an independently scattered random measure such that 
\begin{equation}
 \E\left( \exp(-\theta S_\gamma(A))\right)=\E\left(\exp(-\theta^\frac{\sqrt{2\d}}{\gamma}C (\gamma) M'(A))\right)
\end{equation} for all $\theta\geq 0$.
\end{corollary}
Put in other words, $S_\gamma$ is an integrated $\alpha $-stable Poisson random measure of spatial intensity given by the derivative martingale $M'$. Indeed, the law of $S_\gamma$ may be described as follows. Conditionally on $M'$, consider a Poisson random measure $n_\gamma$ on $\r^\d\times \r_+$ with intensity
$$M'(dx)\otimes \frac{dz}{z^{1+\frac{\sqrt{2\d}}{\gamma}}}.$$
Then the law of $S_\gamma$ is the same as the purely atomic measure ($\Gamma$ stands for the  function gamma)
$$S_\gamma(A)=c\int_A\int_0^\infty z \,n_\gamma(dx,dz)\quad \text{ with }\quad c  =\Big(C (\gamma)\frac{\sqrt{2\d}}{\gamma \Gamma(1-\frac{\sqrt{2\d}}{\gamma})}\Big)^{\frac{\gamma}{\sqrt{2\d}}}.$$

From Theorem \ref{convderiv},  we observe that $M'(\mathcal{O})>0$ almost surely for any open set. By considering this together with Corollary \ref{coro:meas}, it is plain to deduce
\begin{corollary}\label{dirichlet}
For each bounded open set $\mathcal{O}$, the family of random measures $\Big(\frac{M_t(dx\cap \mathcal{O})}{M_t(\mathcal{O})}\Big)_{t}$ converges in law in the sense of weak convergence of measures towards $\frac{S_\gamma(dx)}{S_\gamma(\mathcal{O})}$.
\end{corollary} 

We point out that the size reordered atoms of the measure $\frac{S_\gamma(dx)}{S_\gamma(\mathcal{O})}$ form the Poisson-Dirichlet process studied in \cite{arguin,arguin2}. The interesting point here is that we keep track of the spatial localization of the atoms whereas all this information is lost in the   Poisson-Dirichlet approach. Yet, we stress that the methods used in \cite{arguin,arguin2} rely on spin glass technics and remain thus quite interesting since far different from those used here. 

\begin{Remark} We stress that 
Corollary \ref{dirichlet} also holds for all the examples described below but we will refrain from stating it anymore.
\end{Remark}
 \subsection{Massive Free Field}
%%%%%%%%%%%%%%%%%%%%%%%%%%%%%%%%
In this section, we   extend our  results (Theorem \ref{main}) to kernels with long range correlations, in particular, we will be interested in  the whole plane Massive Free Field (MFF). 
 
The whole plane MFF is a centered Gaussian distribution with covariance kernel given by the Green function of the operator $2\pi (m^2-\triangle)^{-1}$ on $\r^2$, i.e. by:
\begin{equation}\label{MFFkernel}
\forall x,y \in \r^2,\quad G_m(x,y)=\int_0^{\infty}e^{-\frac{m^2}{2}u-\frac{|x-y|^2}{2u}}\frac{du}{2 u}.
\end{equation}
  The real $m>0$ is called the mass. This kernel is of $\sigma$-positive type in the sense of Kahane \cite{cf:Kah} since we integrate a continuous function of positive type with respect to a positive measure. It is furthermore  a star-scale invariant kernel (see \cite{Rnew1}): it can be rewritten as 
\begin{equation}\label{MFF3}
G_m(x,y)=\int_{1}^{+\infty}\frac{k_m(u(x-y))}{u}\,du.
\end{equation}
 for some continuous covariance kernel $k_m=\frac{1}{2}\int_0^\infty e^{-\frac{m^2}{2v}|z|^2-\frac{v}{2}}\,dv$. We consider  a family of centered Gaussian processes $(X_t(x))_{x\in\r^{\d},t\geq 0}$  with covariance kernel given by:
 \begin{equation}
\forall t,s\geq 0,\quad  \E[ X_t(x) X_s(y)   ] 
= G_{m,t \wedge s}  (y-x ),
 \end{equation}

One can construct the derivative martingale $M'$ associated to $(X_t)_{t\geq 0}$ as prescribed in \cite{Rnew7,Rnew12}.   
Now we claim that our result holds in the case of the MFF for any  cut-off family of the MFF uniformly close to $(G_{m,t})_t$: 
\begin{Definition}
A   cut-off family of the MFF is said uniformly close to $(G_{m,t})_t$ if it is a family of stochastically continuous centered Gaussian processes $(X_n(x))_{n\in\N,x\in\r^2}$ with respective covariance kernels $(K_n)_n$ satisfying: \\
- we can find a subsequence $(t_n)_n$ such that
$\lim_{n\to\infty}t_n=+\infty$,\\
- the family $(K_n-G_{m,t_n})_n$ uniformly converges towards $0$ over the compact subsets of $\r^2$.  
\end{Definition}

 Then we claim: 
\begin{theorem}{\bf (Freezing theorem for MFF.)}\label{mainMFF}
For any $\gamma>2$, there exists a constant $C (\gamma)>0$ such that for every   cut-off family $(X_n)_n$ of the MFF uniformly close to $(G_{m,t})_t$, the family of random measures $(  t_n^{\frac{3\gamma}{4}}e^{t_n(\frac{\gamma}{\sqrt{2}}-\sqrt{2})^2}M_n^\gamma(dx))_{t\geq 0}$, where
$$M_n^\gamma(dx)=e^{\gamma X_n(x)-\frac{\gamma^2}{2}\E[X_n(x)^2]}\,dx,$$  
%and for any continuous nonnegative function $f$ with compact support:
% \begin{equation}
%\label{eqLaplMFF}\underset{n\to\infty}{\lim} \E\left( \exp(- t_n^{\frac{3\gamma}{4}}e^{t_n(\frac{\gamma}{\sqrt{2}}-\sqrt{2})^2}M_n^\gamma(f))\right)=\E\left(\exp(-C (\gamma) \int_{\r^2}f(x)^\frac{2}{\gamma}M'(dx))\right)
%\end{equation}
%where
%$$M_n^\gamma(dx)=e^{\gamma X_n(x)-\frac{\gamma^2}{2}\E[X_n(x)^2]}\,dx.$$   
%Also, the family of random measures $(  t_n^{\frac{3\gamma}{2\sqrt{2\d}}}e^{t_n(\frac{\gamma}{\sqrt{2}}-\sqrt{\d})^2}M_n^\gamma(dx))_{t\geq 0}$ 
weakly converges in law towards a purely atomic random measure denoted by $S_\gamma $. The law of $S_\gamma$ can be described as follows: 
%conditionally on $M'$, $S_\gamma$ is an independently scattered random measure such that 
\begin{equation}
 \E\left( \exp(-  S_\gamma(f))\right)=\E\left(\exp\big(-C (\gamma) \int_{\r^2}f(x)^\frac{2}{\gamma}M'(dx)\big)\right)
\end{equation} for all nonnegative continuous  function $f$ with compact support.
\end{theorem}

The above theorem is a bit flexible in the sense that there is some robustness with respect to the chosen cutoff approximation: among the class of cut-off families of the MFF  uniformly close to $(G_{m,t})_t$, the freezing phenomena related to the MFF  do not depend on the structure of the chosen cutoff.   
%For the time being, we believe but cannot prove that this theorem is much more universal regarding the chosen cut-off family of the MFF than the above statement 
%(see the discussion in subsection \ref{discussion}).  
 
\subsection{Gaussian Free Field on planar bounded domains}
%%%%%%%%%%%%%%%%%%%%%%%%%%%%%%%%
Consider a bounded open domain $D$ of $\r^2$. Formally, a GFF on $D$ is a Gaussian distribution with covariance kernel given by the Green function of the Laplacian on $D$ with prescribed boundary conditions (see \cite{She07} for further details). We describe here the case of Dirichlet boundary conditions. The Green function is then given by the formula:
\begin{equation}\label{bounGreen}
G_D (x,y)=  \pi \int_{0}^{\infty}p_D(t,x,y)dt
 \end{equation}
where $p_D$  is the (sub-Markovian) semi-group of a Brownian motion $B$ killed upon touching the boundary of $D$, namely
$$p_D(t,x,y)=P^x(B_{t} \in dy, \; T_D > t)/dy$$ with $T_D=\text{inf} \{t \geq 0, \; B_t\not \in D \}$. Note  the factor $\pi$, which makes sure   that $G_D(x,y)$ takes on the form 
$$G_D(x,y) =\ln_+\frac{1}{|x-y|}+g(x,y)$$ where $\ln_+= \max(\ln,0)$ and for some continuous function $g$ on $D\times D$. The most  direct way to construct a cutoff family of the GFF on $D$ is then to consider  a white noise $W$ distributed on $D\times \r_+$ and define:
 $$X(x)=\sqrt{\pi}\int_{D\times \r_+}p_D(\frac{s}{2},x,y)W(dy,ds).$$
One can check that $\E[X(x)X(x')]=\pi\int_0^{\infty} p_D(s,x,x')\,ds=G_D(x,x') $. The corresponding cut-off approximations are given by:
\begin{equation}\label{defxt}
X_t(x)=\sqrt{\pi}\int_{D\times [e^{-2t},\infty[}p_D(\frac{s}{2},x,y)W(dy,ds),
\end{equation}
which has covariance kernel
$$G_{D,t}(x,y)=\pi \int_{e^{-2t}}^{\infty}p_D(r,x,y)dr.$$
We define the approximating measures
\begin{equation*}
M_t^2(dx)=  e^{2 X_t(x)-2\E[X_t(x)^2]}\,dx \quad\quad \text{ and }\quad \quad M_t'(dx)= (2\E[X_t(x)^2]-X_t(x))e^{2 X_t(x)-2\E[X_t(x)^2]}\,dx. 
\end{equation*}
Let us stress that Theorem \ref{convderiv} holds for this family $(X_t)_t$ (see \cite{Rnew12}).  

\begin{theorem}{\bf (Freezing theorem for GFF on planar domains.)}\label{mainGFF}
For any $\gamma>2$ and every bounded planar domain $D\subset \r^2$, there exists a constant $C (\gamma)>0$ such that for every   cut-off family $(X_n)_n$ of the GFF uniformly close to $(G_{D,t})_t$, the family of random measures $(  t_n^{\frac{3\gamma}{4}}e^{t_n(\frac{\gamma}{\sqrt{2}}-\sqrt{2})^2}M_n^\gamma(dx))_{t\geq 0}$, where 
$$M_n^\gamma(dx)=e^{\gamma X_n(x)-\frac{\gamma^2}{2}t_n}\,dx,$$ weakly converges in law towards a purely atomic random measure denoted by $S_\gamma $.
%and for any continuous nonnegative function $f$ with compact support in $D$:
% \begin{equation}
%\label{eqLaplMFF}\underset{n\to\infty}{\lim} \E\left( \exp(- t_n^{\frac{3\gamma}{4}}e^{t_n(\frac{\gamma}{\sqrt{2}}-\sqrt{2})^2}M_n^\gamma(f))\right)=\E\left(\exp(-C (\gamma) \int_{\r^2}f(x)^\frac{2}{\gamma}C(x,D)^2M'(dx))\right)
%\end{equation}
%$$M_n^\gamma(dx)=e^{\gamma X_n(x)-\frac{\gamma^2}{2}t_n}\,dx.$$   
%Also, the family of random measures $(  t_n^{\frac{3\gamma}{2\sqrt{2\d}}}e^{t_n(\frac{\gamma}{\sqrt{2}}-\sqrt{\d})^2}M_n^\gamma(dx))_{t\geq 0}$ weakly converges in law towards a purely atomic random measure denoted by $S_\gamma $. 
The law of $S_\gamma$ can be described as follows:  
\begin{equation}\label{enoncetheorem}
 \E\left( \exp(-  S_\gamma(f))\right)=\E\left(\exp\big(-C (\gamma) \int_{\r^2}f(x)^\frac{2}{\gamma}C(x,D)^2M'(dx)\big)\right)
\end{equation} for all nonnegative continuous  function $f$ with compact support, where $C(x,D)$ stands for the conformal radius at $x\in D$.
\end{theorem}

\begin{Remark}
The derivative martingale construction of theorem \ref{convderiv} applies to other cut-offs of the GFF than \eqref{defxt}. For instance, one can consider the projection of the GFF on the triangular lattice with mesh going to $0$ along powers of 2 (in this case, the law on the lattice points of this projection is nothing but the discrete GFF on the triangular lattice). Then, the derivative martingale construction of theorem \ref{convderiv} holds in this context by the methods of \cite{Rnew7} since the approximations correspond to adding independent functions: see \cite{She07}. 
%By the universality results in \cite{review}, the approximation \eqref{defxt} and the projection on the triangular lattice yield the same critical %measure $M'$ (in law). 
Unfortunately, the methods of this paper do not enable to prove an analog of theorem \ref{mainGFF} in the context of the projection on the triangular lattice. There are several difficulties to overcome in this context. First, it would be interesting to prove that Seneta-Heyde renormalization of \cite{Rnew12} yields the same limit as the derivative martingale in this setting (this is not obvious from the techniques of \cite{Rnew12}). By the universality results in \cite{review}, this would imply that the approximation \eqref{defxt} and the projection on the triangular lattice yield the same critical measure $M'$ (in law). Proving an analog of theorem \ref{mainGFF} for the the triangular lattice would then imply by the above discussion that the renormalized supercritical measures with the triangular lattice cut-off converge in law to the $S_\gamma$ defined in \eqref{enoncetheorem} (up to some multiplicative constant).         

\end{Remark}

\subsection{Further generalization} 
%%%%%%%%%%%%%%%%%%%%%%%%%%%%%%%%

Our strategy of proofs apply to a more general class of kernels, at least to some extent, in any dimension. There are two main inputs to take care of. 

First we discuss the case of long range correlated star scale invariant kernels. One has to adopt the same strategy as we do for the MFF. Basically, what one really needs is   assumptions [B.1]+[B.2]+[B.3] and the Seneta-Heyde norming, whatever the dimension. However, further conditions on the kernel $k$  are required in order   to make sure that the Seneta-Heyde norming holds (see \cite[Remark 31]{Rnew12}). One may for instance treat in this way the case of covariance kernel given by the Green function of the operator $(m^2-\triangle)^{\d/2}$ in $\r^{\d}$ provided that $m>0$.

One may then wish to treat the case of non translation invariant fields, for instance with correlations given by the Green function of $(-\triangle)^{\d/2}$ in a bounded domain of $\r^\d$ with appropriate boundary conditions. Then one has to adopt the strategy we use for the GFF on planar domains: just replace the conformal radius by the function
$$F(x,D)= \lim_{t\to\infty}e^{\E[X_t(x)^2]-\ln t}.$$

\section{Proofs for star scale invariant fields}
%%%%%%%%%%%%%%%%%%%%%%%%%%%%%%%%%%%%%%%%%%%%%%%%%%%%%%%%%%%%%%%%%%%%%%%%%%%%%%%%%%%%%%%%%%%%%%%%%%%%%%%%%%%%%%%%%%%%
In this section, we carry out the main arguments of the proof of Theorem \ref{main}. Throughout the proofs we will use results that are gathered in a  toolbox   in Appendix \ref{app:toolbox}. Furthermore, from Assumption ${\bf A2}$, the covariance kernel $k$ has compact support. Without loss of generality, we will assume that the support of $k$ is contained in the ball centered at $0$ with radius $1$.
%%%%%%%%%%%%%%%%%%%
\subsection{Some further notations}
%%%%%%%%%%%%%%%%%%% 
We first introduce some further notations and gather some results that we will use and that can be found in the literature.  

\subsubsection*{Processes and measures}
%%%%%%%%%%%%%%%%%%%%%%%%%%%%%%%%%%
Before proceeding with the proof, we introduce some further notations. We   define for all $x\in\r^\d$, $l>0$, $t\geq 0$, all Borelian subset $A$ of $\r^{\d}$:
\begin{equation}\label{not:proof}
Y_t(x):=  X_t(x)-\sqrt{2\d}\,t\quad\text{ and }\quad  Y_{s,t}(x):=Y_{s+t}(x)-Y_s(x).
 \end{equation}
We recall the following scaling property:
\begin{equation}
\label{scalingl}
\left(Y_{s,t}(x)\right)_{t\in \r^+,\,x\in \r^\d}\overset{(law)}{=} (Y_t(xe^{s}))_{t\in \r^+,\,x\in \r^\d},
\end{equation} 
which can be checked with a straightforward computation of covariances. This scaling property is related to the notion of star scale invariance and the reader is referred to  \cite{Rnew1} for more on this.

The main purpose of Theorem \ref{main} will be to establish the convergence of the renormalized measure $t^{\frac{3\gamma}{2\sqrt{2\d}}}e^{t(\frac{\gamma}{\sqrt{2}}-\sqrt{\d})^2}M_t^\gamma(dx)$ and it will thus be convenient to shortcut this expression as:
\begin{equation}\label{defmtilde}
\tilde{M}_t^\gamma(dx):= t^{\frac{3\gamma}{2\sqrt{2\d}}}e^{t(\frac{\gamma}{\sqrt{2}}-\sqrt{\d})^2}M_t^\gamma(dx).
\end{equation}

We will denote by $|A|$ the Lebesgue measure of a measurable set $A\subset\r^\d$.

\subsubsection*{Regularity, spaces of functions}
%%%%%%%%%%%%%%%%%%%%%%%%%%%
We denote by $\mathcal{C}(B,\r^p)$ the space of continuous functions from $B$ (a subset of $\r^q$) into $\r^p$. 

For any continuous function  $f\in \mathcal{C}([0,R]^\d,\r)$ and $\delta>0$, we consider the two following modulus of regularity of $f$
$$w_f(\delta):=\sup_{x,y\in[0,R]^\d,|x-y|\leq \delta}|f(x)-f(y)|\quad \text{ and }\quad w_f^{1/3}(\delta):=\sup_{x,y\in[0,R]^\d,|x-y|\leq \delta}\frac{|f(x)-f(y)|}{|x-y|^{1/3}}.  $$
For any $a,b,t,R>0 $, we define
\begin{equation}
\mathcal{C}_R(t,a,b)=\left\{f:\r^\d\to\r;\,\, w_f^{1/3}(t^{-1})\leq 1,\min_{y\in [0,R]^{\d}}f(y)>a, \quad \text{ and }\max_{y\in [0,R]^{\d}}f(y)<b\right\}.
\end{equation}

\subsubsection*{Constants}
%%%%%%%%%%%%%%%%%%%%%%%%%%%
We also set for $z,t\geq 0$
\begin{equation}
\kappa_\d=\frac{1}{8\sqrt{2\d}}\quad\quad \quad \quad\text{ and }\quad\quad\quad\quad a_t:=-\frac{3}{2\sqrt{2\d}}\ln t
\end{equation}

%\subsubsection*{a virer si possible}
%%%%%%%%%%%%%%%%%%%%%%%%%%%

%
%
% 
%  
%
%$$I_t(z):=[a_t+z-1,a_t+z]$$
%
%For any $g\in \mathcal{C}(\r_+\times [0,R]^\d,\r)$,  $\delta>0$, $y\in [0,R]^\d$ and $t>0$, let:
%$$w_g(\delta,y,t):=\sup_{s\leq t, x\in[0,R]^\d,|x-y|\leq \delta}|g_s(x)-g_s(y)|. $$
%
%For any $g\in \mathcal{C}(\r_+\times [0,R]^\d,\r)$, $[a,b] \subset \r_+$ with $a<b$ and $x \in [0,R]^\d$, let:
%$$  \overline{g}_{[a,b]}(x)=  \sup_{s \in [a,b]} g_s(x), $$
%where we will abbreviate $\overline{g}_{b}:= \overline{g}_{[0,b]}$.
%

\subsection{A decomposition}\label{sub:decomp}
%%%%%%%%%%%%%%%%%%  

\begin{figure}[t] 
\label{croquis}
\begin{center}
\begin{tikzpicture}[xscale=1.3,yscale=1.3] 
\tikzstyle{sommet}=[circle,draw,fill=yellow,scale=0.6]; 
\fill[Red] (0,0) rectangle (5.7,5.7);
\foreach \x in {0,...,5}{
    \foreach \j in {0,...,5}
        \fill[ProcessBlue] (\x,\j) rectangle (\x+0.7,\j+0.7);
}
\foreach \x in {0,...,5}{
    \foreach \j in {0,...,5}
       \node[sommet]  at (\x,\j){};
}
\node[draw=red!70!black,thick,rounded corners=5pt,top color=blue!20,bottom color=blue!10] (Y) at (-3,5) {Points $a_{{\bf i}}$ (yellow)};
\node[draw=red!70!black,thick,rounded corners=5pt,top color=blue!20,bottom color=blue!10] (B) at (-3,3) {Small squares $A_{{\bf i}}$ (blue)};
\node[draw=red!70!black,thick,rounded corners=5pt,top color=blue!20,bottom color=blue!10] (R) at (-3,1) {Buffer zone (red)};
\draw[->,>=latex] (Y) to[bend left] (-0.1,5);
\draw[->,>=latex] (Y) to[bend right] (-0.1,4);
\draw[->,>=latex] (B) to[bend right] (0.2,1.5);
\draw[->,>=latex] (B) to[bend right] (0.2,2.5);
\draw[->,>=latex] (R) to[bend right] (0.2,0.85);
\draw[<->,>=latex] (6,1.7) -- (6,2) node[midway,right] {width $1$};
\draw[<->,>=latex] (6,4) -- (6,4.7) node[midway,right] {width $R$};
\node (C) at (3,-0.4) {$m$ blue squares};
\draw[->,>=latex] (C) -- (5.7,-0.4);
\draw[->,>=latex] (C) -- (0,-0.4);
\end{tikzpicture}
\caption{Decomposition of the cube $[0,e^{t'}]^d$}
\end{center}
\end{figure}
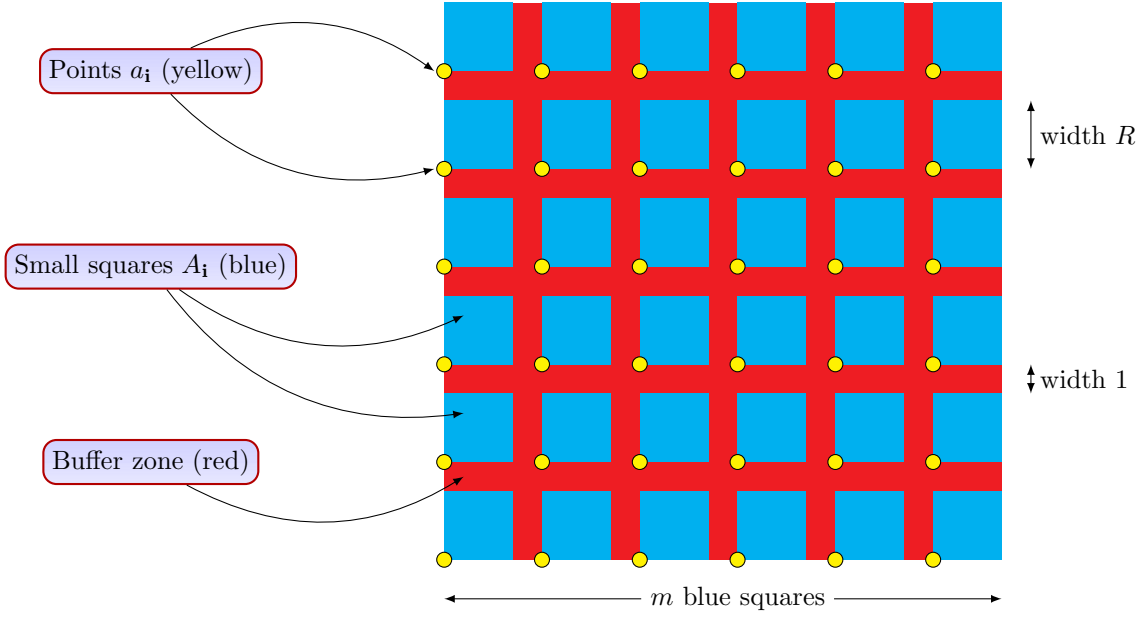

Before proceeding with the proof of Theorem \ref{main}, we first explain a decomposition of the cube $[0,e^{t'}]^\d$ with $t'>0$ that will be used throughout the proof of Theorem \ref{main}. We will divide this cube into several smaller cubes of size $R>0$, all of these smaller cube being at distance greater than $1$  from each other. To understand more easily our notations, the reader may keep in mind the picture of  Figure \ref{croquis}. 
We assume that $R,t'$ are such that
 $$m:=\frac{e^{t'}+1}{(R+1)}\in \N^*.$$  
The integer $m$ stands for the number of small squares of size $R$ that one meets along an edge of the cube. The basis of each small square will be indexed with a $\d$-uplet
$${\bf i}= (i_1,...,i_\d)  \in\{1,...,m\}^\d    .$$ The basis of the square $A_{{\bf i}}$ is then located at 
$$a_{{\bf i}}:= (R+1)\left((i_1-1),...,(i_d-1)\right)\in [0,e^{t'}]^\d$$
in such a way that  
$$A_{{\bf i}}:= a_{\bf i}+ [0,R]^\d.$$
One may observe on Figure \ref{croquis} that all the squares $(A_{{\bf i}}$ are separated from each other by a fishnet shaped buffer zone (red), which is precisely
$${\rm BZ}_{R,t'}:=[0,e^{t'}]^\d\setminus \bigcup_{{\bf i}\in \{1,...,m\}^\d} A_{{\bf i}}.$$

The terminology "buffer zone" is used because this is the minimal area needed to make sure that the values taken by the process $Y_t$ inside each (blue square) $A_{\bf i}$ are independent of its values on all other $A_{\bf j}$ for ${\bf j}\neq {\bf i}$.

\subsection{Main frame of the proof of Theorem \ref{main}}
%%%%%%%%%%%%%%%%%%%%%%%%%%%%%%%%%  

This subsection is devoted to the proof of Theorem \ref{main} up to admitting a few auxiliary results, which will be proved later. 

We fix $\epsilon>0$ and $\theta>0$. For $R>0$ and $t'>0$ such that $\frac{e^{t'}+1}{R+1}\in\N^*$, we define the set (recall the decomposition in subsection \ref{sub:decomp})
\begin{eqnarray}
\nonumber\mathcal{Y}_{R,\theta}(t'):=\big\{    w_{Y_{t'}(\cdot)}^{1/3}(\frac{1}{t'}e^{-t'})\leq e^{\frac{t'}{3}},\,  |\gamma^{-1}\ln \theta|M_{t'}^{\sqrt{2\d}}([0,1]^\d)+ |M_{t'}'(e^{-t'}{\rm BZ}_{R,t'})|  \leq \epsilon \theta^{-\frac{\sqrt{2\d}}{\gamma}},\,\qquad
\\
\label{defYRz}  \forall x\in [0,1]^\d ,\,-10\sqrt{2\d}t'\leq Y_{t'}(x)\leq - \kappa_\d \ln t'\}.  
\end{eqnarray} 

Now we consider $t,t'$ such that $t\geq e^{t'}$. We have
\begin{eqnarray}
\label{firstinequ} &&\E\left(e^{-\theta \tilde{M}_t^\gamma([0,1]^\d) } ; \mathcal{Y}_{R,\theta}(t') \right)  \leq \E\left(e^{-\theta \tilde{M}_t^\gamma([0,1]^\d) } \right) \leq   \E\left(e^{-\theta \tilde{M}_t^\gamma([0,1]^\d) } ; \mathcal{Y}_{R,\theta}(t') \right) +\P(\mathcal{Y}_{R,\theta}(t')^c).
\end{eqnarray}
We estimate now the left-hand side of this relation. Because
$$\tilde{M}_t^\gamma([0,1]^\d)=\tilde{M}_t^\gamma(e^{-t'}{\rm BZ}_{R,t'}) +\tilde{M}_t^\gamma(e^{-t'} \cup A_{\bf i}) ,$$
we can use the relation  $uv\geq u+v-1$ for $u,v\in[0,1]$ to get $$e^{-\theta \tilde{M}_t^\gamma([0,1]^\d) } = e^{-\theta \tilde{M}_t^\gamma(e^{-t'}{\rm BZ}_{R,t'}) }-1+ e^{-\theta \tilde{M}_t^\gamma(e^{-t'} \cup A_{\bf i}) } .$$  
We deduce from \eqref{firstinequ}
\begin{multline}\label{secineq}
 \E\left(e^{-\theta \tilde{M}_t^\gamma(e^{-t'}{\rm BZ}_{R,t'}) }   -1; \mathcal{Y}_{R,\theta}(t')\right) + \E\left(e^{-\theta \tilde{M}_t^\gamma(e^{-t'} \cup A_{\bf i}) } ; \mathcal{Y}_{R,\theta}(t') \right)   \\
 \leq  \E\left(e^{-\theta \tilde{M}_t^\gamma([0,1]^\d) } \right)\\
 \leq   \E\left(e^{-\theta \tilde{M}_t^\gamma(e^{-t'} \cup A_{\bf i}) } ; \mathcal{Y}_{R,\theta }(t') \right) +\P(\mathcal{Y}_{R,\theta}(t')^c).
\end{multline}
Now we claim
\begin{Lemma}\label{new1}
The following convergences hold:
\begin{align}
&\limsup_{t\to\infty} \E\Big[1-e^{-\theta \tilde{M}_t^\gamma(e^{-t'}{\rm BZ}_{R,t'}) }  ; \mathcal{Y}_{R,\theta}(t')\Big]\leq \epsilon,\label{eq:new4}\\
&\limsup_{t\to\infty} \E\Big[e^{-\theta \tilde{M}_t^\gamma(e^{-t'}\cup A_{\bf i}) } ; \mathcal{Y}_{R,\theta}(t') \Big]%\\&\leq \E\left( \exp(- C  \theta^{\frac{\sqrt{2\d}}{\gamma}}(1-c\epsilon)M'(e^{-l} \underset{{\bf i}}{\cup} A_{\bf i}) +2\epsilon C    )\right)
 \leq   \E\Big[\exp(- (C (\gamma) -\epsilon)\theta^{\frac{\sqrt{2\d}}{\gamma}} M'([0,1]^\d) + 2\epsilon    (C (\gamma) -\epsilon) )\Big]\label{eq:new1}
\end{align}
and   lower bound similar to \eqref{eq:new1} with a $\liminf_{t\to\infty}$ in the left-hand side.
\end{Lemma}

By taking the $\limsup_{t\to\infty}$ in \eqref{secineq} and by using Lemma \ref{new1}, we get
$$\limsup_{t\to\infty} \E\Big[e^{-\theta \tilde{M}_t^\gamma([0,1]^\d) } \Big]\leq \E\Big[\exp(- (C (\gamma) -\epsilon)\theta^{\frac{\sqrt{2\d}}{\gamma}} M'([0,1]^\d) + 2\epsilon    (C (\gamma) -\epsilon) )\Big]+\P(\mathcal{Y}_{R,\theta}(t')^c).$$
From Lemma \ref{lto0} we have  $\limsup_{t'\to\infty}\P(\mathcal{Y}_{R,\theta}(t')^c)\leq \epsilon$. We deduce
$$\limsup_{t\to\infty} \E\Big[e^{-\theta \tilde{M}_t^\gamma([0,1]^\d) } \Big]\leq \E\Big[\exp(- (C (\gamma) -\epsilon)\theta^{\frac{\sqrt{2\d}}{\gamma}} M'([0,1]^\d) + 2\epsilon    (C (\gamma) -\epsilon) )\Big]+\epsilon.$$
We can proceed in the same way for the lower bound. Since $\epsilon$ can be chosen arbitrarily close to $0$, the proof of Theorem \ref{main} follows, provided that we prove the  above lemma.\qed

To prove Lemma \ref{new1} we need the following proposition, which can actually be seen as the key tool of this subsection. Its proof requires some additional material and is carried out in Section \ref{sec:prop}.
\begin{Proposition}
\label{trlaplace}
Let $\gamma>\sqrt{2\d}$. There exists a constant $C (\gamma)>0$ such that for all $R\geq1$, $\theta>0$ and $ \epsilon>0$, we can find $t_0>0$ such that for all $t'>t_0$ satisfying $\frac{e^{t'}+1}{R+1}\in\N^*$, there exists $T>0$, such that
\begin{equation}
\label{eqtrlaplace} \left|\E\Big[ \exp\Big(-\theta \int_{[0,R]^\d}e^{ \gamma[Y_t(x)-a_t-\chi(x)]+\d t}dx \Big) -1\Big]+C (\gamma)  \mathtt{I}(\chi,\theta,\gamma)\right|\leq \epsilon \, \mathtt{I}(\chi,\theta,\gamma),
\end{equation}
for any  $t>T$ and for any function $\chi(\cdot) \in \mathcal{C}_R(t', \kappa_\d \ln t',\ln t)$, where 
\begin{equation}\label{def:ixi}
 \mathtt{I}(\chi,\theta,\gamma) =\theta^{\frac{\sqrt{2\d}}{\gamma}} \int_{[0,R]^{\d}}\big(\chi(x)-\frac{\ln\theta}{\gamma}\big)e^{-\sqrt{2\d}\chi(x)}\,dx.
 \end{equation}
\end{Proposition}

\vspace{2mm}
\noindent {\it Proof of Lemma \ref{new1}.} We first prove the first relation \eqref{eq:new4}. By the Markov property at time $t'$ and the scaling property \ref{scalingl} applied on the set ${\rm BZ}_{R,t'}$ we get that
\begin{align}
\nonumber \E\Big[1-&e^{-\theta \tilde{M}_t^\gamma(e^{-t'}{\rm BZ}_{R,t'}) }   ; \mathcal{Y}_{R,\theta}(t')\Big]\\
\nonumber &=  \E\Big[ \E\Big[1-\exp\big(-\theta \int_{e^{-t'}{\rm BZ}_{R,t'}}e^{\gamma[Y_{t',t-t'}(x)-a_t-\chi(x)]+\d t}dx   \big) \Big]_{\big| \chi(x)=-Y_{t'}(x)}; \mathcal{Y}_{R,\theta }(t')\Big]
\\
\label{expectation0} &= \E\Big[ \E\Big[1-\exp\big(-\theta \int_{ {\rm BZ}_{R,t'}}e^{\gamma[Y_{t-t'}(x)-a_t-\chi(x)]+\d (t-t')}dx   \big)\Big]_{\big| \chi(x)=-Y_{t'}(\frac{x}{e^{t'}})}; \mathcal{Y}_{R,\theta}(t')\Big].
\end{align}
We can find a finite collection of points in $ [0,e^{t'}]^\d $, call it $(y_j)_{j\in J}$,  such that

- for any distinct $j_1,..., j_{\d+2}\in J $, $ \underset{k=1}{\overset{\d+2}{\cap}} (y_{j_k}+[0,1]^\d)=\emptyset$,

- the closure $\overline{{\rm BZ}_{R,e^{t'}}}$ of ${\rm BZ}_{R,e^{t'}}$ is contained in $\underset{j\in J}{\cup}(y_j+[0,1]^\d)$.  

We do not detail the construction of these points but this is rather elementary: basically, you have to cover the red area in Figure \ref{croquis} with closed squares of side length $1$ (which corresponds to the width of the red strips). Of course, the squares that you choose  may overlap but if this covering is made efficiently enough, they will not overlap too much in such a way  that  any intersection of $\d+2$ such squares will be  empty. 

By using in turn  the elementary inequality  $1-\prod_{j\in J}u_j\leq \sum_{j\in J}1-u_j $ for $(u_j)_j\in [0,1]^J$ and then   invariance by translation, we get
\begin{align}
\E\Big[1-\exp(-\theta &\int_{ {\rm BZ}_{R,t'}}e^{\gamma[Y_{t-t'}(x)-a_t-\chi(x)]+\d (t-t')}dx )   \Big]_{\big| \chi(x)=-Y_{t'}(e^{-t'}x)}\nonumber
\\
&\leq \E\Big[1-\prod_{j\in J}\exp(-\theta \int_{ y_j+[0,1]^{\d}}e^{\gamma[Y_{t-t'}(x)-a_t-\chi(x)]+\d (t-t')}dx )   \Big]_{\big| \chi(x)=-Y_{t'}(e^{-t'}x)}\nonumber\\
&   \leq \sum_{j\in J}  \E\Big[1-\exp(-\theta \int_{[0,1]^\d}e^{\gamma[Y_{t-t'}(x)-a_t-\chi(x)]+\d (t-t')}dx ) \Big]_{\big| \chi(x)=-Y_{t'}(y_j+e^{-t'}x)}.\label{eq:new2}
\end{align}

Moreover, on $\mathcal{Y}_{R,\theta}(t')$, the function  $x\in [0,1]^{\d}\mapsto -Y_{t'}(y_j+e^{-t'}x) $ belongs to $ \mathcal{C}_1(t',\kappa_\d \ln t',\ln t)$ as soon as $\ln t>10\sqrt{2\d} t'$. So, by Proposition \ref{trlaplace}, we can find $t_0$  such that for any $t'> t_0$ satisfying $\frac{e^{t'}+1}{R+1}\in\N^*$ there exists $T>0$ such that for any $t>T$ and on $\mathcal{Y}_{R,\theta}(t')$  
\begin{align*}
 \E\Big[1-\exp(-\theta& \int_{[0,1]^\d}e^{\gamma[Y_{t-t'}(x)-a_t-\chi(x)]+\d (t-t')}dx ) \Big]_{\big| \chi(x)=-Y_{t'}(y_j+e^{-t'}x)}\\
\leq & (C (\gamma) +\epsilon)\int_{ y_j+[0,1]^\d} (-\frac{1}{\gamma}\ln \theta -Y_{t'}(xe^{-t'}))e^{\sqrt{2\d}(Y_{t'}(xe^{-t'}) +\frac{1}{\gamma}\ln \theta)}dx
\end{align*}
for any $j\in J$. Plugging this estimate into \eqref{eq:new2} yields
\begin{align}
\E\Big[1-\exp(-\theta &\int_{ {\rm BZ}_{R,t'}}e^{\gamma[Y_{t-t'}(x)-a_t-\chi(x)]+\d (t-t')}dx )   \Big]_{\big| \chi(x)=-Y_{t'}(e^{-t'}x)}\nonumber
\\
&   \leq \sum_{j\in J}  (C (\gamma) +\epsilon)\int_{ y_j+[0,1]^\d} (-\frac{1}{\gamma}\ln \theta -Y_{t'}(xe^{-t'}))e^{\sqrt{2\d}(Y_{t'}(xe^{-t'}) +\frac{1}{\gamma}\ln \theta)}dx.\label{eq:new3}
\end{align}
Now we may assume that $\kappa_\d\ln t'>\frac{1}{\gamma}\ln \theta$ in such a way that, on $\mathcal{Y}_{R,\theta}(t')$  , we have $(-\frac{1}{\gamma}\ln \theta -Y_{t'}(xe^{-t'}))\geq 0$ for $x\in [0,e^{t'}]^\d$. Furthermore, the relation $ \underset{k=1}{\overset{\d+2}{\cap}} (y_{j_k}+[0,1]^\d)=\emptyset$ (valid for all family of distinct indices) entails that $\sum_{j\in J} \1_{\{y_{j_k}+[0,1]^\d\}}\leq (d+2)\1_{\{\overline{{\rm BZ}_{R,t'}}\}} $. Hence, on $\mathcal{Y}_{R,\theta}(t')$
\begin{align}
\E\Big[1-\exp(-\theta &\int_{ {\rm BZ}_{R,t'}}e^{\gamma[Y_{t-t'}(x)-a_t-\chi(x)]+\d (t-t')}dx )   \Big]_{\big| \chi(x)=-Y_{t'}(e^{-t'}x)}\nonumber
\\
&   \leq   (\d+2)( {C (\gamma)}+\epsilon) \int_{ {\rm BZ}_{R,t'}} (-\frac{1}{\gamma}\ln \theta -Y_l(xe^{-t'}))e^{\sqrt{2\d}(Y_{t'}(xe^{-t'}) +\frac{1}{\gamma}\ln \theta)}dx 
\nonumber\\
\label{3.7} &\leq (\d+2)( {C (\gamma)}+1) \int_{e^{-t'} {\rm BZ}_{R,t'}} (-\frac{1}{\gamma}\ln \theta -Y_{t'}(x))e^{\sqrt{2\d}(Y_{t'}(x) +\frac{1}{\gamma}\ln \theta)+\d t'}dx  .
\end{align}
The last inequality results from the change of variables $xe^{-t'}\to x$. We recognize the expressions of the  martingales $M_{t'}^{\sqrt{2\d}}$ and $M'_{t'}$ as in \eqref{deuxdef}. By gathering \eqref{expectation0} and the above relation, we deduce
\begin{align}
\nonumber \E\Big[1-&e^{-\theta \tilde{M}_t^\gamma(e^{-t'}{\rm BZ}_{R,t'}) }   ; \mathcal{Y}_{R,\theta}(t')\Big] \\&\leq (\d+2)( {C (\gamma)}+1)\E\Big[\theta^{\frac{\sqrt{2\d}}{\gamma}}\big(-\gamma^{-1}\ln \theta M^{\sqrt{2\d}}_{t'}(e^{-t'}{\rm BZ}_{R,t'})+ M'_{t'}(e^{-t'}{\rm BZ}_{R,t'})\big) ; \mathcal{Y}_{R,\theta}(t')\Big].
\end{align}
By using the definition \eqref{defYRz} of $\mathcal{Y}_{R,\theta}(t')$, we see that this latter quantity is less than $\epsilon  (\d+2)( {C (\gamma)}+1)$. By choosing $\epsilon$ as small as we please, we complete the proof of the first relation  \eqref{eq:new4}.

\vspace{2mm}
Now we prove \eqref{eq:new1}. As previously, we first apply the Markov property at time $t'$ and the scaling property \eqref{scalingl}. 
\begin{align}
\nonumber \E\Big[1-&e^{-\theta \tilde{M}_t^\gamma(e^{-t'} \cup A_{\bf i}) }   ; \mathcal{Y}_{R,\theta}(t')\Big]\\
\label{eq:new5} &= \E\Big[ \E\Big[1-\exp\big(-\theta \int_{ \cup A_{\bf i}}e^{\gamma[Y_{t-t'}(x)-a_t-\chi(x)]+\d (t-t')}dx   \big)\Big]_{\big| \chi(x)=-Y_{t'}( x e^{-t'})}; \mathcal{Y}_{R,\theta}(t')\Big].
\end{align}
The important point here is to see that for any $t\geq 0$, the process $(Y_t(x))_{x\in\r^d}$ is decorrelated at distance $1$ (recall that $k$ has compact support in the ball $B(0,1)$). Therefore, the random variables $\big( \int_{  A_{\bf i}}e^{\gamma[Y_{t-t'}(x)-a_t-\chi(x)]+\d (t-t')}dx\big)_{\bf i}$ appearing in the latter expectation are independent since ${\rm dist}(A_{\bf i},A_{\bf j})\geq 1$ for any ${\bf i } \neq {\bf j}$. We deduce that
\begin{align}
\E\Big[1-e^{-\theta \tilde{M}_t^\gamma(e^{-t'}\cup A_{\bf i}) } &; \mathcal{Y}_{R,\theta}(t') \Big]\nonumber\\
=&\E\Big[ \underset{{\bf i}\in \{1,...,m\}^\d}{\prod} \E\Big[ \exp(-\theta \int_{A_{\bf i}}e^{\gamma[Y_{t-t'}(x)-a_t-\chi(x)]+\d (t-t')}dx  ) \Big]_{\big | \chi(x)=-Y_{t'}(xe^{-t'})};\mathcal{Y}_{R,\theta}(t') \Big].\label{eq:new6}
\end{align}
As previously, we can choose $t$ sufficiently large so that, on $\mathcal{Y}_{R,\theta}(t')$ and for any $j\in J$,  the function   $x\in [0,R]^\d\mapsto -Y_{t'}(e^{-t'}(x+a_{\bf i}))$ belongs to $ \mathcal{C}_R(t',\kappa_\d\ln t',\ln t) $. We can then apply Proposition \ref{trlaplace} once again and get some $t_0>0$ such that for all $t'>t_0$ (with $\frac{e^{t'}+1}{R+1}\in\N^*$) there exists $T>0$ such that for all $t\geq T$ and all ${\bf i}$,
\begin{eqnarray*}
\Big| \E\Big[ \exp(-\theta \int_{A_{\bf i}}e^{\gamma[Y_{t-t'}(x)-a_t-\chi(x)]+\d (t-t')}dx  \Big]-1+   C  (\gamma) \mathtt{I}_\d (\chi ( \cdot) -\frac{1}{\gamma}\ln \theta) \Big| \leq \epsilon \mathtt{I}_\d (\chi ( \cdot)-\frac{1}{\gamma}\ln \theta),
\end{eqnarray*}
 with $\chi (x)=-Y_{t'}(e^{-t'}x)$. By plugging this estimate into \eqref{eq:new6} and by making a change of variables $xe^{-t'}\to x$, we obtain (once again by identifying $M^{\sqrt{2\d}}_{t'}$ and $M_{t'}'$)
\begin{align*}
\E\Big[e^{-\theta \tilde{M}_t^\gamma(\cup A_{\bf i}) } ; \mathcal{Y}_{R,\theta}(t') \Big]\leq \E\Big[ \prod_{{\bf i}\in \{1,...,m\}^\d} \Big( 1-(C (\gamma)-\epsilon)\theta^{\frac{\sqrt{2\d}}{\gamma}}\big[-\gamma^{-1}\ln \theta M^{\sqrt{2\d}}_{t'}(e^{-t'}A_{\bf i})+M_{t'}'(e^{-t'} A_{\bf i})\big]   \Big)  ;\mathcal{Y}_{R,\theta}(t') \Big].
\end{align*} 
 On $\mathcal{Y}_{R,\theta}(t')$, $\forall{\bf i} \in \{1,...,m\}^\d$ we  have $ |\gamma^{-1}\ln \theta M_{t'}^{\sqrt{2\d}}(e^{-t'}A_{\bf i})|+|M_{t'}'(e^{-t'}A_{\bf i})| \leq c\frac{\ln t'}{(t')^a}\leq \epsilon$ for $t'$ large enough.  Indeed, this is clear for  $M_{t'}^{\sqrt{2\d}}$ and, for $M_{t'}'(e^{-t'}A_{\bf i})$, it suffices to observe that:
$$M_{t'}'(e^{-t'}A_{\bf i})=\int_{e^{-t'}A_{\bf i}}-Y_{t'}(x)e^{\sqrt{2\d}Y_{t'}(x)+\d t'}\,dx\leq R^\d\sup_{u\geq \kappa_\d\ln t'} u e ^{-\sqrt{2\d}u}.$$
Then, by using the inequality $\prod_{i\in I}(1-u_i)\leq e^{-\sum_{i\in I}u_i}$ for $u_i\in[0,1]$, we obtain
\begin{eqnarray*}
\E\Big[e^{-\theta \tilde{M}_t^\gamma(\cup A_{\bf i}) } ; \mathcal{Y}_{R,\theta}(t') \Big]\leq  \E\Big[\exp\Big(-(C (\gamma)-\epsilon)\theta^{\frac{\sqrt{2\d}}{\gamma}} [M'_{t'}(e^{-t'} \underset{{\bf i}}{\cup} A_{\bf i}) -\gamma^{-1}\ln \theta M^{\sqrt{2\d}}_{t'}(e^{-t'}\underset{\bf i}{\cup} A_{\bf i})]   \Big) ;\mathcal{Y}_{R,\theta}(t')\Big].
\end{eqnarray*}
%\begin{eqnarray*}
%\E\left(e^{-\theta \tilde{M}_t^\gamma(\cup A_{\bf i}) } ; \mathcal{Y}_{R,\theta}(l) \right)\leq \E\left(\exp(-(C (\gamma)-\epsilon)\theta^{\frac{\sqrt{2\d}}{\gamma}} [M'(e^{-l} \underset{{\bf i}}{\cup} A_{\bf i}) -\gamma^{-1}\ln \theta M^{\sqrt{2\d}}_l(e^{-l}\underset{\bf i}{\cup} A_{\bf i})]   ) ;\mathcal{Y}_{R,\theta}(l)\right)
%\end{eqnarray*} 
Recall that, onn $\mathcal{Y}_{R,\theta}(t')$,
\begin{equation*}
|   M'_{t'}(e^{-t'} \underset{{\bf i}}{\cup} A_{\bf i})-M_{t'}'([0,1]^\d)| +|\gamma^{-1}\ln \theta M^{\sqrt{2\d}}_{t'}(e^{-t'}\underset{\bf i}{\cup} A_{\bf i}) |\leq |\gamma^{-1}\ln \theta M^{\sqrt{2\d}}_{t'}([0,1]^\d) | + | M_{t'}'(e^{-t'}{\rm BZ}_{R,t'})| \leq \epsilon \theta^{-\frac{\sqrt{2\d}}{\gamma}}.
\end{equation*}
so $\E\left(e^{-\theta \tilde{M}_t^\gamma(\cup A_{\bf i}) } ; \mathcal{Y}_{R,\theta}(t') \right)\leq     \E\left( \exp(- (C (\gamma) -\epsilon)\theta^{\frac{\sqrt{2\d}}{\gamma}} M'([0,1]^\d) + 2\epsilon    (C (\gamma) -\epsilon) )\right)$. 
The lower bound of (\ref{eqLapl}) can be derived in the same way.\qed

%%%%%%
\subsection{Proof of Corollary \ref{coro:meas}}
%%%%%%%%%%%%%%%%%%%%%%%%%%%%
Here we assume that Theorem \ref{main} holds and we show that this implies   convergence in law in the sense of weak convergence of measures. For  $a>0$, let us denote by $C_a$ the cube $[-a,a]^\d$.   Since for all bounded continuous function $f$ compactly supported in $C_R$, we have
$$0\leq \int_{C_a} f(x)\tilde{M}_t^\gamma(dx)\leq \|f\|_\infty  \tilde{M}_t^\gamma(C_a)$$
and since the right-hand side is tight, this ensures that the family of random measures $\big(\tilde{M}_t^\gamma(dx)\big)_t$ is tight for the weak convergence of measures on $C_a$. Since we can find  a sequence $(f_n)_n$ of smooth strictly positive functions on $C_a$ that is dense in the set of nonnegative continuous compactly supported functions in $C_a$ for the uniform topology, uniqueness in law then results from Theorem \ref{main}. As it is   rather a standard argument of functional analysis, we let the reader check the details, if need be.\qed

%%%%%%%%%%%%%%%%%%%%%%%%%%%%%%%%%%%%%%%%%%%%%%%%%%%%%%%%
%%%%%%%%%%%%%%%%%%%%%%%%%%%%%%%%%%%%%%%%%%%%%%%%%%%%%%%%
\section{Estimation on the tail of distribution of $\tilde{M}_t^\gamma([0,1]^\d)$}
%%%%%%%%%%%%%%%%%%%%%%%%%%%%%%%%%%%%%%%%%%%%%%%%%%%%%%%%
%%%%%%%%%%%%%%%%%%%%%%%%%%%%%%%%%%%%%%%%%%%%%%%%%%%%%%%%
In this section, we will identify the path configuration $t\mapsto Y_t(x)$ that really contribute to the behaviour of the measure $\tilde{M}_t^\gamma$. We will show that, for these paths, $Y_t(x)$ typically goes faster than $a_t=-\frac{3}{2\sqrt{2\d}}\ln t$.

To quantify the above rough claim, we will establish %(recall the expression \eqref{def:ixi})
\begin{Proposition}
\label{partieneg}
Let $R,\,  \epsilon>0$. There exists a constant $A>0$ such that for any $t',T$ large enough we have   
\begin{equation}
\label{eqpartieneg}
\E\Big[1- \exp\Big(- \int_{[0,R]^\d}e^{ \gamma[Y_t(x)-a_t-\rho(x)]+\d t}\1_{\{Y_t(x)\leq a_t+\chi(x)-A \}}dx \Big) \Big]\leq \epsilon\int_{[0,R]^\d} \chi(x)  e^{-\sqrt{2\d}\chi(x)}dx,%\mathtt{I}(\chi,1,1) 
\end{equation}
for any $t \geq T$ and $\chi(\cdot) \in \mathcal{C}_R(t', \kappa_\d \ln t',\ln t)$.
\end{Proposition}
 
Then we  focus on the    shape of the tail distribution of $\tilde{M}_t$. For instance, it is well known in Tauberian theory that an estimate of the type
\begin{equation}\label{tauber}
C^{-1}xe^{-\sqrt{2\d} x}\leq 1- \E[e^{-e^{-\gamma x}\tilde{M}_t([0,R]^\d)}]\leq Cxe^{-\sqrt{2\d} x}
\end{equation}
valid for $x>0$ gives you a tail estimate for $\tilde{M}_t([0,R]^\d)$ of the type
$$\P\big(\tilde{M}_t([0,R]^\d)>e^{\gamma x}\big)\asymp  xe^{-\sqrt{2\d} x}$$ as $x\to \infty$. Basically, the following proposition is a functional version of \eqref{tauber}, meaning that we will replace the variable $x$ by some function $\chi$. Thus we claim

\begin{Proposition}
\label{uplowbound}
There exist $c_1,\, c_2,$ such that for any $ t'>0$ there exists $T>0$  such that for any $R\in [1,\ln t']$:
\begin{itemize}
\item for any $t\geq T$ and  any  $\chi(\cdot) \in \mathcal{C}_R(t', \kappa_\d \ln t',\ln t)$,
\begin{equation*}
 \E\Big[1- \exp\Big(- \int_{[0,R]^\d}e^{ \gamma[Y_t(x)-a_t-\chi(x)]+\d t}dx \Big) \Big]\leq c_2  \int_{[0,R]^\d} \chi(x)  e^{-\sqrt{2\d}\chi(x)}dx .\label{equpbound}
\end{equation*}
\item for any $\chi(\cdot) \in \mathcal{C}_R(t', \kappa_\d \ln t',+\infty)$
\begin{equation*}
 c_1  \int_{[0,R]^\d} \chi(x)  e^{-\sqrt{2\d}\chi(x)}dx\leq \liminf_{t\to\infty} \E\Big[1- \exp\Big(- \int_{[0,R]^\d}e^{ \gamma[Y_t(x)-a_t-\chi(x)]+\d t}dx \Big) \Big] .\label{eqlowbound}
\end{equation*}
\end{itemize}
\end{Proposition}

\subsection{Proof of Proposition \ref{partieneg}}
%%%%%%%%%%%%%%%%%%%%%%%%%%%%%%%%%%%%%%%%%%

Fix $\epsilon>0$. We consider $t'>0$ et $R\geq 1$ such that $\frac{e^{t'}+1}{R+1}\in\N^*$. We have for $t>e^{t'}$
\begin{align*}
\nonumber    \E\Big[1- \exp\Big(- &\int_{[0,R]^\d}e^{ \gamma[Y_t(x)-a_t-\chi(x)]+\d t}\1_{\{Y_t(x)\leq a_t+\chi(x)-A \}}dx \Big) \Big]
 \\
   \leq &\E\Big[1- \exp\Big(- \int_{[0,R]^\d}e^{ \gamma[Y_t(x)-a_t-\chi(x)]+\d t}\1_{\{ \sup_{s\in[\ln t',t]}Y_s(x)\leq \chi(x),\,  Y_t(x)\leq a_t+\chi(x)-A \}}dx \Big)  \Big]\nonumber\\
&+  \P\Big( \sup_{ x\in [0,R]^\d} \sup_{s\in [\ln t',\infty[}Y_s(x)\geq \chi(x)\Big).
\end{align*}
If $\chi(\cdot) \in \mathcal{C}_R(t', \kappa_\d \ln t',\ln t)$ (with $t'$ large enough so as to make $\kappa_\d\ln t'>10$), we can estimate the probability in the right-hand side with the help of Lemma \ref{maxfx}. If $t'$ is again large enough, we have
$$(\ln t')^\frac{3}{8}+  \chi(x)^{\frac{3}{4}}\leq \frac{\epsilon}{2} \kappa_\d\ln t'+\chi(x)^{\frac{3}{4}}\leq  \epsilon\chi(x)$$ in such a way that 
\begin{align} 
&\E\Big[1- \exp\Big(- \int_{[0,R]^\d}e^{ \gamma[Y_t(x)-a_t-\chi(x)]+\d t}\1_{\{Y_t(x)\leq a_t+\chi(x)-A \}}dx \Big) \Big]\label{frisette} \\
&\leq   \E\Big[1- \exp\Big(- \int_{[0,R]^\d}e^{ \gamma[Y_t(x)-a_t-\chi(x)]+\d t}\1_{\{ \sup_{s\in[\ln t',t]}Y_s(x)\leq \chi(x),\,  Y_t(x)\leq a_t+\chi(x)-A \}}dx \Big)  \Big]\nonumber\\
&\quad+\epsilon  \int_{[0,R]^\d} \chi(x)  e^{-\sqrt{2\d}\chi(x)}dx .\nonumber
\end{align}
So we need to bound the first term in the right-hand side of \eqref{frisette}.  To this purpose, we will use 
Lemma \ref{LemL2}. We consider the constants $c_4,c_5$ of this lemma and we decompose the event
$\{\sup_{s\in [\ln t',t]}Y_s(x)\leq \chi(x),\,   Y_t(x)\leq a_t+\chi(x)-A \}$ for some constant $A>0$ as follows
\begin{eqnarray*}
\1_{\{\sup_{s\in [\ln t',t]}Y_s(x)\leq \chi(x),\,   Y_t(x)\leq a_t+\chi(x)-A \}} \leq \1_{E^1_{t',t}(x)}+\1_{E^2_{t',t}(x)}+\1_{E^3_{t',t}(x)},
\end{eqnarray*}
where the set $E^1_{t',t}(x)$, $E^2_{t',t}(x)$, $E^3_{t',t}(x)$ as follows.  For any $ j\geq 1$, we define $a_{j}:= e^{\frac{c_5}{2}j}$ ($c_5$ is defined by Lemma \ref{LemL2}). Then we set
\begin{align*}
E^1_{t',t}(x):=&  \{\sup_{s\in [\ln t',t]}Y_s(x)\leq \chi(x),\, \sup_{s\in [\frac{t}{2},t]}Y_s(x)\leq a_t+\chi(x)+L,\,  Y_t(x)\leq a_t+\chi(x)-A \}, 
\\
E^2_{t',t}(x):=&  \bigcup_{j\geq L+ 1}   \{\sup_{s\in[\ln t',t]}Y_s(x)\leq \chi(x),\,  \sup_{s\in[\frac{t}{2},t-a_{j}]}Y_s(x)-a_t-\chi(x)\in [j,j+1],
\\&\quad \quad\quad  \,\sup_{s\in[t-a_{j},t]}Y_s(x)\leq a_t+\chi(x)+j  ,\,  Y_t(x)\leq a_t+\chi(x)-A \},  
\\
E^3_{t',t}(x):= &\bigcup_{j\geq L+ 1}   \{\sup_{s\in[\ln t',t]}Y_s(x)\leq \chi(x),\,\sup_{s\in[\frac{t}{2},t-a_{j}]}Y_s(x)\leq a_t+\chi(x)+j-1,\,   \sup_{s\in[t-a_{j},t]}(x)-a_t-\chi(x)\in[j,j+1]) \}.
\end{align*}
According to Lemma \ref{LemL2}, there exists $T>0$ such that for all $t>T$ and $\chi(\cdot)\in \mathcal{C}_R(t', 10, +\infty)$
\begin{align}
\nonumber
\P\Big(\sup_{x\in[0,R]^\d} \1_{E^3_{t',t}(x)}=1\Big)  \leq &\underset{j\geq L+ 1}{\sum}c_4(1+a_j) e^{-c_5 j}\int_{[0,R]^\d} (\sqrt{\ln t'}+\chi(x)) e^{-\sqrt{2\d}\chi(x)}dx\nonumber\\
\leq &c \,e^{-\frac{c_5}{2} L}\int_{[0,R]^\d} (\sqrt{\ln t'}+\chi(x)) e^{-\sqrt{2\d}\chi(x)}dx .\label{eqLemL2Zguin} 
\end{align}
If we further impose $\chi(\cdot) \in \mathcal{C}_R(t', \kappa_\d \ln t',+\infty)$ while choosing $t'$ large enough so as to make the term $\sqrt{\ln t'}$ smaller than $\epsilon \kappa_\d \ln t'$ (and therefore less than $\epsilon \chi$) as well as choosing $L$ large enough to have $c \,e^{-\frac{c_5}{2} L}\leq \epsilon$, we deduce  
 \begin{equation}
\label{eq:rem1}   \E\Big[1- \exp\Big(- \int_{[0,R]^\d}e^{ \gamma[Y_t(x)-a_t-\chi(x)]+\d t} \1_{E^3_{t',t}(x)}dx \Big) \Big]\leq  2\epsilon\int_{[0,R]^\d} \chi(x)  e^{-\sqrt{2\d}\chi(x)}dx.
\end{equation}

Now we focus on $E^1_{t,t'}(x)$. By partitioning the event $\{Y_t(x)\leq a_t+\chi(x)-A\}$ as
 $$\{Y_t(x)\leq a_t+\chi(x)-A\}=\bigcup_{p\geq 0}\{Y_t(x)- a_t-\chi(x)+A\in [-p-1,-p]\}$$
and by using the relation $1-e^{-u}\leq u$ for $u\geq 0$, we obtain
\begin{align}
  \E\Big[1- &\exp\Big(- \int_{[0,R]^\d}e^{ \gamma[Y_t(x)-a_t-\chi(x)]+\d t} \1_{E^1_{t',t}(x)}dx \Big)\Big] \leq   \E\Big[  \int_{[0,R]^\d}e^{ \gamma[Y_t(x)-a_t-\chi(x)]+\d t}\1_{E^1_{t,t'}(x)}dx \Big]
\nonumber\\
 \leq& e^{-\gamma A} \sum_{p\geq 0} e^{-\gamma p}e^{\d t} \int_{[0,R]^\d}\P\Big(  E^1_{t,t'}(x), Y_t(x)- a_t-\chi(x)+A\in [-p-1,-p]\Big) dx.  \label{combine1}
\end{align}
By the Girsanov's transform  (with density $e^{\sqrt{2\d}Y_t(x)+\d t}$), we obtain for any $x\in [0,R]^\d$ and $p\geq 0$,
\begin{align}
&\P\Big(  E^1_{t,t'}(x), Y_t(x)- a_t-\chi(x)+A\in [-p-1,-p]\Big) \nonumber
\\
&\leq e^{-\sqrt{2\d}[a_t+\chi(x)-A-p-1 ] -\d t}\P_{-\chi(x)}\Big(\sup_{s\in [\ln t',t]}B_s\leq 0,\, \sup_{s\in [\frac{t}{2},t]}B_s\leq a_t+L,\,   B_t-a_t-A\in [-p-1,-p] \Big),\label{combine2}
\end{align}
where, under $\P_{-\chi(x)}$, the process $B$ is a standard Brownian motion starting from $-\chi(x)$. At this step, we observe that similar quantities have been treated in \cite{Mad13}. More precisely, a combination of (B.5) and (B.6) in \cite{Mad13} shows that, for some constant $\bar{c}>0$ (which does not depend on relevant quantities) 
\begin{align}\label{combine3}
\P_{-\chi(x)}\Big(\sup_{s\in [\ln t',t]} B_s\leq 0,\, \sup_{s\in [\frac{t}{2},t]}B_s\leq a_t+L&,\,   B_t-a_t-A\in [-p-1,-p] \Big)\\
\leq &t^{-3/2}\bar{c}(L+A+p)\E\Big[(\chi(x)-B_{\ln t'})\1_{\{ \chi(x)-B_{\ln t'}\geq 0 \}}\Big].\nonumber
\end{align}
Finally by combining \eqref{combine1}+\eqref{combine2})+\eqref{combine3}) we get:
\begin{align}
  \E\Big[1- \exp\Big(-& \int_{[0,R]^\d}e^{ \gamma[Y_t(x)-a_t-\chi(x)]+\d t} \1_{E^1_{t',t}(x)}dx \Big)\Big]\nonumber\\
\leq &e^{-(\gamma-\sqrt{2\d})A}\sum_{p\geq 0}(L+A+p) e^{-(\gamma-\sqrt{2\d})p}\int_{[0,R]^\d}\chi(x)e^{-\sqrt{2\d}\chi(x)} dx\nonumber\\
\leq & (L+A)e^{-(\gamma-\sqrt{2\d})A}c \int_{[0,R]^\d}\chi(x)e^{-\sqrt{2\d}\chi(x)} dx.\label{eq:rem2}
\end{align}
where we took for instance  $c=2e^{\sqrt{2\d}\bar{c}}  \sum_{p\geq 0} (1+p) e^{-(\gamma-\sqrt{2\d})p} $.

Finally we treat the contribution of the term $E^2_{t,t'}(x)$. First, we can follow the same argument as for $E^1_{t,t'}(x)$ to get
\begin{align}
  \E\Big[1- \exp\Big(-& \int_{[0,R]^\d}e^{ \gamma[Y_t(x)-a_t-\chi(x)]+\d t} \1_{E^2_{t',t}(x)}dx \Big)\Big]\nonumber\\
\leq &\sum_{j\geq L+1}   e^{-\gamma A} \sum_{p\geq 0} e^{-\gamma p}e^{\d t} \int_{[0,R]^\d}\P\Big( E^2_{t',t}(x), Y_t(x)- a_t-\chi(x)+A\in [-p-1,-p]\Big)\,dx.\label{eq:om1}
\end{align}
By the Girsanov's transform again (with density $e^{\sqrt{2\d}Y_t(x)+\d t}$), we can estimate the probability in \eqref{eq:om1} by
\begin{align}
\label{eq:om2}  \P\Big( &E^2_{t',t}(x),  Y_t(x)- a_t-\chi(x)+A\in [-p-1,-p]\Big)\\
 & \leq   e^{-\sqrt{2\d}[a_t+\chi(x)-A-p-1 ] -\d t}\P_{-\chi(x)}\Big(\sup_{s\in [\ln t',t]}B_s\leq 0, \sup_{s\in [\frac{t}{2},t-a_{j}]}B_s   -a_t\in [j,j+1],\nonumber\\
 &\hspace{5.5cm}\sup_{s\in[t-a_{j},t]}B_s\leq a_t+j ,\,  B_t-a_t+A\in [-p-1,-p]  \Big).\nonumber
\end{align}
Once again, we use (B.3) in \cite{Mad13} to see that this latter quantity is smaller than
\begin{equation}
\label{eq:om3} ce^{-\d t} e^{\sqrt{2\d} (A+p+1)}e^{-\sqrt{2\d}\chi(x)}(1+j+A+p) a_j^{-\frac{1}{2}}\chi(x).
\end{equation}
By recalling that $a_j=e^{\frac{c_5}{2}j}$ and by combining \eqref{eq:om1}+\eqref{eq:om2}+\eqref{eq:om3}, we get:
\begin{align}
  \E\Big[1- \exp\Big(-& \int_{[0,R]^\d}e^{ \gamma[Y_t(x)-a_t-\chi(x)]+\d t} \1_{E^2_{t',t}(x)}dx \Big)\Big]\nonumber \\
  \leq & ce^{-(\gamma-\sqrt{2\d}) A} \sum_{j\geq L+1} e^{-\frac{c_5}{4}j}   \sum_{p\geq 0} (1+j+A+p)  e^{-(\gamma-\sqrt{2\d}) p} \int_{[0,R]^\d}    e^{-\sqrt{2\d}\chi(x)} \chi(x) dx 
\nonumber\\
\label{eq:rem3}\leq&  ce^{-\frac{c_5}{8}L}  Ae^{-(\gamma-\sqrt{2\d})A} \int_{[0,R]^\d}    e^{-\sqrt{2\d}\chi(x)} \chi(x) dx   .
\end{align}
Now recall that our purpose is to estimate the right-hand side in  \eqref{frisette}. The expectation in this right-hand is estimated by combining  \eqref{eq:rem1}+\eqref{eq:rem2}+\eqref{eq:rem3} in such a way that
\begin{align}
  \E\Big[1- \exp\Big(- &\int_{[0,R]^\d}e^{ \gamma[Y_t(x)-a_t-\chi(x)]+\d t}\1_{\{Y_t(x)\leq a_t+\chi(x)-A \}}dx \Big)\Big]\nonumber\\ &\leq c (e^{-[\gamma-\sqrt{2\d}]A}[ (L+A)+ e^{-\frac{c_5}{8}L} A ]+ 2\epsilon )\int_{[0,R]^\d}    e^{-\sqrt{2\d}\chi(x)} \chi(x) dx    .\label{fauve2} 
\end{align}
So it suffices to choose $A$ large enough such that $ ce^{-[\gamma-\sqrt{2\d}]A}[ (L+A)+ e^{-\frac{c_5}{8}L} A ]\leq \epsilon $ to conclude  the proof of Proposition \ref{partieneg}.\qed

\subsection{Proof of Proposition \ref{uplowbound}}
%%%%%%%%%%%%%%%%%%%%%%%%%%%%%%%%%%%%%%%%%
%\subsubsection{Trash}
%\begin{equation*} 
% \Squ_{t }^{\rho(x),0}= \lbrace f: \: \overline{f}_{t/2} \leq \rho(x) , \:  \overline{f}_{[t/2,t]} \leq a_t + \rho(x), \: f_t \in I_t(\rho(x))   \rbrace.
%\end{equation*}
%\subsubsection{go on}
%%%%%%%%%%%%%%%%%
The first relation of Proposition \ref{uplowbound} is an easy consequence of Lemma \ref{tightmad} and (\ref{fauve2}). Indeed by using the relation $1-e^{-(u+v)}\leq (1-e^{-u})+(1-e^{-v}) $ for $u,v\geq 0$ and by applying (\ref{fauve2}) with $\epsilon=1$, we obtain
\begin{align}
 \E\Big[1- & \exp\Big(-  \int_{[0,R]^\d}e^{ \gamma[Y_t(x)-a_t-\chi(x)]+\d t}dx\Big)  \Big]\nonumber\\
 \leq &\E\Big[1- \exp\Big(- \int_{[0,R]^\d}e^{ \gamma[Y_t(x)-a_t-\chi(x)]+\d t}\1_{\{ Y_t(x)\geq a_t+\chi(x)-1\}}\Big)  \Big]\nonumber\\
 &  +\E\Big[1- \exp\Big(- \int_{[0,R]^\d}e^{ \gamma[Y_t(x)-a_t-\chi(x)]+\d t}\1_{\{ Y_t(x)\leq a_t+\chi(x)-1\}}  dx\Big)\Big]\nonumber
\\
   \leq &\P\Big(\exists x\in [0,R]^\d,\, Y_t(x)\geq a_t+\chi(x)-1\Big) +  c (e^{-[\gamma-\sqrt{2\d}]A}[ (L+A)+ e^{-\frac{c_5}{8}L} A ]+ 2  )\int_{[0,R]^\d} \chi(x)   e^{-\sqrt{2\d}\chi(x)}  dx \nonumber  \\
 \leq & c' \int_{[0,R]^\d}  \chi(x)  e^{-\sqrt{2\d}\chi(x)}  dx ,
\end{align}
with $c':= c_2  + c( (L+1)+ e^{-\frac{c_5}{8}L}   + 2 )$ where $c_2$ is the constant appearing in Lemma \ref{tightmad}.

\medskip
Now we  prove the second inequality.  Let us define the event
$$E(R,t,\chi)=\{\exists x_0\in[0,R]^\d,   Y_t(x_0)\geq a_t + \chi(x_0)  \}. $$
From Lemma \ref{tightmad}, there exists   $c_2>0$ such that for any $t'\geq 2$, there exists $T>0$ such that for any $R\in [1,\ln t']$ and $t\geq T$ 
\begin{equation}\label{eq:psg1}
\P\Big(E(R,t,\chi)  \Big)=\P\Big(\exists x_0\in[0,R]^\d,   Y_t(x_0)\geq a_t + \chi(x_0)  \Big)\geq c_2\int_{[0,R]^\d}   \chi(x)  e^{-\sqrt{2\d}\chi(x)} dx
\end{equation}
for any function $\chi \in \mathcal{C}_R(t', \kappa_\d\ln t',\ln t)$. Then we  observe that
\begin{align*}
 \E\Big[1- \exp\Big(- &\int_{[0,R]^\d}e^{ \gamma[Y_t(x)-a_t-\chi(x)]+\d t}dx\Big) \Big]\\
 \geq  &\E\Big[\Big(1- \exp\big(- \int_{[0,R]^\d}e^{ \gamma[Y_t(x)-a_t-\chi(x)]+\d t}dx\big)\Big) \1_{E(R,t,\chi)}\Big]\\
\geq  &\E\Big[\Big(1- \exp\big(- \int_{B(x_0,e^{-t})}e^{ \gamma[Y_t(x)-a_t-\chi(x)]+\d t}dx\big)\Big) \1_{E(R,t,\chi)}\Big].
\end{align*}
Now we use the Jensen inequality to get for some fixed constant $c$ (which is precisely the Lebesgue volume of the unit ball in such a way that $|B(x,e^{-t})|=ce^{\d t}$)
\begin{align*}
 \E\Big[1- \exp\Big(- &\int_{[0,R]^\d}e^{ \gamma[Y_t(x)-a_t-\chi(x)]+\d t}dx\Big) \Big]\\
 \geq   &\E\Big[\Big(1- \exp\big(-  ce^{ c^{-1}e^{\d t}\int_{B(x_0,e^{-t})} \gamma[Y_t(x)-a_t-\chi(x)]dx}\big)\Big) \1_{E(R,t,\chi)}\Big]\\
  \geq   &\E\Big[\Big(1- \exp\big(-  ce^{ \inf_{x_0\in[0,R]^\d}H_t(x_0)}\big)\Big) \1_{E(R,t,\chi)}\Big],
\end{align*}
where we have set 
$$H_t(x_0)=c^{-1}e^{\d t}\int_{B(x_0,e^{-t})} \gamma[X_t(x)-X_t(x_0)-\chi(x)+\chi(x_0)]dx.$$
Let us consider $M>0$. We have 
\begin{align}
 \E\Big[1- \exp\Big(- &\int_{[0,R]^\d}e^{ \gamma[Y_t(x)-a_t-\chi(x)]+\d t}dx\Big) \Big]\nonumber\\
  \geq   &\E\Big[\Big(1- \exp\big(-  ce^{-M}\big)\Big) \1_{E(R,t,\chi)}\1_{\{ \inf_{x_0\in[0,R]^\d}H_t(x_0)>-M\}}\Big]\nonumber\\
\geq & \Big(1- \exp\big(-  ce^{-M}\big)\Big) \Big( \P\big(E(R,t,\chi)\big)-\P\big( \inf_{x_0\in[0,R]^\d}H_t(x_0)\leq -M\big)\Big)\nonumber\\
\geq & \Big(1- \exp\big(-  ce^{-M}\big)\Big) \Big(c_2\int_{[0,R]^\d}   \chi(x)  e^{-\sqrt{2\d}\chi(x)} dx-\P\big( \inf_{x_0\in[0,R]^\d}H_t(x_0)\leq -M\big)\Big).\label{eq:psg2}
\end{align}
In the last equality, we have used   \eqref{eq:psg1}. Now we claim
\begin{Lemma}\label{kolm}
For each fixed $R$, we have
$$\liminf_{t\to\infty}\P\big( \inf_{x_0\in[0,R]^\d}H_t(x_0)\leq -M\big)=0.$$
\end{Lemma}
Taking the $\liminf_{t\to\infty}$ in \eqref{eq:psg2} then completes the proof of Proposition \ref{uplowbound}.\qed

\medskip
\noindent {\it Proof of Lemma \ref{kolm}.} Since $\chi(\cdot)\in \mathcal{C}_R(t', \kappa_\d\ln t',+\infty)$, we have 
$$|\chi(x)-\chi(y)|\leq |x-y|^{1/3}$$ for $x,y\in[0,R]^\d$ such that $|x-y|\leq 1/t'$. If $t$ is large enough so as to become larger than $ \ln t'$ then the above relation is valid for $x,y\in[0,R]^\d$ such that $|x-y|\leq e^{-t}$. Thus we have
\begin{align*}
c^{-1}e^{\d t}\int_{B(x_0,e^{-t})} |\chi(x)-\chi(x_0)|dx\leq &c^{-1}e^{\d t}\int_{B(x_0,e^{-t})} |x-x_0|^{1/3}dx=e^{\d t}\int_0^{e^{-t}}r^{1/3+\d-1}\,dr\leq e^{-t/3} .
\end{align*}
This proves that the process $x_0\mapsto c^{-1}e^{\d t}\int_{B(x_0,e^{-t})} (\chi(x)-\chi(x_0))dx$ converges uniformly over $[0,R]^\d$  towards $0$. It remains to treat the (Gaussian) process
$$I_t(x_0)=c^{-1}e^{\d t}\int_{B(x_0,e^{-t})}  (X_t(x)-X_t(x_0))dx.$$
We will use the Kolmogorov criterion to prove the uniform convergence of this process towards $0$. Let us first compute the variance  \begin{align*}
 \E[I_t(x_0)^2 ]&=c^{-2}e^{2\d t}\int_{B(x_0,e^{-t})}\int_{B(y_0,e^{-t})}\E[(X_t(x)-X_t(x_0))(X_t(x')-X_t(x_0))]\,dxdx'\\
 &=c^{-2}e^{2\d t}\int_{B(x_0,e^{-t})}\int_{B(y_0,e^{-t})}\int_1^{e^t}\frac{k(u(x-x'))-k(u(x-x_0))-k(u(x_0-x'))+1}{u} \,du\,dxdy
\end{align*}
 in such a way that we get by using assumption [A.3] (for some irrelevant constant $C$ that may change along lines) and the relation $|x-x'|\leq |x'-x_0|+|x-x_0|$
\begin{align*}
 \E[I_t(x_0)^2]
 &\leq Ce^{2 \d t}\int_{B(x_0,e^{-t})}\int_{B(x_0,e^{-t})}\int_1^{e^t}\frac{u|x-x_0|  }{u} \,du\,dxdx' =  Ce^{ \d t} \int_{B(x_0,e^{-t})}\int_1^{e^t} |x-x_0|  \,dx \\
 &\leq Ce^{-t}.
\end{align*}
Following similar computations, one can also establish that $ \E[(I_t(x_0)-I_t(y_0))^2]\leq C|x_0-y_0|$. This entails, as the process is Gaussian, that for $q>2$
$$ \E[(I_t(x_0)-I_t(y_0))^q]\leq C \E[(I_t(x_0)-I_t(y_0))^2]^{q/2}\leq C|x_0-y_0|^{q/2}$$ where $C$ is a constant that does not depend on $t$. The Kolmogorov criterion then ensures that the family $(I_t(\cdot))_t$ is tight in the space of continuous functions on $[0,R]^\d$ equipped with the topology of uniform convergence. As $ \E[I_t(x_0)^2]\to 0$ as $t\to\infty$, we deduce the convergence in law of this family in the same space towards $0$. The statement of the lemma is then a straightforward consequence.\qed

%%%%%%%%%%%%%%%%%%%%%%%%%%%%%%%%%%%%%%%%%%%%%%%%%%%%%%%%
%%%%%%%%%%%%%%%%%%%%%%%%%%%%%%%%%%%%%%%%%%%%%%%%%%%%%%%%
\section{Proof of Proposition \ref{trlaplace}}\label{sec:prop}
%%%%%%%%%%%%%%%%%%%%%%%%%%%%%%%%%%%%%%%%%%%%%%%%%%%%%%%%%%%%%%%%%%%%%%%%%%%%%%%%%%%%%%%%%%%%%%%%%%%%%%%%%%%%%% 
Our aim is to study for $t,\, t'$ large and $\chi(\cdot)\in \mathcal{C}_{R}(t',\kappa_\d \ln t',\ln t)$,
\begin{equation}
\E\Big[ \exp\big(- \int_{[0,R]^\d}e^{ \gamma[Y_t(x)-a_t-\chi(x)]+\d t} dx \big) \Big].
\end{equation}
According to the Proposition \ref{partieneg}, for $A$ large enough, we can restrain our study to the expectation of 
\begin{equation}
  \Phi^{(A)}(\chi(\cdot),t ):=   \exp\Big(- \int_{[0,R]^\d}e^{ \gamma[Y_t(x)-a_t-\chi(x)]+\d t}\1_{\{  Y_t(x)-a_t-\chi(x)\geq -A \}} dx \Big) .
\end{equation}
Throughout this section, keep in mind that the function $  \Phi^{(A)}(\chi(\cdot),t )$ is bounded by $1$. We fix $R,A,\, \epsilon>0$. We stick to the notations introduced in \cite{Mad13} and we define\begin{align}
\label{def1}
M_{t,\chi} :=&\underset{y\in [0,R]^\d}{\sup} (Y_t(y)-\chi(y)), & &\o_{t,\chi}:=\{y\in [0,R]^\d,\, Y_t(y)\geq a_t+\chi(y)-1\},\\
\label{def2}
M_{t,\chi}(x, b):= &\underset{y\in B(x,e^{b-t})}{\sup}(Y_t(y)-\chi(y)),& &\o_{t,\chi}(x,b):=\{y\in B(x,e^{b-t}),\, Y_t(y)\geq a_t+\chi(y)-1\},\\
\label{def3}  \r_t:=& [e^{-t/2},R-e^{-t/2}]^\d.&&
\end{align}

Observe that on the set $\{ M_{t,\chi-A}< a_t\}$, $1-\Phi^{(A)}(\chi(\cdot),t )=0 $. Moreover for any $t>0$, because of the continuity of the function $x\mapsto Y_t(x)-\chi(x)$, the random variables $|\o_{t,\chi-A}| $  and $|\o_{t,\chi-A}(x,b)|$  are strictly positive respectively on $\{ M_{t,\chi-A}\geq a_t\} $ and $ \{M_{t,\chi-A}(x,b)\geq a_t\} $ (recall that $|B|$ stands for the Lebesgue measure of the set $B\subset \r^\d$). Therefore for any $L\geq 1$, 
\begin{eqnarray}
\nonumber &&\E\Big[ 1-\exp\Big(- \int_{[0,R]^\d}e^{ \gamma[Y_t(x)-a_t-\chi(x)]+\d t}\1_{\{  Y_t(x)-a_t-\chi(x)\geq -A \}} dx \Big) \Big]
\\
\nonumber &&=\E\Big[1-\phi^{(A)}(\chi(\cdot),t)];  M_{t,\chi-A}\geq a_t \Big]
\\
\label{harendt} &&= \E\Big[ \int_{[0,R]^\d}\frac{\1_{\{ m\in \o_{t,\chi-A}\}}\1_{\{M_{t,\chi-A}\geq a_t \}}}{ |\o_{t,\chi-A} |} \big[1-\phi^{(A)}(\chi(\cdot),t)\big]  dm  \Big]:=\E_{(\ref{harendt})}.
\end{eqnarray}
Now we want to exclude the particles $m\in \o_{t,\chi-A}$ such that their paths $Y_\cdot(m)$ are unlikely. We set 
\begin{equation}
\label{squareblan}
 \square_t^{\alpha,A,L}:= \left\{ (f_s)_{s\geq 0},\, \sup_{s\in [0,\frac{t}{2}]}f(s) \leq \alpha,\, \sup_{s\in[\frac{t}{2},t]} f(s) \leq a_t+ \alpha+L,\, f_t \geq a_t+\alpha-A-1 \right\},\quad \forall L,\,\alpha,\, t>0.
\end{equation} 

\begin{Lemma}\label{Barca}
For any $A,\, \epsilon>0$ there exists $L>0$ such that for any $t',T>0$ large enough we have for any $t\geq T$, $\chi \in \mathcal{C}_R(t',\kappa_\d \ln t', \ln t)$,
\begin{equation}
\label{eqBarca} \P\left( \exists m\in \o_{t,\chi-A}\cap [0,R]^\d,\, Y_\cdot(m)\in  \square_t^{\chi(m),A,L}\right) \leq \epsilon I_\d(\chi),\qquad 
\P\left( \exists m\in [0,R]^\d/_{\r_t},\, m\in \mathcal{D}_{t,\chi-A}\right) \leq \epsilon I_\d(\chi).
\end{equation}
\end{Lemma}
In \cite{Mad13} the inequalities of (\ref{eqBarca}) are proved for $A=0$ (via Proposition 4.4 \cite{Mad13} and the arguments of \cite{Mad13} to obtain (5.11)) but it doesn't make any difficulties to extend for any fixed $A>0$, thus we do not detail the proof of Lemma \ref{Barca}.
\\

Going back to (\ref{harendt}), from Lemma \ref{Barca}, we deduce that, for any $A$, there exist $L>0$, $t_0>0$  such that for any $t'\geq t_0$ there exists $T>0$ such that $\forall t\geq T,\, \chi(\cdot) \in \mathcal{C}_R(t', \kappa_\d\ln t',\ln t)$, 
\begin{equation}
\label{dbleharen}
\Big| \E\Big[\int_{\r_t}  \frac{\1_{\{m\in \o_{t,\chi-A},\, Y_\cdot(m)\in \square_t^{\chi(m),A,L}\}}\1_{\{M_{t,\chi-A}\geq a_t \}}}{ |\o_{t,\chi-A} |} [1-\phi^{(A)}(\chi(\cdot),t)]     dm\Big]-\E_{(\ref{harendt})}\Big|\leq \epsilon I_\d(\chi).
\end{equation}
Now the constant $L$ is also fixed.
\medskip

For any $t>b\geq 0$, let us introduce:
\begin{equation}
\label{defXi}
\Xi_{\chi-A,t}(b,m)=\{ \exists y\in [0,R]^\d,\, |y-m|\geq e^{b-t},\, Y_t(y)\geq a_t+\chi(y)-A-1\}.
\end{equation}
 On the complement of $   \Xi_{\chi-A,t}(b,m)  $, we have (just observe that everything happens inside the ball $B(m,e^{b-t})$)
 $$\frac{\1_{\{M_{t,\chi-A}\geq a_t\}}}{|\o_{t,\chi-A}|}=\frac{\1_{\{M_{t,\chi-A}(m,b)\geq a_t\}}}{|\o_{t,\chi-A}(m,b)|} .$$
Also, still on the complement of $   \Xi_{\chi-A,t}(b,m)  $, the function $ [1-\phi^{(A)}(\chi(\cdot),t)] $ is equal to 
 \begin{equation}
 1- \exp\Big(- \int_{B(m,e^{b-t})}e^{ \gamma[Y_t(x)-a_t-\chi(x)]+\d t}\1_{\{ Y_t(x)-a_t-\chi(x)\geq -A \}} dx \Big) :=1-\phi^{(A,b)}(\chi,t,m).
\end{equation}
Therefore for any $b\geq 1,\, m\in \r_t$ we can write,
\begin{eqnarray*}
&&[1-\phi^{(A)}(\chi(\cdot),t)]\frac{\1_{\{M_{t,\chi-A}\geq a_t\}}}{|\o_{t,\chi-A}|}=   [1-\phi^{(A)}(\chi(\cdot),t)]\frac{\1_{\{M_{t,\chi-A}\geq a_t\}}}{|\o_{t,\chi-A}|}(\1_{\{ \Xi_{\chi-A,t}(b,m)^c\}}+ \1_{\{\Xi_{\chi-A,t}(b,m) \}})
\\
&&= (1-\phi^{(A,b)} (\chi,t,m) )  \frac{\1_{\{M_{t,\chi-A}(m,b)\geq a_t\}}}{|\o_{t,\chi-A}(m,b)|}\1_{\{\Xi_{\chi-A,t}(b,m)^c \}} + [1-\phi^{(A)}(\chi(\cdot),t)]\frac{\1_{\{M_{t,\chi-A}\geq a_t\}}}{|\o_{t,\chi-A}|} \1_{\{ \Xi_{\chi-A,t}(b,m)\}}
\\
&&=  (1-\phi^{(A,b)} (\chi,t,m) )  \frac{\1_{\{M_{t,\chi-A}(m,b)\geq a_t\}}}{|\o_{t,\chi-A}(m,b)|} - (1-\phi^{(A,b)} (\chi,t,m) )  \frac{\1_{\{M_{t,\chi-A}(m,b)\geq a_t\}}}{|\o_{t,\chi-A}(m,b)|}\1_{\{ \Xi_{\chi-A,t}(b,m)\}}
\\
&&\qquad \qquad\qquad \qquad \qquad\qquad \qquad \qquad\qquad \qquad  +[1-\phi^{(A)}(\chi(\cdot),t)]\frac{\1_{\{M_{t,\chi-A}\geq a_t\}}}{|\o_{t,\chi-A}|} \1_{\{ \Xi_{\chi-A,t}(b,m)\}}.
\end{eqnarray*}
Following this decomposition, the first expectation in (\ref{dbleharen}) is equal to the sum of 
\begin{eqnarray}
\label{(1)}\qquad(1)_{b}&:=&\E\Big[\int_{\r_t}\frac{\1_{\{Y_\cdot(m) \in \square_t^{\chi(m),A,L}\}}\1_{\{M_{t,\chi-A}(m,b)\geq a_t \}} }{  |\o_{t,\chi-A}(m,b) |} [1-\phi^{(A,b)}(\chi(\cdot),t,m)]dm\Big],
\\
\label{(2)}(2)_{b}&:=& \E\Big[\int_{\r_t} \1_{\Xi_{\chi-A,t}(b,m)}  \frac{\1_{\{  Y_\cdot(m)\in \square_t^{\chi(m),A,L}\}}  \1_{\{M_{t,\chi-A}\geq a_t \}}}{ |\o_{t,\chi-A}|} [1-\phi^{(A)}(\chi(\cdot),t)] dm\Big],
\\
\label{(3)}(3)_{b}&:=&-\E\Big[\int_{\r_t} \1_{\Xi_{\chi-A,t}(b,m)}  \frac{\1_{\{  Y_\cdot(m)\in \square_t^{\chi(m),A,L}\}}\1_{\{M_{t,\chi-A}(m,b)\geq a_t  \}}}{ |\o_{t,\chi-A}(m,b) |} [1-\phi^{(A,b)}(\chi(\cdot),t,m)] dm\Big].
\end{eqnarray}

%We shall show, via two lemmas, that $(2)_{L,A,b}$ and $(3)_{L,A,b}$ are negligible. 
%\begin{Lemma}
%\label{avantXi}
%There exists $c_8>0$ such that for any $L,\,  t\geq b>1$,\nomenclature[i5]{ $ {\bf r}_t(x)$ }{$ :=\sup\{ r>0,\, w_{Y_\cdot(\cdot)}(r,x ,t)\leq \frac{1}{4}\}  \wedge  e^{-t} $}
%\begin{eqnarray}
%\label{eqavantXi}(2)_{L,A,b}+(3)_{L,A,b}\leq c_8\int_{[0,R]^\d}\E\left(\frac{\1_{\{Y_\cdot(m) \in \DDi_t^{\chi(m),L}\}}}{{\bf r}_t(m)^\d}\1_{\{ \Xi_{\chi-A,t}(b-\ln 2,m)\}}\right)dm, \text{ with}
%\\
%\label{defrx}{\bf r}_t(m):=\sup\{ r>0,\, w_{Y_\cdot(\cdot)}(r,m ,t)\leq \frac{1}{4}\}  \wedge  e^{-t}.
%\end{eqnarray}
%\end{Lemma} 

\begin{Lemma}
\label{Barca2}
For any $A,L,\epsilon>0$, there exists $b_0,t_0$ large enough such that for any $t'\geq t_0$, $b\geq b_0$, $\exists T>0$ such that for any $t\geq T$, $\chi\in \mathcal{C}_R(t',\kappa_\d \ln t',\ln t)$ we have
\begin{equation}
 |(2)_{b}|+|(3)_{b}|\leq \epsilon I_\d(\chi).
\end{equation}
\end{Lemma}
We do not detail the proof of Lemma \ref{Barca2} but, recalling that $|1-\phi^{(A,b)}(\chi(\cdot),t,m)|$ and $|1-\phi^{(A)}(\chi(\cdot),t)|$ are bounded by $1$, we just remark that the amounts $(2)_{b}$ and $(3)_{b}$ are very similar to $(2)_{L,b}$ and $(3)_{L,b}$ defined in (5.15) and (5.16) of \cite{Mad13}. Then Lemma \ref{Barca2} is a minor adaptation  of the proofs of Lemmas 5.1 and 5.2 in \cite{Mad13} (in \cite{Mad13} $A=0$, whereas here $A$ is a fixed positive constant).
\\

Thus combining Lemma \ref{Barca2} and (\ref{dbleharen}), we deduce that there exist $b$ and $t_0 >0$, such that for any  $t'>t_0$ there exists $T>0$ such that $\forall\, t\geq T,\, \chi(\cdot) \in \mathcal{C}_R(t',\kappa_\d \ln t',\ln t)$, 
\begin{equation}
\label{inter2}
\left| \E_{(\ref{harendt})}-(1)_{b}\right| \leq  2\epsilon \int_{[0,R]^\d}   \chi(x)  e^{-\sqrt{2\d}\chi(x)} dx.
\end{equation}

Therefore we can  restrain our study to $(1)_{b}$ (with $A,L,b$ fixed). The Markov property at time $t_b=t-b$ and the invariance by translation of $(Y_s(x))_{s\geq 0,\, x\in \r^\d}$ give:
\begin{eqnarray}
\nonumber (1)_{b}&=&\E\Big[\int_{\r_t}\frac{\1_{\{  Y_\cdot(m) \in \square_t^{\chi(m),A,L},\, m\in \o_{t,\chi-A}\}} \1_{\{ M_{t,\chi-A}(m,b)\geq a_t\}}  }{ |\o_{t,\chi-A}(m,b)|}[1-\phi^{(A,b)}(\chi(\cdot),t,m)]  d m\Big]
\\
\label{tito}&=&   \int_{\r_t} \E\Big[\1_{\{ \sup_{s\in[0,t_b]}Y_s(m)\leq \chi(m),\, \sup_{s\in[\frac{t}{2},t_b]}Y_s(m)\leq a_t+\chi(m)+L\}} D^{b}_{m,t}\Big]dm,
\end{eqnarray}
%
%\E_{Y_{t_b}(x)-a_t(\chi(x)+L)}(...)_{\Big|{\begin{matrix}
%&g(y)=\Y_{t_b,x}(y)-f_x(y)&
%\\
%&y\in B(0,e^{b-t})&
%\end{matrix}          }}
where
\begin{eqnarray}
 \nonumber D^{b}_{m,t}&:=&\E \Big[ \frac{\1_{\{ \sup_{s\in[0,t_b]}Y^{(t_b)}_{s}(0)+ \bar{z}\leq 0,\, Y_{t_b,b}(0)\geq -L-A-1,\, \exists y\in B(0,e^{b-t}),\, Y^{(t_b)}_{b}(y)+\bar{z}\geq -L-A - g (y) \}}}{|B(0,e^{b-t}) \cap \{y: Y^{(t_b)}_b(y)+\bar{z}\geq -L-A-1-g(y)\} |}  
  \\
\nonumber && \quad\quad\quad \quad\Big( 1-\exp\{ -\int_{B(0,e^{b-t})}e^{ \gamma[Y_b^{(t_b)}(y)+\bar{z}+g (y)+L ]+\d t}\1_{\{ Y_b^{(t_b)}(y) \geq -A-L -g(y)\}} dy \} \Big)\Big] , 
\end{eqnarray}
with
$$ g(y)=Y_{t_b}(m+y)-Y_{t_b}(m)-(\chi(m+y)-\chi(m)),\hspace{2cm} \bar{z} = Y_{t_b}(m)-a_t -\chi(m)-L .$$

\medskip
 In the following we will denote
\begin{equation}
\forall m\in \r_t,\, \chi_m(.):=\chi(m+.)-\chi(m) .
 \end{equation} 
According to the scaling property $\left(Y_s^{(t_{b})}(y)\right)_{s\leq b,\,y\in B(0,e^{b-t})}\overset{(d)}{=} (Y_s(ye^{t-b}))_{s\leq b,\,y\in B(0,e^{b-t})}$, thus we can rewrite $ D_{m,t}^{b} $ as
\begin{eqnarray*}
&&e^{\d t_b}\E_{\bar{z}}\Big[ \frac{\1_{\{ \sup_{s\in[0,b]}Y_s(0)\leq 0,\, Y_b(0)\geq -L-A-1\}} \1_{\{\exists y\in B(0,1),\, Y_{b}(y)\geq -L-A- g(ye^{b-t}) \}}}{|B(0,1) \cap \{y: Y_b(y)\geq -L-A-1-g(ye^{b-t})\} |} 
 \\
\nonumber && \hspace{1.5cm}\Big(1- \exp\Big(-\int_{B(0,1)}e^{ \gamma[Y_b(y)+g  (ye^{b-t}) +L]}\1_{\{ Y_b(y) \geq -A-L -g (ye^{b-t})\}} dy \Big)  \Big)\Big],
\end{eqnarray*}
where we have used the convention: for any $z\in \r$, $(Y_s(x))_{s\geq 0,x\in \r^\d}$ under $\P_z$ has the law $(z+Y_s(x))_{s\geq 0,\, x\in \r^\d}$ under $\P$.  Then Lemma \ref{rRhodesdec} and the Girsanov transformation lead to 
\begin{align*}
 (1)_{b}&=&\int_{\r_t} \E\Big[e^{\sqrt{2\d}Y_{t_b}(m)+\d t_b}\1_{\{ \sup_{s\in[0,t_b]}Y_s(m)\leq \chi(m),\, \sup_{s\in[\frac{t}{2},t_b]}Y_s(m)\leq a_t+\chi(m)+L\}}   e^{-\sqrt{2\d}Y_{t_b}(m)-\d t_b}e^{\d t_b}  D^{b}_{m,t} \Big]dm
\\
&=&\int_{\r_t} {e^{-\sqrt{2\d}\chi(m)}}{t^\frac{3}{2}}\E_{-\chi(m)}\Big[\1_{\{ \sup_{s\in[0,t_b]}B_s\leq 0,\sup_{s\in[\frac{t}{2},t_b]}B_s\leq a_t+L\}}F\big(B_{t_b}-a_t -L,\mathfrak{G}_{t,b}^{\chi_m}\big)\Big]dm,%\label{egal(1)L} 
\end{align*}
where   $B$ a standard Brownian motion and, for $ g\in \mathcal{C}(B{(0,1)},\r),\, z\in \r  $,
\begin{eqnarray}
&&\label{defFLb}F(z,g):=e^{-\sqrt{2\d}(z+L)}\E_{z}\left[ \frac{\1_{\{ \sup_{s\in [0,b]}Y_s(0)\leq 0,\, Y_b(0)\geq -L-A-1\}} \1_{\{\exists y\in B(0,1),\, Y_{b}(y)\geq -L-A- g(ye^{b}) \}}}{|B(0,1) \cap \{y: Y_b(y)\geq -L-A-1-g(ye^b)\} |}\right.
 \\
\nonumber &&\left. \times \Big(1- \exp\big(- \int_{B(0,1)}e^{ \gamma[Y_b(y)+g (ye^{b})+L ]}\1_{\{ Y_b(y) \geq -A-L -g(ye^{b})\}} dy \big)  \Big)\right],  
\end{eqnarray}
and for any $\Psi \in \mathcal{C}_R(B(0,e^b),\r)$, 
 \begin{equation}
\label{defG} \mathfrak{G}_{t,b}^\Psi:B{(0,e^{b})}\ni y\mapsto   -\int_{0}^{t_b}(1-k(e^{s-t} y))dB_{s} - \zeta_t(ye^{-t})+ Z_{t_b}^{0}(ye^{-t})+\Psi(ye^{-t}).
 \end{equation}
For $\Psi=0$ we denote $\mathfrak{G}_{t,b}^0=\mathfrak{G}_{t,b}$. In passing we take the opportunity to define for any $\sigma \in [0,t_b]$, 
\begin{equation}
\label{defGsigma}
\mathfrak{G}_{t,b,\sigma} :B{(0,e^{b})}\ni y\mapsto   -\int_{t_b-\sigma}^{t_b}(1-k(e^{s-t} y))(e^{s-t} y)dB_{s} -  \zeta_t(ye^{-t})+ Z_{t_b}^{0}(ye^{-t})
\end{equation}
and the processes $\zeta,Z $ are defined in Lemma \ref{rRhodesdec}. Note that $Z^0_{t_b}(\cdot)  $ is a centered Gaussian process, independent of $B$, which has the covariances as in  \cite[equation (2.6)]{Mad13}. Furthermore by  \cite[Proposition 2.4]{Mad13} for  any $b>0$, the  Gaussian process $B{(0,e^b)}\ni y\mapsto Z_{t-b}^0({ye^{-t}})$, converges in law to $B(0,e^b)\ni y\mapsto Z(y)$ when $t$ goes to infinity.

\medskip

%Now we want to get (via a renewal theorem): for any $L,A,b>0$, uniformly in $x\in \r_t $,
Finally with our new notations, we have to study for any $m\in \r_t$, 
$$\E_{-\chi(m)}\left(\1_{\{ \sup_{s\in [0,t_b]}B_s\leq 0,\,  \sup_{s\in [\frac{t}{2},t_b]}B_s\leq a_t+L\}}F\left(B_{t_b}-a_t-L,\mathfrak{G}_{t,b}^{\chi_m}\right)\right).$$
%The crucial point is to have a constant $C^*$ which does not depend of $x$ or $\chi$. 
Recalling Theorem \ref{trlaplace}, our goal is to prove that this quantity is equivalent to a constant times $t^{-\frac{3}{2}}\chi(m)$, when $t$ goes to infinity. We will prove this by the renewal theorem below. To state our result we need the following definition.

\begin{Definition}
A continuous function $F:\r\times \mathcal{C}(B(0,e^b),\r)\to \r^+$ is ''{\bf $b$ regular}" if there exists two functions $h:\r\to \r_+$ and $F^*: \mathcal{C}(B(0,e^b))\to \r^+$ satisfying 

(i) 
\begin{equation}
\label{et2}
\underset{x\in \r }{\sup}\, h(x)<+\infty ,\quad  \text{and } h(x)\underset{x\to -\infty}{=}o(e^x).
\end{equation}

(ii) There exists $c>0$ such that for any $ \delta\in (0,1)$, $g\in \mathcal{C}(B{(0,e^b)},\r)$ with $w_{g(\cdot e^b)}(\delta)\leq \frac{1}{4} $, 
\begin{equation}
\label{et1}
 F^*(g)\leq c{\delta^{-10}}.
\end{equation}

(iii) For any $z\in \r,\, g\in \mathcal{C}(B(0,e^b),\r)$,  $F(z,g)\leq h(z)  F^*(g) $.

(iv) There exists $c>0$ such that for any $z\in \r,\, g_1,g_2\in \mathcal{C}(B(0,e^b),\r)$ with  $||g_1-g_2||_\infty\leq \frac{1}{8} $,   
\begin{equation}
\label{et3bibi}
|F(z,g_1)-F(z,g_2)|\leq c_9 ||g_1-g_2||_\infty^\frac{1}{4} h(z) F^*(g_1).
\end{equation}

\end{Definition}

We will prove the following two results at the end of the section.

\begin{Lemma}[Control of $F$]
\label{controlF} $F$ is {\bf $b$  regular}.
\end{Lemma} 

For any $M\geq 0$ and $F$ a function $b$ regular, we define
\begin{equation}
\label{defMreload}F^{(M)}(x,g):=(F(x,g)\wedge M)  \1_{\{x\geq -M\}},\quad   \text{ and } \, \overline{F}^{(M)}=F-F^{(M)}.
\end{equation}
and 
\begin{theorem}
\label{RESUME}
Let $b>0$ and $F:\r\times \mathcal{C}(B(0,e^b),\r)\to \r^+$ be a function {\bf $b$ regular}. For any $\epsilon>0$,  there exists $M,\sigma, T >0$ such that for any $t\geq T$, $\chi(\cdot) \in \mathcal{C}_R(t',\kappa_d \ln t',\ln t) $, $\mathtt{z}\in [1,\ln t)^{30}$,
\begin{equation}
\label{eqRESUME}
\left| \int_{\r_t} e^{-\sqrt{2\d}\chi(x)} \E_{-\chi(x)} \left(\1_{\{ \sup_{s\in [0,t_b]}B_s\leq 0,\,  \sup_{s\in [\frac{t}{2},t_b]}B_s\leq -\mathtt{z} \}}  F\left(B_{t_b}+\mathtt{z},\mathfrak{G}^{\chi_x}_{t,b}\right)\right)dx- C_{M,\sigma}(F) I_\d(\chi)\right|\leq \epsilon I_\d(\chi).
\end{equation}
with
\begin{eqnarray}
\nonumber C_{M,\sigma}(F):= C\int_{0}^M\int_0^u\E\Big(  F^{(M)}\Big(-u, y\mapsto Z (ye^{-b})-\zeta(ye^{-b}) \qquad\qquad\qquad\qquad\qquad
\\
\label{defCLBMS}  -\int_{0}^{T_{-\gamma} \wedge \sigma }(1-k(e^{-s}ye^{-b})) dB_s - \int_{T_{-\gamma}\wedge \sigma}^\sigma (1-k(e^{s}ye^{-b}) )d \mathtt{R}_{s-T_{-\gamma}}\Big)\Big)d\gamma dU.
\end{eqnarray}
\end{theorem}

Now we are in position to complete the proof Theorem \ref{trlaplace}.  Indeed by combining Proposition \ref{partieneg}, inequalities (\ref{dbleharen}), (\ref{inter2}), Lemma \ref{controlF} and Theorem \ref{RESUME} we deduce that: {\it $\forall \epsilon>0$ there exist $A,\, L,\, b,\,M,\, \sigma>0$ such that for $t',T>0$ large enough we have : {\bf for any $t\geq T,\, \chi(\cdot)\in \mathcal{C}_R(t',\kappa_\d \ln t',\, \ln t)$,
\begin{equation}
\label{eqinterfinal}
\Big|\E\Big[ 1-\exp\big(- \int_{[0,R]^\d}e^{ \gamma[Y_t(x)-a_t-\chi(x)]+\d t}dx \big) \Big]  -C_{M,\sigma}(F)\int_{[0,R]^\d}   \chi(x)  e^{-\sqrt{2\d}\chi(x)} dx\Big|\leq \epsilon\int_{[0,R]^\d}   \chi(x)  e^{-\sqrt{2\d}\chi(x)} dx.
\end{equation}}}

In addition by Proposition \ref{uplowbound}:  there exist $c_2>0$ and $t',T>0$ large enough such that  for any $t\geq T$ and   $\rho(\cdot) \in \mathcal{C}_R(t', \kappa_\d \ln t',\ln t)$, 
\begin{equation}
\label{utillower}   \E\Big[1- \exp\big(- \int_{[0,R]^\d}e^{ \gamma[Y_t(x)-a_t-\rho(x)]+\d t}dx \big) \Big]\leq c_2  \int_{[0,R]^\d}   \chi(x)  e^{-\sqrt{2\d}\chi(x)} dx .
\end{equation}
For any $ n>0$, let $(L_n,b_n,M_n,\sigma_n)$ such that (\ref{eqinterfinal}) is true with $\epsilon=\frac{1}{n}$. Clearly $C_n:= C_{M_n,\sigma_n}(F)\in [0, 2c_2]$ for any $n\in \N$ (notice that $F$ depends on $L_n,A_n,b_n$ though it does not appear through the notations). Let $\phi:\N\to \N$ strictly increasing such that $C_{\phi(n)}\to C^*\in [0,2c_2]$ as $n\to\infty$. 

Now we fix $\epsilon>0$. Let $N_0>0$ such that for any $n\geq N_0$, $|C_{\phi(n)}-C(\gamma) |\leq \epsilon$. Then we choose $N_1> N_0$ 
such that $n\geq N_1$ implies $\frac{1}{\phi(n)}\leq \epsilon$.   Finally there exist (according to Proposition \ref{uplowbound}) $t'(=t'(N_1))$ and $T(=T(N_1))>0$ such that for any $t\geq T,\, \chi (\cdot)\in \mathcal{C}_R(t',\kappa_\d \ln t',\ln t)$,
\begin{eqnarray*}
\Big|   \E\Big[ 1-\exp\big(- \int_{[0,R]^\d}e^{ \gamma[Y_t(x)-a_t-\chi(x)]+\d t}dx \big) \Big]-C(\gamma) \int_{[0,R]^\d}   \chi(x)  e^{-\sqrt{2\d}\chi(x)} dx \Big|\leq \epsilon \int_{[0,R]^\d}   \chi(x)  e^{-\sqrt{2\d}\chi(x)} dx .
\end{eqnarray*} 
To complete the proof of Theorem \ref{trlaplace} it remains to prove that $C(\gamma)>0$. It is a consequence of Proposition \ref{uplowbound}. Indeed let $t'>0$ large and $\chi\in \mathcal{C}_R(t',\kappa_d\ln t',+\infty)$ such that for any $t>T$, 
\begin{eqnarray*}
\Big|   \E\Big[ 1-\exp\big(- \int_{[0,R]^\d}e^{ \gamma[Y_t(x)-a_t-\chi(x)]+\d t}dx \big) \Big]-C(\gamma)  \int_{[0,R]^\d}   \chi(x)  e^{-\sqrt{2\d}\chi(x)} dx|\leq \frac{c_1}{2} \int_{[0,R]^\d}   \chi(x)  e^{-\sqrt{2\d}\chi(x)} dx,
\end{eqnarray*}
with $c_1$ the constant defined in Proposition \ref{uplowbound}. From  Proposition \ref{uplowbound}, we have $$\liminf_{t\to\infty} \E\Big[ 1-\exp\big(- \int_{[0,R]^\d}e^{ \gamma[Y_t(x)-a_t-\chi(x)]+\d t}dx \big) \Big]\geq c_1 \int_{[0,R]^\d}   \chi(x)  e^{-\sqrt{2\d}\chi(x)} dx,$$ then it is plain to deduce $C(\gamma) \geq \frac{c_1}{2}>0$, which completes the proof of Theorem \ref{trlaplace}.\qed

%%%%%%%%%%%%%%%%%%%%%%%%%%%%%%%%%%%%%%%%%%%%%%%%%%%%
%%%%%%%%%%%%%%%%%%%%%%%%%%%%%%%%%%%%%
\subsection{Proof of Lemma \ref{controlF}}
Recall the convention: for any $z\in \r$, $(Y_s(x))_{s\geq 0,x\in \r^\d}$ under $\P_z$ has the law $(z+Y_s(x))_{s\geq 0,\, x\in \r^\d}$ under $\P$.
\\

 \noindent{\it Proofs of Lemma \ref{controlF}} Fix $L,b>1$, recall (\ref{defFLb}) for the definition of $F$. We shall prove that $F$ is {\bf $b$  regular} with  
\begin{eqnarray}
\label{defhLb}
h=h_{L,b}(z)&:=&  e^{-\sqrt{2\d}(z+L)}\P_{z+L+1}\Big( Y_b(0)\geq 0\Big)^{\frac{1}{2}} ,
\\
\label{defF}F^*=F^*_{b}(g)&:=& \underset{z\in \r}{\sup}\, \E_z\Big[\frac{\1_{\{\exists y\in B(0,1),\, Y_b(y)\geq  -g(ye^b)  \}} }{ |B(0,1) \cap\{ y,\, Y_b(y)\geq -g(ye^b)  -\frac{1}{2}\}|^8}\Big]^{\frac{1}{4}} .
\end{eqnarray}

Check (i) is an elementary computation whereas (iii) stems from the Cauchy-Schwartz inequality. Let us start by showing that $F^*_b$ satisfies (\ref{et1}).
Let $g\in \mathcal{C}(B(0,e^b),\r)$ such that $w_{g(\cdot  e^b)}(\delta)\leq \frac{1}{4}$. We define 
$$\Lambda=|B{(0,1)}\cap \{y, Y_b(y)\geq -g(ye^b)-\frac{1}{2}\}  | .$$
%\begin{eqnarray}
%\Lambda=\lambda_{B{(0,1)}}(\{y, \Y_b(y)\geq g(ye^b)-\frac{1}{2}\}  ).
%\end{eqnarray} 
On the set $\{ \exists y\in B(0,1),\, Y_b(y)\geq -g(ye^{b}) \}$, we introduce
$${\bf r}=\sup\Big\{{\bf s};\exists x_{ \mathbf{s}} \text{ with }B(x_{ \mathbf{s}} , { \mathbf{s}} )\subset B(0,1), \exists z_{\bf s}\in B(x_{\bf s}, {\bf s})\text{ with }Y_b(z_{\bf s})\geq g(z_{\bf s}e^b),  \forall y\in B(x_{\bf s}, {\bf s}), Y_b(y)\geq -g(ye^b)-\frac{1}{2} \Big\}.$$
 Observe that
\begin{eqnarray*}
F^*_{b}(g)^4 &=&  \underset{x\in \r}{\sup}\,\E_x\Big[\frac{\1_{\{\exists y\in B(0,1),\, Y_b(y)\geq -g(ye^{b}) \}}}{\Lambda^8}\Big]
   \\
 &&  \leq  (e^{b}/\delta)^8+ \underset{k=e^{b}/\delta}{\overset{\infty}{\sum}}  (k+1)^{8}  \underset{x\in \r}{\sup}\,\E_{x}\Big[   \1_{\{\exists y\in B(0,1),\, Y_{b}(y)\geq -g(ye^b) \}} \1_{\{\frac{1}{k+1}\leq \Lambda \leq \frac{1}{k}\}}\Big] .
\end{eqnarray*}
Clearly, $\Lambda\leq (\frac{1}{k})^\d$ implies $  {\bf r}\leq \frac{1}{k}$, moreover $\{ {\bf r} \leq \frac{1}{k}<\delta \}$ implies  $ \{\underset{ \underset{|x-y|\leq \frac{1}{k}}{x,y\in B(0,1)}}{\sup}|Y_b(x)-Y_b(y)|\geq \frac{1}{2}-w_{g(.e^b)}(\delta) \} $. It follows that
\begin{eqnarray*}
\P\left( {\bf r} \leq \frac{1}{k}<\delta\right)& \leq &\P\Big( \underset{ \underset{|x-y|\leq \frac{1}{k}}{x,y\in B(0,1)}}{\sup}|Y_b(x)-Y_b(y)|\geq \frac{1}{2}-w_{g(.e^b)}(\delta)\Big) 
\\
&\leq & \P\Big(\underset{ \underset{|x-y|\leq \frac{1}{k}}{x,y\in B(0,1)}}{\sup}|Y_b(x)-Y_b(y)|\geq \frac{1}{4}\Big),\qquad \text{(recall that }w_{g(\cdot  e^b)}(\delta)\leq \frac{1}{4}).
\end{eqnarray*}
From  in \cite[equation (3.10)]{Mad13} (with $h=\frac{1}{k},\, m=2k,\, p=2,\, t'=b$ and $x=ce^{-b}k $), we have 
$$\underset{z\in \r}{\sup}\,\P_z\Big( \underset{\underset{|x-y|\leq \frac{1}{k}}{x,y\in B(0,1)}}{\sup}|Y_b(x)-Y_b(y)|\geq \frac{1}{4}  \Big)=\P_0\Big( \underset{\underset{|x-y|\leq \frac{1}{k}}{x,y\in B(0,1)}}{\sup}|Y_b(x)-Y_b(y)|\geq \frac{1}{4}  \Big)\leq c'e^{-\frac{1}{c''}e^{-b} k}.$$ 
Finally
$
F^*_{b}(g)^4\leq e^{8b}/\delta^8 + \underset{k=1+e^b/\delta}{\overset{\infty }{\sum}}  (k+1)^{8}ce^{-\frac{1}{c''}e^{-b} k} \leq  c(b) \delta^{-8}$, which proves (\ref{et1}).

Now it remains to prove (\ref{et3bibi}).  Let $g_1,\, g_2$ two continuous functions from $B(0,e^b)\to \r$ such that $||g_1-g_2||_\infty= \delta<\frac{1}{8}$. Let us define (uniquely for this proof) $\forall g\in \mathcal{C}(B(0,e^b),\r)$ and $\gamma>0$:
\begin{eqnarray*}
&&M(g):=\underset{y\in B(0,1)}{\sup}(Y_b(y)+g(ye^b)), \qquad \Lambda_g(\gamma) :=|B(0,1)\cap\{y,\, Y_b(y)\geq -g(ye^b)+\gamma  \}|.
\end{eqnarray*}
With   two notations  and twice the Cauchy-Schwartz inequality we obtain 
 \begin{eqnarray}
\nonumber |F (z,g_1)-F (z,g_2)|\leq e^{-\sqrt{2\d}(z+L)} \E_{z+L+1}\Big[\1_{\{  Y_b(0)\geq 0\}} \Big|  \frac{\1_{\{ M({g_1})\geq 1 \}}}{\Lambda_{g_1}(0)} - \frac{\1_{\{  M({g_2})\geq 1 \}}}{\Lambda_{g_2}(0)}\Big| \Big]+
 \\
\label{bb2}  h_{L,b}(z)F_b^*(g_1) \E_{z+L+1}\Big[ \Big( \Delta(g_1,g_2) \Big)^8    \Big].
\end{eqnarray}
with
\begin{equation}
\label{diffexpectation}      \Delta(g_1,g_2):= e^{-\int_{B(0,1)} e^{\gamma [Y_b(y)+g_1(ye^b)-1]}\1_{\{Y_b(y)+A \geq1- g_1(ye^b)\}} dy  } -e^{-\int_{B(0,1)}e^{\gamma [Y_b(y)+g_2(ye^b)-1]}\1_{\{Y_b(y)+A \geq 1-g_2(ye^b)\}} dy }.
\end{equation}
Let us treat the first term of \eqref{bb2}. From   \cite[Theorem 3.1]{PittTran}, as $Var(Y_b(y))=b\geq 1,\,\forall y\in B(0,1)  $,  we can affirm that  there exists $c>0 $ such that for any $\delta\in(0,1) $, $g\in  \mathcal{C}(B(0,e^b),\r)$,
\begin{equation}
\label{PTa}
\underset{z\in \r}{\sup}\, \P\Big(M(g)\in [z-\delta,z+\delta]\Big)\leq c \delta.
\end{equation}
Thus the first term in (\ref{bb2}) is smaller than
\begin{eqnarray*}
&&\leq  \E_{z+L+1}\Big[\frac{\1_{\{ Y_b(0)\geq 0,\, M(g_1)\in [1-\delta,1+\delta]\}}}{\Lambda_{g_1}(0)}  \Big] +  \E_{z+L+1}\Big[\1_{\{ Y_b(0)\geq 0,\, M(g_2)\geq 1 \}}  \frac{\Lambda_{g_1}(-\delta)-\Lambda_{g_1}(\delta)}{\Lambda_{g_1}(0) \Lambda_{g_2}(0)} \Big]
\\
&& := (A)+(B).
\end{eqnarray*}
By applying twice the Cauchy-Schwartz inequality to $(A)$ we get that
\begin{eqnarray*}
(A)&\leq & \P_{z+L+1}\Big( Y_b(0)\geq 0\Big)^{\frac{1}{2}}\times\E_{z+L+1}\Big[\frac{\1_{\{M(g_1)\geq 1  -\delta \}} }{\Lambda_{g_1}(0)^4}\Big]^{\frac{1}{4}} \times \P_{z+L+1}\Big( M(g_1)\in [1-\delta ,1+\delta]\Big)^{\frac{1}{4}} .
\end{eqnarray*}
Now by applying (\ref{PTa}) to the last term we obtain
\begin{eqnarray*}
(A)&\leq& c \P_{z+L+1}\Big( Y_b(0)\geq 0\Big)^{\frac{1}{2}}\times\E_{z+L+\delta}\Big[\frac{\1_{\{ M(g_1)\geq 0  \}} }{\Lambda_{g_1}^4(\delta-1)}\Big]^{\frac{1}{4}}\delta^\frac{1}{4}
\\
&\leq& c ||g_1-g_2||_\infty^\frac{1}{4} h_{L,b}(z) F_b^*(g_1),\qquad (\text{as }\delta -1\leq -\frac{1}{2}).
\end{eqnarray*}
Similarly, observing that  $\min(\Lambda_{g_1}(0),\Lambda_{g_2}(0))\geq \Lambda_{g_1}(\frac{1}{4})$, we deduce that
\begin{eqnarray*}
(B)&=&\int_{B(0,1)}\E_{z+L+1}\Big[\frac{\1_{\{ Y_b(0)\geq 0,\,  M(g_2)\geq 1 \}}}{\Lambda_{g_1}(0)\Lambda_{g_2}(0)}  \1_{\{Y_b(x)+ g_{1}(xe^b) \in [-\delta ,\delta]\}} \Big]dx
\\
&&\leq \P_{z+L+1}\Big( Y_b(0)\geq 0\Big)^{\frac{1}{2}}                     \E_{z+L+1}\Big[\frac{\1_{\{ M(g_1)\geq 1-\delta\}} }{[\Lambda_{g_1}(\frac{1}{4})]^8}\Big]^{\frac{1}{4}} \int_{B(0,1)}\P_{z+L+1+ g_{1}(xe^b) }\Big(  Y_b(x)\in [-\delta,\delta]  \Big)^\frac{1}{4}dx
\\
&&\leq c\P_{z+L+1}\Big( Y_b(0)\geq 0\Big)^{\frac{1}{2}}                     \E_{z+L+1+\delta}\Big[\frac{\1_{\{\exists y\in B(0,1),\, Y_b(y)\geq - g_{1}(ye^b)+1  \}} }{[\Lambda_{g_1}(\frac{1}{4}+\delta)]^8}\Big]^{\frac{1}{4}} \delta^{\frac{1}{4}}
\\
&&\leq c||g_1-g_2||_\infty^\frac{1}{4} h_{L,b}(z) F_b^*(g_1).
\end{eqnarray*}
So we are done with the study of the first term of (\ref{bb2}). Now we treat the second term. By the triangular inequality, $|\Delta(g_1,g_2)| $ is smaller than $(1)+(2)$ with
\begin{eqnarray}
\label{ouah}(1):= \Big| \exp\Big(-\int_{B(0,1)} e^{\gamma[ Y_b(y)+g_1(ye^b)-1]}\1_{\{Y_b(y)+A\geq  1-g_1(ye^b)\}}dy  \Big)   \qquad\qquad 
 \\
\nonumber  -\exp\Big(-\int_{B(0,1)} e^{\gamma [Y_b(y)+g_2(ye^b)-1]}\1_{\{Y_b(y)+A \geq 1-g_1(ye^b)\}} dy \Big) \Big|  ,
\\
\label{ouah2} (2):=  \Big| \exp\Big(-\int_{B(0,1)} e^{\gamma [Y_b(y)+g_2(ye^b)-1]}\1_{\{Y_b(y)+A \geq 1-g_1(ye^b)\}}dy  \Big)  \qquad\qquad  
 \\
\nonumber -\exp\Big(-\int_{B(0,1)} e^{\gamma [Y_b(y)+g_2(ye^b)-1]}\1_{\{Y_b(y)+A\geq 1-g_2(ye^b)\}} dy \Big)  \Big|  .
\end{eqnarray}
Recalling that $||g_1-g_2||_\infty = \underset{x \in B(0,e^b)}{\sup} |g_1(x)-g_2(x)|:=\delta$, in (\ref{ouah}) by forcing the factorization by $\exp (-U*):=  \exp\Big(-   \int_{B(0,1)} e^{\gamma [ Y_b(y)+g_1(ye^b)-1]}\1_{\{Y_b(y)+A\geq  1-g_1(ye^b)\}}dy  \Big) $, we have
\begin{eqnarray}
\label{deltag1} (1)\leq    e^{-U^*}( e^{ (e^{\gamma \delta }-e^{-\gamma \delta})   U^*}  -1)  \leq  e^{\gamma \delta }-e^{-\gamma \delta }.
\end{eqnarray}
Similarly by some elementary computations we get,
\begin{eqnarray*}
(2)\leq \Big[1- \exp\Big(-\int_{B(0,1)} e^{\gamma[Y_b(y)+g_2(ye^b)-1]} \1_{\{ \min(-g_1(ye^b),-g_2(ye^b))\geq Y_b(y)+A\leq 1- g_1(y)   \}} dy \Big)\Big]
\\
+\Big[1- \exp\Big( -\int_{B(0,1)} e^{\gamma[Y_b(y)+g_2(ye^b)-1]} \1_{\{ \min(-g_1(ye^b),-g_2(ye^b))\geq Y_b(y)+A\leq 1-g_1(y)   \}} dy \Big)\Big]
\\
\leq 2\Big[1- \exp\Big( -\int_{B(0,1)} e^{-\gamma(A-\delta)}  \1_{\{ \min(-g_1(ye^b),-g_2(ye^b))\geq Y_b(y)+A\leq 1-g_1(y)+\delta  \}} dy \Big)\Big]
\\
\leq  \int_{B(0,1)} \1_{\{ \min(-g_1(ye^b),-g_2(ye^b))\leq Y_b(y)+A+1\leq \min (-g_1(ye^b),-g_2(ye^b))+\delta  \}} dy,\qquad (\text{for }\delta\leq 1)
\end{eqnarray*}
By the Jensen inequality and recalling that $\underset{y\in B(0,1)}{\sup}\underset{z\in \r}{\sup}\,\P\Big( Y_b(y)\in [z,z+\delta]\Big)\leq \frac{\delta}{\sqrt{b}}$, we deduce that the expectation of $[(1)+(2)]^8$ is smaller than $ c\delta $. Combining this inequality with (\ref{deltag1}) yields
\begin{equation}
\underset{z\in \r}{\sup}\,\E_{z+L+1}\Big[|\Delta(g_1,g_2)|^8\Big]\leq c||g_1-g_2 ||_{\infty}.
\end{equation}

Finally the condition (iv) follows with   
\begin{eqnarray}
h=h_{L,b}(z)&:=&  e^{-\sqrt{2\d}(z+L)}\P_{z+L+1}\Big( Y_b(0)\geq 0\Big)^{\frac{1}{2}} ,
\\
F^*=F^*_{b}(g)&:=& \underset{z\in \r}{\sup}\, \E_z\Big[\frac{\1_{\{\exists y\in B(0,1),\, Y_b(y)\geq  -g(ye^b)  \}} }{[\lambda_{B(0,1)}\Big(\{ y,\, Y_b(y)\geq -g(ye^b)  -\frac{1}{2}\}\Big)]^8}\Big]^{\frac{1}{4}} .
\end{eqnarray}
\qed
%%%%%%%%%%%%%%%%%%%%%%%%%%%%%%%%%%%%%%%%%%%%%%%
\subsection{Proof of Theorem \ref{RESUME}.}

Now we proceed with the proof of Theorem \ref{RESUME}, which is a  slight extension of  \cite[Proposition 5.4 ]{Mad13}. We divide it in two lemmas. The first one is 
 \begin{Lemma}
\label{sansrho}
 Let $b>0$ and $F:\r\times \mathcal{C}(B(0,e^b),\r)\to \r^+$ be a function {\bf $b$ regular}. For any $\epsilon>0$,  there exists $ t',\, T >0$ such that for any $t\geq T$, $\chi(\cdot) \in \mathcal{C}_R(t',\kappa_d \ln t',\ln t) $, $\mathtt{z}\leq (\ln t)^{30}$,
\begin{equation}
\label{eqsansrho}
\Big| \int_{\r_t} e^{-\sqrt{2\d}\chi(x)} \E_{-\chi(x)} \Big[\1_{\{ \sup_{s\in[0,t_b]}B_s\leq 0,\, \sup_{s\in[\frac{t}{2},t_b]}B_s\leq -\mathtt{z} \}}\Big(F(B_{t_b}+\mathtt{z},\mathfrak{G}^{\chi_x}_{t,b} )-F(B_{t_b}+\mathtt{z},\mathfrak{G}_{t,b})\Big)\Big] dx\Big|\leq \epsilon \int_{[0,R]^\d}   \chi(x)  e^{-\sqrt{2\d}\chi(x)} dx.
\end{equation}
 \end{Lemma}

\noindent{\it  Proof of Lemma \ref{sansrho}.} For $t\geq \ln t'+b$,  as $\chi(\cdot) \in \mathcal{C}_R(t',\kappa_\d \ln t',\ln t)$,
\begin{eqnarray*}
||\mathfrak{G}_{t,b}^{\chi_x}-\mathfrak{G}_{t,b}||_\infty\leq  \underset{x\in \r_t,\,y \in B{(0,e^b)}}{\sup}|\chi(x+ye^{-t})-\chi(x)=w_{\chi(\cdot)}(e^{b-t})  \leq e^{-\frac{t-b}{3}}.
\end{eqnarray*}
Recalling (\ref{et3bibi}), the quantity in \eqref{eqsansrho} is smaller than:
\begin{eqnarray*}
 \int_{\r_t} {e^{-\sqrt{2d}\chi(x)}}{t^\frac{3}{2}}    \E_{-\chi(x)}   \Big[ \1_{\{ \sup_{s\in[0,t_b]}B_s\leq 0,\, \sup_{s\in[\frac{t}{2},t_b]}B_s\leq -\mathtt{z} \}} c_9 e^{-\frac{t-b}{12}} h(B_{t_b}+\mathtt{z})F^*(\mathfrak{G}_{t,b})  \Big]dx.
\end{eqnarray*}
with $h$ and $F^*$ two functions associated to the {\bf $b$ regular} function $F$. Finally Lemma \ref{sansrho} follows  provided that we prove:  there exists a constant $c(F)>0 $ (possibly depending  on $F$) such that, for some  $T>0 $ and for any $\alpha\in [1,\ln t],\, \sigma \in [0,t_b] $,
\begin{eqnarray}
\label{sansrhofirst} \E_{-\alpha}\Big[  \1_{\{ \sup_{s\in[0,t_b]}B_s\leq 0,\, \sup_{s\in[\frac{t}{2},t_b]}B_s\leq -\mathtt{z} \}}  h(B_{t_b}+\mathtt{z})F^*(\mathfrak{G}_{t,b,\sigma})  \Big]\leq  {c(F)}\alpha t^{-\frac{3}{2}}.
\end{eqnarray} 

Let us prove (\ref{sansrhofirst}). From (C.23) in \cite{Mad13} we can affirm that   for any $t\geq b>0$ large enough, $\alpha \in [1,\ln t]$, $k,\, j,\, \geq 1$ and $\sigma\in[0,t_b]$, 
 \begin{equation} 
\label{Gcond1}  \E_{-\alpha}\Big[  \1_{\{ \sup_{s\in[0,t_b]}B_s\leq 0,\, \sup_{s\in[\frac{t}{2},t_b]}B_s\leq -\mathtt{z},\, B_{t_b}+\mathtt{z}\in [-(k+1),k] \}}\1_{\{   w_{\mathfrak{G}_{t,b,\sigma}(\cdot)}(\frac{1}{j})\geq \frac{1}{4}\}}   \Big]    \leq  c_{23}(1+L+k)\frac{ \alpha}{t^{\frac{3}{2}}} e^{-c_{24}(b)j}.
\end{equation}  

 According to ({\ref{et2}), there exists $c_1(F)>0$ such that 
\begin{eqnarray*}
h(B_{t_b}+\mathtt{z})\leq \Big\{ \begin{array}{ll} e^{ B_{t_b}+\mathtt{z}} & \text{if  } B_{t_b}+\mathtt{z}\leq -c_1(F), \\
c_1(F) &\text{if   }  B_{t_b}+\mathtt{z} \geq -c_1(F).
\end{array}\Big.
\end{eqnarray*}
 From ({\ref{et1}) we get that
\begin{align*}
\E_{-\alpha}&\Big[  \1_{\{ \sup_{s\in[0,t_b]}B_s\leq 0,\, \sup_{s\in[\frac{t}{2},t_b]}B_s\leq -\mathtt{z} \}}   h(B_{t_b}+\mathtt{z})F^*(\mathfrak{G}_{t,b,\sigma})  \Big]
\\
\leq    &  c_1(F)\sum_{j=1}^\infty  (j+1)^{10}  \Big(\sum_{k=1}^{  c_1(F)} \E_{\P} \Big[ \1_{\{  \sup_{s\in[0,t_b]}B_s\leq 0,\, \sup_{s\in[\frac{t}{2},t_b]}B_s\leq -\mathtt{z},\, B_{t_b}+\mathtt{z}\in [-(k+1),k] \}} \1_{\{   w_{\mathfrak{G}_{t,b,\sigma}(\cdot)}(j^{-1})\geq \frac{1}{4}  \}}   \Big]  
\\
& +\sum_{k= c_1(F)}^{\infty} e^{- k} \E_{\P} \Big[ \1_{\{ \sup_{s\in[0,t_b]}B_s\leq 0,\, \sup_{s\in[\frac{t}{2},t_b]}B_s\leq -\mathtt{z},\, B_{t_b}+\mathtt{z}\in [-(k+1),k] \}} \1_{\{  w_{\mathfrak{G}_{t,b,\sigma}(\cdot)}(j^{-1})\geq \frac{1}{4}  \}}   \Big]\Big).
\end{align*}
Finally by (\ref{Gcond1}) we have for any $t>0$ large enough, $\alpha \in [1,\ln t]$, $k,\, j,\, \geq 1$ and $\sigma\in[0,t_b]$, 
\begin{eqnarray}
\nonumber \E_{-\alpha}\Big[   \1_{\{\sup_{s\in[0,t_b]}B_s\leq 0,\, \sup_{s\in[\frac{t}{2},t_b]}B_s\leq -\mathtt{z}\}} h(B_{t_b}+\mathtt{z})F^*(\mathfrak{G}_{t,b,\sigma})  \Big]&\leq &  \frac{c \alpha}{t^{\frac{3}{2}}}\sum_{j=1}^\infty \frac{(j+1)^{10} }{e^{c_{24}(b)j}}\Big[  c_1(F)^2  + \sum_{k= c(F)}^{\infty} \frac{(1+k)}{e^{ k}} \Big]
\\
\label{ineqPrat}&\leq& \frac{ c(F)\alpha}{t^\frac{3}{2}},
\end{eqnarray}
which ends the proof of Lemma \ref{sansrho}. \qed 

\begin{Remark}As a by product we have also shown the following affirmation.
 Fix $F$ a function {\bf $b$ regular}. For any $\epsilon>0$ there exists $M,\, T>0$  such that for any $t\geq T$, $\alpha\in [1, \ln t]$, $\mathtt{z}\leq (\ln t)^{30}$ and $\sigma\in [0,t_b]$  we have
\begin{equation}
\label{0.3}
\E_{-\alpha}\Big(  \1_{\{\sup_{s\in[0,t_b]}B_s\leq 0,\, \sup_{s\in[\frac{t}{2},t_b]}B_s\leq -\mathtt{z} \}} h(B_{t_b}+\mathtt{z})  F^*\Big(\mathfrak{G}_{t,b,\sigma} \Big)  \1_{\{ w_{\mathfrak{G}_{t,b,\sigma}(\cdot)}(\frac{1}{M})\geq \frac{1}{4}\}} \Big)\leq \frac{\epsilon \alpha}{t^{\frac{3}{2}}},
\end{equation}
and 
\begin{equation}
\label{0.4}
\E_{-\alpha}\Big(  \1_{\{ \sup_{s\in[0,t_b]}B_s\leq 0,\, \sup_{s\in[\frac{t}{2},t_b]}B_s\leq -\mathtt{z} \}} h(B_{t_b}+\mathtt{z})  F^*\Big(\mathfrak{G}_{t,b,\sigma} \Big)  \1_{\{  B_{t_b} +\mathtt{z}\leq -M\}} \Big)\leq \frac{\epsilon \alpha}{t^{\frac{3}{2}}}.
\end{equation} 
Indeed for (\ref{0.3}) as well as for (\ref{0.4}), it suffices to choose $M$ large enough in the inequality  \eqref{ineqPrat}.\end{Remark}

\medskip

For any $M>0$, recall the definition (\ref{defMreload}). 
 
\begin{Lemma}
\label{passageaueps}
(i) Let $b>0$ and $F:\r\times \mathcal{C}(B(0,e^b),\r)\to \r^+$ be a function {\bf $b$ regular}. For any $\delta>0$,  there exist $M,\,\sigma,\, T >0$ such that for any $t\geq T$, $\alpha\in[1,\ln t]$, $\mathtt{z}\leq (\ln t)^{30}$,
\begin{equation}
\label{eqpassageaueps}
\frac{t^{\frac{3}{2}}}{\alpha }\Big|\E_{-\alpha} \Big[\1_{\{ \sup_{s\in[0,t_b]}B_s\leq 0,\, \sup_{s\in[\frac{t}{2},t_b]}B_s\leq -\mathtt{z}\}}\Big(F^{(M)}(B_{t_b}+\mathtt{z},\mathfrak{G}_{t,b,\sigma} )-F(B_{t_b}+\mathtt{z},\mathfrak{G}_{t,b})\Big) \Big] \Big|\leq \delta.
\end{equation}

(ii) Let $b>0$ and $F:\r\times \mathcal{C}(B(0,e^b),\r)\to \r^+$ be a function {\bf $b$ regular}. Fix $M,\sigma > 0$. There exists $C_{M,\sigma}(F)>0$ such that for any $\delta > 0$, there exists $T\geq 0$ such that for any  $t\geq T$, $\alpha \in [1,\ln t]$, $\mathtt{z}\leq (\ln t)^{30}$,
\begin{equation}
\label{eqpassageaueps2}
\Big|\frac{t^{\frac{3}{2}}}{\alpha } \E_{-\alpha}\Big[\1_{\{\sup_{s\in[0,t_b]}B_s\leq 0,\, \sup_{s\in[\frac{t}{2},t_b]}B_s\leq -\mathtt{z}\}}F^{(M)}\big(B_{t_b}+\mathtt{z},\mathfrak{G}_{t,b,\sigma}\big)\Big]-C_{M,\sigma}(F)\Big|\leq \delta.
\end{equation}
\end{Lemma}

\medskip

 \noindent{\it  Proof of  (\ref{eqpassageaueps}).}  Let    $ b,\, \delta>0  $ and $F$ {\bf $b$ regular}. We have to study the expectation under  $\E_{-\alpha} $ of 
 \begin{eqnarray*}
  \1_{\{ \sup_{s\in[0,t_b]}B_s\leq 0,\, \sup_{s\in[\frac{t}{2},t_b]}B_s\leq -\mathtt{z} \}}\Big|F^{ (M)}(B_{t_b}+\mathtt{z},\mathfrak{G}_{t,b,\sigma})-F (B_{t_b} +\mathtt{z},\mathfrak{G}_{t,b})  \Big|.
 \end{eqnarray*}
Thanks to (\ref{0.3}) and (\ref{0.4}) we can choose  $M_1  $ large enough to restrain our study to the expectation of 
\begin{equation}
\label{restrain}
\1_{\{ w_{\mathfrak{G}_{t,b,\sigma}(\cdot)}(\frac{1}{M_1})\leq \frac{1}{4},\,  B_{t_b}\geq -\mathtt{z}-M_1\}}  \1_{\{ \sup_{s\in[0,t_b]}B_s\leq 0,\, \sup_{s\in[\frac{t}{2},t_b]}B_s\leq -\mathtt{z}\}}  \Big| F^{ (M)}(B_{t_b} +\mathtt{z},\mathfrak{G}_{t,b,\sigma})-F (B_{t_b} +\mathtt{z},\mathfrak{G}_{t,b} ) \Big|,
\end{equation}
%Attention restrain est vrai mais demande de bien reinvestir ineqPrat
with $M>M_1,\, t>b  $. On   $ \{ w_{\mathfrak{G}_{t,b,\sigma}(\cdot)}(\frac{1}{M_1})\leq \frac{1}{4},\,  B_{t_b} +\mathtt{z}\geq  -M_1\} $, by (i) of Lemma \ref{controlF}, (\ref{et1}) and (\ref{et2}) 
\begin{equation}
\label{useful}
F ( B_{t_b}+\mathtt{z},\mathfrak{G}_{t,b,\sigma} )\leq h(B_{t_b} +\mathtt{z})F^*(\mathfrak{G}_{t,b,\sigma}) \leq c M_1^{10}:=M .
\end{equation}
 Then (\ref{restrain}) is equal to 
\begin{equation}
\label{restrain2}
\1_{\{ w_{\mathfrak{G}_{t,b,\sigma}(\cdot)}(\frac{1}{M_1})\leq \frac{1}{4},\,  B_{t_b}\geq -\mathtt{z}-M_1\}}  \1_{\{ \sup_{s\in[0,t_b]}B_s\leq 0,\, \sup_{s\in[\frac{t}{2},t_b]}B_s\leq -\mathtt{z}\}}  \Big|F (B_{t_b} +\mathtt{z},\mathfrak{G}_{t,b,\sigma})-F (B_{t_b} +\mathtt{z},\mathfrak{G}_{t,b} ) \Big|\wedge 2M.
\end{equation}
We denote  $ || \Delta \mathfrak{G}_{\sigma}||_{\infty}:= \underset{y\in B(0,e^b)}{\sup} \Big|\mathfrak{G}_{t,b} (y)-\mathfrak{G}_{t,b,\sigma}(y) \Big| $, by Lemma \ref{controlF} (via the property (\ref{et3bibi})) then (\ref{useful}) we deduce that (\ref{restrain2}) is smaller than
\begin{equation}
\1_{\{ B_{t_b}\geq -\mathtt{z}-M_1\}}  \1_{\{\sup_{s\in[0,t_b]}B_s\leq 0,\, \sup_{s\in[\frac{t}{2},t_b]}B_s\leq -\mathtt{z} \}}   2M\Big(  \1_{\{  || \Delta \mathfrak{G}_{\sigma}||_{\infty}\geq \delta^4 \}}   +      \1_{\{  || \Delta \mathfrak{G}_{\sigma}||_{\infty}\leq \delta^4 \}} \delta  \Big).
\end{equation}
Taking the expectation the proof of (\ref{eqpassageaueps}) stems directly from the following two assertions:

-For any $L,\, b,\delta,\,M_1,\,  \epsilon >0$ there exists $ T>0$ such that for any $t\geq T$, $\alpha \in [1,\ln t]$, $\mathtt{z}\leq (\ln t)^{30}$ we have
\begin{equation}
\label{eqcontroldiff01}
 \P_{-\alpha} \Big( \sup_{s\in[0,t_b]}B_s\leq 0,\, \sup_{s\in[\frac{t}{2},t_b]}B_s\leq -\mathtt{z},\,  B_{t_b} +\mathtt{z}\geq  -M_1  \Big)\leq c_{12} \frac{ \alpha}{t^{\frac{3}{2}}}(1+M_1)^2,
\end{equation}

 and 
\\

-For any $L,\, b,\delta,\,M_1,\,  \epsilon >0$ there exists $\sigma,\, T>0$ such that for any $t\geq T$, $\alpha \in [1,\ln t]$, $\mathtt{z}\leq (\ln t)^{30}$ we have
\begin{equation}
\label{eqcontroldiff}
 \P_{-\alpha} \Big( || \Delta \mathfrak{G}_{\sigma}||_{\infty}\geq \delta,\, \sup_{s\in[0,t_b]}B_s\leq 0,\, \sup_{s\in[\frac{t}{2},t_b]}B_s\leq -\mathtt{z},\,  B_{t_b} +\mathtt{z}\geq  -M_1  \Big)\leq c \frac{ \alpha}{t^{\frac{3}{2}}}(1+M_1)^2\exp(- \frac{c_{19}}{2}\delta^2 e^{2\sigma}).
\end{equation} 
The assertion (\ref{eqcontroldiff01}) comes from \cite[equation (B.7)]{Mad13}, whereas  (\ref{eqcontroldiff}) is a consequence of \cite[equation (5.55)]{Mad13}.
By choosing $\delta $ small enough and applying  \eqref{eqcontroldiff01}, then by choosing $\sigma>0$ large enough and applying (\ref{eqcontroldiff}), we obtain \eqref{eqpassageaueps}.\qed
\\

Now we tackle the proof of (\ref{eqpassageaueps2}). Let us introduce some notations from \cite{Mad13}:

- Let $(\mathtt{R}_{s})_{s \geq 0}$ be a three dimensional Bessel process starting from $0$.

- Let $(B_s)_{s\geq 0}$ be real Brownian motion and for any $\sigma>0 $ we denote $(B_s^{(\sigma)})_{s\geq 0}:= (B_{s+\sigma}-B_\sigma)_{s\geq 0} $.

- Let  $g$, $h$ be two processes, for any $ t_0\in \r^+$ the process $\mathtt{X}_\cdot(t_0,g,h)$ is defined by  
\begin{eqnarray}
\label{defXX} \mathtt{X}_s(t_0,g,h)=\Big\{\begin{array}{ll} g_s ,&\qquad \text{if  } s\leq t_0,
\\
g_{t_0}+ h_{t-t_0},&\qquad \text{if  } s\geq t_0.
\end{array}
 \Big. 
\end{eqnarray}

- Let $\sigma>0 $ for any process $(g_s)_{s\leq \sigma}$ we set 
\begin{equation}
\label{defarrw} (\overset{\leftarrow \sigma}{g_s})_{s\leq \sigma}:= (g_{\sigma-s}-g_\sigma)_{s\leq \sigma}.
\end{equation}

- We set  $\mathfrak{H}_{m,\sigma}$ the set of continuous function $F:\r\times \mathcal{C}([0,\sigma],\r)\to \r^+$ with $\underset{u\in \r,\, g\in \mathcal{C}([0,\sigma ],\r))}{\sup} F(u,g)\leq m$. For $ g\in \mathcal{C}^1(\r^\d,\r) $ we denote by $ \nabla_{y}(g) $ the gradient of $g$ at $y\in \r^\d $. Finally, we denote by $ \langle\cdot,\cdot\rangle $  the inner product in $\r^\d $.
 \\

We will derive \eqref{eqpassageaueps2}  from the following Proposition, which is proven in \cite{Mad13} (cf Proposition 5.5 in \cite{Mad13}).

 \begin{Proposition}[\cite{Mad13}]
\label{renouv}
Let $B$ be a Brownian motion and let $ \mathtt{R}$ be a three dimensional Bessel process starting from $0$ independent of $B$. Let $m,\, \sigma \geq 0$ be two constants. For any $\epsilon>0$ there exists $T(m,\sigma ,\epsilon)\geq 0$ such that for any $t\geq T$, $1\leq \alpha,\,  \mathtt{z}\leq (\ln t)^{30} $ and $F\in \mathfrak{H}_{m,\sigma}$
\begin{eqnarray}
\label{eqrenouv}
 \Big| \frac{t^\frac{3}{2}}{\alpha}\E_\alpha\Big[\1_{\{ \inf_{s\in[0,t]}B_s \geq 0,\, \inf_{s\in[\frac{t}{2},t]}B_s\geq \mathtt{z},\, B_t-\mathtt{z} \leq m\}}F\Big( B_t-\mathtt{z}, (B^{(t-\sigma)}_{s})_{s\leq \sigma}\Big)\Big]- \qquad\qquad\qquad 
 \\
\nonumber    C  \int_{0}^m  \int_{0}^u \E\Big[     F(u, (\overset{\leftarrow \sigma}{\mathtt{X}_s}(T_{-\gamma}, B, \mathtt{R}))_{s\leq \sigma} )\Big]d\gamma du    \Big|\leq \epsilon    ,
\end{eqnarray}
where $C:=  \frac{4}{\sqrt{\pi}} \E[e^{- {\mathtt{R}^2_{1}}/2 }]=\sqrt{\frac{2}{\pi}}$ and $T_\gamma:=\inf\{ s\geq 0,\, B_s=\gamma\}$, $\gamma \in \r $.
\end{Proposition}
%\sup_{s\in[0,t_b]}B_s\leq 0,\, \sup_{s\in[\frac{t}{2},t_b]}B_s\leq -\mathtt{z}

\medskip 

\noindent{\it  Proof of  \eqref{eqpassageaueps2}.} Fix $b,\, M,\, \sigma>0$ and $F$ a function {\bf $b$ regular}. Let us explicit the expectation in (\ref{eqpassageaueps2}). As $(B_s)_{s\geq 0}\overset{law}{=}(-B_s)_{s\geq 0} $ we have,
\begin{eqnarray*}
\E(\ref{eqpassageaueps2})&:=&\E_{-\alpha}\Big[\1_{\{ \sup_{s\in[0,t_b]}B_s\leq 0,\, \sup_{s\in[\frac{t}{2},t_b]}B_s\leq -\mathtt{z}\}}F^{(M)}\Big(B_{t_b} +\mathtt{z},\mathfrak{G}_{t,b,\sigma}\Big)\Big]
\\
&=& \E_{\alpha}\Big[\1_{\{ \inf_{s\in[0,t]}B_s \geq 0,\, \inf_{s\in[\frac{t}{2},t]}B_s\geq \mathtt{z},\, B_{t_b}-\mathtt{z}\leq M \}} F^{(M)}\Big(-[B_{t_b}-\mathtt{z}], 
\\
&&\qquad\qquad\qquad \qquad\qquad\qquad  y\mapsto \int_{t_b-\sigma}^{t_b}(1-k(e^{s-t}y))dB_s-\zeta_{t_b}(ye^{-t})+Z^0_{t_b}(ye^{-t})    \Big)\Big].
\end{eqnarray*}
Moreover by integration by parts, the second argument of the function in $F^{(M)}$ can be rewritten as:
\begin{eqnarray*}
y\mapsto \big(1-k(e^{-b}y)\big)\Big[B_{t_b}-B_{t_b-\sigma}  \Big] +   \int_{t_b-\sigma}^{t_b} [B_s-B_{t_b-\sigma}]\langle \nabla_{ye^{s-t}} k\cdot ye^{s-t}\rangle ds    -\zeta_{t_b}(ye^{-t})+Z^0_{t_b}(ye^{-t}),
\end{eqnarray*}
and we recall that the processes $B $ and $Z$ are independent. So $ \E(\ref{eqpassageaueps2}) $ is equal to 
\begin{eqnarray*}
\E_{\alpha}\Big[\1_{\{  \inf_{s\in[0,t]}B_s \geq 0,\, \inf_{s\in[\frac{t}{2},t]}B_s\geq \mathtt{z},\, B_{t_b}- \mathtt{z}\leq M\}}\Phi_{t_b}(B_{t_b}-\mathtt{z},\, (B_s^{(t-\sigma)})_{s\leq \sigma})\Big],
\end{eqnarray*}
with $ \Phi_{t_b}: \r\times \mathcal{C}([0,\sigma],\r) \to \r^+ $ , a continuous function, bounded by $M$ and defined by
\begin{eqnarray*}
(u,h) \mapsto \E\Big[ F^{(M)}(-u,y\mapsto   \big(1-k(e^{t_b-t}y)\big)[h_{\sigma}-h_{0}]+\int_{t_b-\sigma}^{t_b} [h_{s-(t_b-\sigma)}-h_{0}]  \langle \nabla_{ye^{s-t}}k. ye^{s-t}\rangle  ds \Big.\qquad
\\
\Big. - \zeta_t(ye^{-t})+ Z^0_{t_b}(ye^{-t})\Big].
\end{eqnarray*}
 Now we can apply Proposition \ref{renouv}, with  $t \leftrightarrow t_b>0$,  $\alpha    \leftrightarrow  \alpha ,\,  \mathtt{z} \leftrightarrow  \mathtt{z}\leq (\ln t)^{30},\, \sigma \leftrightarrow  \sigma,\, m  \leftrightarrow M$ and  $F\leftrightarrow \Phi_{t_b}  $. Then for any $\epsilon>0$ there exists $T\geq 0$ such that for any $t\geq T$, $1\leq \alpha\leq (\ln t)^{30}$
\begin{eqnarray}
\label{1(ii)}
\Big[|\frac{t^\frac{3}{2}}{\alpha}\E(\ref{eqpassageaueps2})-C  \int_{0}^m  \int_{0}^u \E\Big[(     \Phi_{t_b}(u, (\overset{\leftarrow \sigma}{\mathtt{X}_s}(T_{-\gamma}, B, \mathtt{R}))_{s\leq \sigma} )\Big)d\gamma du    \Big|\leq \epsilon   .
\end{eqnarray}
Moreover, we observe that for any  $u>0,\gamma\leq u $,  
\begin{eqnarray*} 
&& \E\Big[(\Phi_{t_b}(u, \overset{\leftarrow \sigma}{\mathtt{X}}_s(T_{-\gamma}, B, R))_{s\leq \sigma} )\Big]= \E\Big[ F^{(M)}\Big(-u,  y\mapsto Z_{t_b}(ye^{-t})-\zeta_{t_b}(ye^{-t})
 \\
&& \qquad\qquad\qquad\qquad \qquad\qquad\qquad   -\int_{0}^{T_{-\gamma} \wedge \sigma }\big(1-k(e^{-s}ye^{-b})\big) dB_s- \int_{T_{-\gamma}\wedge \sigma}^\sigma\big(1-k(e^{-s}ye^{-b})\big) d \mathtt{R}_{s-T_{-\gamma}}\Big)\Big].
\end{eqnarray*}
Finally as $ (Z_{t_b}(ye^{-t})-\zeta_{t_b}(ye^{-t}))_{y\in B(0,e^{b})} $ is independent of $(B,R)$ and converges in law, when $t$ goes to infinity, (see (2.6) in \cite{Mad13}) to $ (Z(ye^{-b})-\zeta(ye^{-b}))_{y\in B(0,e^{b})} $, by combining with (\ref{1(ii)}) we deduce that: for any $\epsilon>0$ there exists $T\geq 0$ such that for any $t\geq T$, $1\leq \alpha\leq (\ln t)^{30}$
\begin{eqnarray}
\Big|\frac{1}{\alpha t^{\frac{3}{2}}} \E_{-\alpha}\Big[ \1_{\{\sup_{s\in[0,t_b]}B_s\leq 0,\, \sup_{s\in[\frac{t}{2},t_b]}B_s\leq -\mathtt{z}\}} F^{(M)}\big(B_{t_b}+\mathtt{z},\mathfrak{G}_{t,b,\sigma}\big)\Big]-C_{M,\sigma}(F)\Big|\leq \epsilon,
\end{eqnarray}
with 
\begin{eqnarray}
\nonumber C_{M,\sigma}(F):= C\int_{0}^M\int_0^u\E\Big[  F^{(M)}\Big(-u, y\mapsto Z (ye^{-b})-\zeta(ye^{-b}) \qquad\qquad\qquad\qquad\qquad
\\
  -\int_{0}^{T_{-\gamma} \wedge \sigma }\big(1-k(e^{-s}ye^{-b})\big) dB_s - \int_{T_{-\gamma}\wedge \sigma}^\sigma \big(1-k(e^{-s}ye^{-b})\big) d \mathtt{R}_{s-T_{-\gamma}}\Big)\Big]d\gamma dU.
\end{eqnarray}
This completes the proof of (\ref{eqpassageaueps2})  \hfill $\Box $

Theorem \ref{RESUME} is a combination of (\ref{eqsansrho}), (\ref{eqpassageaueps}) and (\ref{eqpassageaueps2}). \hfill$\Box$

%%%%%%%%%%%%%%%%%%%%%%%%%%%%%%%%%%%%%%%%%%%%%%%%%%%%%%%%%%%%%%%%%%%%%%%%%%%%%%%%%%%%%%%%%%%%%%%%%%%%%%%%%%%%%%%%%%%%
\section{Proofs for  two dimensional Free Fields}\label{proofGFF}
%%%%%%%%%%%%%%%%%%%%%%%%%%%%%%%%%%%%%%%%%%%%%%%%%%%%%%%%%%%%%%%%%%%%%%%%%%%%%%%%%%%%%%%%%%

%We can for instance consider:\\
%$\bullet$ a discrete MFF on a square grid of $\z^2$ with mesh size going to $0$ or even discrete MFF on planar isoradial graphs (with rhombi half-angles uniformly bounded away from $0$ and $\frac{\pi}{2}$, see \cite{chelkak} for isoradial graphs and \cite[section 4.3]{review} to see the relation with our framework).\\
%$\bullet$ the convolution of a full MFF $X$ with a mollifying sequence $(\chi_n)_n$ (see \cite{vincent}), i.e. $(X\star \chi_n)_n$ or circle averages of the MFF as described in \cite{She07}.\\
%$\bullet$ the standard white noise decomposition of the MFF with covariance kernel $\int_{1/n}^{+\infty}p(t,x,y)\,dt$ where $p(t,x,y)$ stands for the transition densities of the planar Brownian motion.\\
%$\bullet$ and many others...

\subsection{Proof of Theorem \ref{mainMFF}.}  
%%%%%%%%%%%%%%%%%%%%%%%%%%%
Before proceeding with the proof let us make a few observations. First we stress that the kernel $k_m$ satisfies:
 \begin{description}
\item[B.1] $k_m$ is nonnegative, continuous and $k(0)=1$.
\item[B.2] $k_m$ is Lipschitz at $0$, i.e. $|k_m(0)-k_m(x)|\leq C|x|$ for all $x\in\r^2$
\item[B.3] $k_m$ satisfies the integrability condition $\sup_{|e|=1}\int_1^{\infty}\frac{k_m(ue)}{u}\,du<+\infty$.
\end{description}
We stick to the notations of Section \ref{setup} so that we set for $t\geq 0$ and $x\in\r^{\d}$
 \begin{equation}\label{kernelMFF}
G_{m,t}(x) =  \int_1^{e^t}  \frac{k_m(xu)}{u}du.
 \end{equation}
In \cite{Rnew12}, it is proved that

 \begin{theorem}\label{seneta}
We set $M^\gamma_t(dx)=e^{\gamma X_t(x)-\frac{\gamma^2}{2}\E[X_t(x)^2]}\,dx$. The family  $(\sqrt{ t}M_t^{\sqrt{2d}})_t$ weakly converges in probability as $t\to \infty$ towards a non trivial limit, which turns out to be the same, up to a multiplicative constant, as the limit of the derivative martingale. More precisely, we have
$$\sqrt{ t}M_t^{\sqrt{2d}}(dx)\to \sqrt{\frac{2}{\pi}} M'(dx),\quad\text{in probability as }t\to\infty.$$
 \end{theorem}

Now we begin the proof and we first theat the case when the cut-off  family of the MFF is $(X_t)_t$. We will see thereafter that the general case (i.e. any other cut-off uniformly close to $(G_{m,t})_t$) is a straightforward consequence. The problem to face is that the covariance kernel of the family $(X_t)_t$ does not possess a compact support so that Theorem \ref{main} does not apply as it is. We split the proof into two levels: a main level along which we explain the main steps of the proof relying on a few lemmas and a second level in which we prove these auxiliary lemmas.

\vspace{2mm}
{\bf Main level:} Let us consider any nonnegative smooth function $\rho:\r^2\to\r_+$ such that: $\int_{\r^2}\rho^2(y)\,dy=1$, $\rho$ is isotrop and has compact support. We set $$\varphi(x)=\int_{\r^2}\rho(y+x)\rho(y)\,dy.$$
Under the above assumptions on $\rho$, it is plain to see that $\varphi$ is nonnegative, smooth, positive definite, with compact support, $\varphi(0)=1$ and isotrop. For each $\epsilon>0$, let us define the function
$$\forall x\in\r^2,\quad \varphi_\epsilon(x)=\varphi(\epsilon x).$$
It is  straightforward to check that the family $(\varphi_\epsilon)_\epsilon$ uniformly converges towards $1$ over the compact subsets of $\r^2$ as $\epsilon\to 0$.
For $\epsilon<0$, we further define
$$k_\epsilon(x)=k_m( x)\varphi_\epsilon( x),\quad K^\epsilon_t(x)=\int_1^{e^t}\frac{k_\epsilon(ux)}{u}\,du.$$ 
Observe that $k_\epsilon$ satisfies Assumption (A) (it is positive definite because it is the Fourier transform of the convolution of the spectrum of $k_m$ and that of $\varphi_\epsilon$). For each $\epsilon>0$, we follow Section \ref{setup} to introduce all the objects related to the kernel $k_\epsilon$ and add an extra superscript $\epsilon$ in the notation to indicate the relation to $k_\epsilon$ (i.e. we introduce $(X^\epsilon_t(x))_{t,x}$, $M^{',\epsilon}_t$, $M^{\gamma,\epsilon}_t$, $M^{',\epsilon}$).

Now we claim:
\begin{Lemma}\label{lem:diffk} 
For each $\delta>0$, there exists $\epsilon_0>0$ such that for all $0<\epsilon<\epsilon_0$ and all $x\in \r^2$ and all $t\geq 0$:
\begin{equation}
K^\epsilon_t(x)-\delta\leq G_{m,t}(x)\leq K^\epsilon_t(x)+\delta.
\end{equation}
\end{Lemma}

This lemma  will help us to see the family $(X^\epsilon_t)_{t}$ as a rather good approximation of the family $(X_t)_t$ as $\epsilon\to 0$. Because $k_\epsilon$ satisfies Assumption (A), Theorem \ref{main} holds for the family $(X^\epsilon_t)_{t}$ for any $\gamma>2$. The conclusion of this theorem involves some constant $C _\epsilon(\gamma)$, which may depend on $\epsilon$. Fortunately, we claim:

\begin{Lemma}\label{univC}
For each fixed $\gamma>2$, the family $(C _\epsilon(\gamma))_\epsilon$ converges as $\epsilon\to 0$ towards some constant denoted by $C (\gamma)$.
\end{Lemma}

Then we use Lemmas \ref{lem:diffk} and \ref{univC} to prove
\begin{Lemma}\label{uniflimitMFF}
For any $\gamma>2$ and for any continuous nonnegative function $f$ with compact support, we have
\begin{equation*}
\underset{t\to\infty}{\lim}\E\Big[ \exp(- t^{\frac{3\gamma}{4}}e^{t(\frac{\gamma}{\sqrt{2}}-\sqrt{2})^2}M_t^\gamma(f))\Big] =\underset{\epsilon\to0}{\lim} \,\underset{t\to\infty}{\lim}\E\Big[\exp(-C (\gamma) \int_{\r^2}f(x)^\frac{2}{\gamma}M^{',\epsilon}(dx))\Big]  .
\end{equation*}
\end{Lemma}
It thus remains to compute the above double limit:
\begin{Lemma}\label{invert}
For any $\gamma>2$ and for any continuous nonnegative function $f$ with compact support, we have
\begin{equation*}
\underset{\epsilon\to0}{\lim}  \E\Big[\exp(-C (\gamma) \int_{\r^2}f(x)^\frac{2}{\gamma}M^{',\epsilon}(dx))\Big]=\E\Big[\exp(-C (\gamma) \int_{\r^2}f(x)^\frac{2}{\gamma}M^{'}(dx))\Big].
\end{equation*}
\end{Lemma}
We are now done with the proof of Theorem \ref{mainMFF}. It just remains to prove  the four above  lemmas.\qed

\vspace{2mm}
{\bf Proofs of auxiliary lemmas.} Most of the forthcoming proofs will heavily rely on Kahane's convexity inequalities (KCI for short) so that the reader is referred to  Kahane's original paper \cite{cf:Kah} (or \cite{review} for an english statement).

\noindent {\it Proof of Lemma \ref{lem:diffk}.} Let us fix $\delta>0$.  In what follows and when considering an isotrop function $f$, we will identify, with a slight abuse of notations, the function $f:\r^2\to \r$ with the function $f:\r_+\to\r$ through the relation $f(y)=f(|y|)$ for $y\in\r^2$. Observe that:
\begin{align*}
|K^\epsilon_t(x)- G_{m,t}(x)|= &\big|\int_1^{e^t}\frac{k_m(u|x|)-k_\epsilon(u|x|)}{u}\,du\big| 
\leq \int_1^{\infty}\frac{|k_m(u|x|)-k_\epsilon(u|x|)|}{u}\,du 
\leq   \int_{0}^{\infty}\frac{|k_m(v)-k_\epsilon(v)|}{v}\,dv.%\int_{\min(1,|x|)}^1\frac{|k(v)-k_\epsilon(v)|}{v}\,dv+\int_1^{\infty}\frac{|k(v)-k_\epsilon(v)|}{v}\,dv
\end{align*}
We prove now that we can get the   above quantities arbitrarily close to $0$. We fix $R>1$ such that $\int_R^{\infty}\frac{|k_m(v)|}{v}\,dv\leq \delta/4$. Since $\varphi_\epsilon(y)\leq \varphi_\epsilon(0)=1$ (by positive definiteness), we also have 
$$\int_R^{\infty}\frac{|k_\epsilon(v)|}{v}\,dv|\leq  \int_R^{\infty}\frac{|k_m(v)|}{v}\,dv\leq \delta/4.$$
On $[1,R]$, we use the fact that the family $(\varphi_\epsilon)_\epsilon$ uniformly converges towards $1$ over compact sets to deduce that for some $\epsilon_0$ and all $\epsilon<\epsilon_0$, we have
$$\int_1^R \frac{|k(v)-k_\epsilon(v)|}{v}\,dv\leq \delta/4.$$
It remains to treat the interval $[0,1]$. Since $\varphi$ is smooth, it is locally Lipschitz at $0$, meaning that we can find a constant $C$ such that $|1-\varphi(x)|\leq C|x|$ for all $x$ belonging to some ball centered at $0$, say $B(0,1)$. Furthermore $|k_m(v)|\leq k_m(0)=1$. We deduce:
\begin{align*}
\int_{0}^{1}\frac{|k_m(v)-k_\epsilon(v)|}{v}\,dv&=\int_{0}^{1}\frac{|k_m(v)||1-\varphi(\epsilon v)|}{v}\,dv\\
&\leq C\epsilon.
\end{align*}
For $\epsilon$ small enough, this quantity can be made less than $\delta/4$.\qed

\vspace{2mm}
 
\noindent {\it Proof of Lemma \ref{invert}.} Let us fix $\delta>0$. From Lemma \ref{lem:diffk}, we have for $\epsilon$ small enough and all $t$ and $x$ the inequality $K^\epsilon_t(x)\leq G_{m,t}(x)\leq K^\epsilon_t(x)+\delta$ (the $\delta$ has been omitted from the left-hand side because this inequality is obvious and results from $\varphi\leq 1$). By applying KCI to the convex function $x\mapsto e^{-x}$ combined with the inequalities  $K^\epsilon_t(x)\leq G_{m,t}(x)\leq K^\epsilon_t(x)+\delta$, we get for all $\theta>0$ and some standard Gaussian random variable $\mathcal{N}$ independent from everything:
\begin{align*}
\E\Big[\exp(-\theta\sqrt{t}\int_{\r^2} f(x)^\frac{2}{\gamma}M^{2,\epsilon}_t(dx))\Big]\leq &\E\Big[\exp(-\theta\sqrt{t} \int_{\r^2}f(x)^\frac{2}{\gamma}M^{2}_t(dx))\Big]\\
\E\Big[\exp(-\theta\sqrt{t} \int_{\r^2}f(x)^\frac{2}{\gamma}M^{2}_t(dx))\Big]\leq  &\E\Big[\exp(-e^{\sqrt{\delta}\mathcal{N}-\delta/2}\theta\sqrt{t}\int_{\r^2} f(x)^\frac{2}{\gamma}M^{2}_t(dx))\Big].
\end{align*}
By taking the limit as $t\to\infty$ and by using Theorem \ref{seneta}, we obtain for all $\theta\geq 0$:
\begin{align*}
\E\Big[\exp(-\theta M^{',\epsilon}(f^\frac{2}{\gamma}))\Big]\leq &\E\Big[\exp(-\theta M'(f^\frac{2}{\gamma}))\Big]\leq 
\E\Big[\exp(-\theta e^{\sqrt{\delta}\mathcal{N}-\delta/2} M^{',\epsilon}(f^\frac{2}{\gamma}))\Big] 
\end{align*}
It is then straightforward to deduce that:
$$\lim_{\epsilon\to 0}\E\left[\exp(-\theta M^{',\epsilon}(f^\frac{2}{\gamma}))\right]=\E\left[\exp(-\theta M'(f^\frac{2}{\gamma}))\right],$$
which concludes the proof of the Lemma.\qed

\vspace{2mm}
\noindent {\it Proof of Lemmas \ref{univC} and \ref{uniflimitMFF}.} First recall that the family $\big(t^{\frac{3\gamma}{4}}e^{t(\frac{\gamma}{\sqrt{2}}-\sqrt{2})^2}M_t^{\gamma,\epsilon}(f)\big)_t$ is tight and every possible converging limit (in law) is non trivial \cite[section 4.3]{Rnew7} provided that $f$ is non trivial.

Then, for any $\epsilon>0$, we have from Theorem \ref{main}
\begin{equation*}
\underset{t\to\infty}{\lim}\E\left[ \exp(- t^{\frac{3\gamma}{4}}e^{t(\frac{\gamma}{\sqrt{2}}-\sqrt{2})^2}M_t^{\gamma,\epsilon}(f))\right] =\E\left[\exp(-C _\epsilon(\gamma) \int_{\r^2}f(x)^\frac{2}{\gamma}M^{',\epsilon}(dx))\right]  .
\end{equation*}
Furthermore, for each $\delta>0$ and $\epsilon$ small enough, we have at our disposal the inequality $K^\epsilon_t(x)\leq G_{m,t}(x)\leq K^\epsilon_t(x)+\delta$ and a convex function $x\mapsto e^{-x}$. So, denoting by $\mathcal{N}$ a standard Gaussian random variable,  we can apply KCI to get for all $\delta>0$,  $\epsilon$ large enough and all $\theta\geq 0$:
\begin{align}
\label{KCI1}\E\Big[ \exp(- \theta t^{\frac{3\gamma}{4}}e^{t(\frac{\gamma}{\sqrt{2}}-\sqrt{2})^2}M_t^{\gamma,\epsilon}(f))\Big]\leq & \E\Big[\exp(- \theta t^{\frac{3\gamma}{4}}e^{t(\frac{\gamma}{\sqrt{2}}-\sqrt{2})^2}M_t^{\gamma}(f))\Big]\\
\label{KCI2}\E\Big[ \exp(-\theta  t^{\frac{3\gamma}{4}}e^{t(\frac{\gamma}{\sqrt{2}}-\sqrt{2})^2}M_t^{\gamma}(f))\Big] \leq &\E\Big[ \exp(- \theta e^{\sqrt{\delta}\mathcal{N}-\delta/2}t^{\frac{3\gamma}{4}}e^{t(\frac{\gamma}{\sqrt{2}}-\sqrt{2})^2}M_t^{\gamma,\epsilon}(f))\Big].
\end{align}
Consider a possible limit $Z$ of some subsequence of the family $\big(t^{\frac{3\gamma}{4}}e^{t(\frac{\gamma}{\sqrt{2}}-\sqrt{2})^2}M_t^{\gamma }(f)\big)_t$. By taking the limit as $t\to \infty$ along the proper subsequence in \eqref{KCI1}+\eqref{KCI2}, we get for all $\theta\geq 0$:
\begin{align}
\label{KCIlim1} \E\Big[\exp(-\theta C _\epsilon(\gamma) \int_{\r^2}f(x)^\frac{2}{\gamma}M^{',\epsilon}(dx))\Big] \leq \E\Big[ \exp(- \theta Z)\Big]
\end{align} and for $\eta>0$:
\begin{align}
 \E\Big[ \exp(- \theta Z)\Big]\leq &\E\Big[\exp(-\theta e^{\sqrt{\delta}\mathcal{N}-\delta/2} C _\epsilon(\gamma) \int_{\r^2}f(x)^\frac{2}{\gamma}M^{',\epsilon}(dx))\Big]\nonumber\\
 \label{KCIlim2} \leq & \E\Big[\exp(-\theta (1-\eta)C _\epsilon(\gamma) \int_{\r^2}f(x)^\frac{2}{\gamma}M^{',\epsilon}(dx))\Big]+\P\big(e^{\sqrt{\delta}\mathcal{N}-\delta/2}\leq 1-\eta\big).
\end{align}
By taking the $\limsup_{\epsilon\to 0}$ and $\liminf_{\epsilon\to 0}$ and then $\lim_{\delta\to 0}$, we deduce that for all $\theta\geq 0$ and $\eta>0$:
\begin{align}
\limsup_{\epsilon\to 0} \E\Big[\exp(-\theta C _\epsilon(\gamma) \int_{\r^2}f(x)^\frac{2}{\gamma}M^{',\epsilon}(dx))\Big] \leq &\E\Big[ \exp(- \theta Z)\Big] \\
\liminf_{\epsilon\to 0}   \E\Big[\exp(-\theta C _\epsilon(\gamma) \int_{\r^2}f(x)^\frac{2}{\gamma}M^{',\epsilon}(dx))\Big]  \geq &  \E\Big[ \exp(- \theta (1-\eta)^{-1}Z)\Big].
\end{align}
Now we can take the limit as $\eta\to 0$ and get:
\begin{equation}
\lim_{\epsilon\to 0} \E\Big[\exp(-\theta C _\epsilon(\gamma) \int_{\r^2}f(x)^\frac{2}{\gamma}M^{',\epsilon}(dx))\Big] =\E\Big[ \exp(- \theta Z)\Big].
\end{equation}
Therefore the family $\big(C _\epsilon(\gamma) \int_{\r^2}f(x)^\frac{2}{\gamma}M^{'}(dx)\big)_\epsilon$  converges in law towards $Z$. As a by-product, this shows that the law of $Z$ does not depend on the chosen subsequence along which the family $\big(t^{\frac{3\gamma}{4}}e^{t(\frac{\gamma}{\sqrt{2}}-\sqrt{2})^2}M_t^{\gamma }(f)\big)_t$ converges in law. Thus the whole family converges in law towards a non trivial random variable $Z$. Furthermore, Lemma \ref{invert} shows that the family $\big(\int_{\r^2}f(x)^\frac{2}{\gamma}M^{',\epsilon}(dx)\big)_\epsilon$ converges in law as $\epsilon\to 0$ towards $\int_{\r^2}f(x)^\frac{2}{\gamma}M^{'}(dx)$ which is almost surely strictly positive because $f$ is not trivial. This comes the fact that $M'$ has full support \cite{Rnew7}. It is then straightforward to deduce that the family  $\big(C _\epsilon(\gamma)\big)_\epsilon$ converges as $\epsilon\to 0$.\qed

\vspace{2mm}
{\bf General case:} Now,  we consider a general  cut-off family $(X_n)_n$ of the MFF uniformly close to $(G_{m,t})_t$. By assumption, this family satisfies Lemma \ref{lem:diffk} with $K_n$ instead of $K^\epsilon_t$ and $G_{m,t_n}$ instead of $G_{m,t}$. We can then control the kernel $K_n$ in terms of $G_{m,t_n}$. Furthermore we now that the freezing theorem holds for the family  $(G_{m,t_n})_n$ with some fixed constant $C (\gamma)$: this was the difficult part that we have handled above. Now we can use the same strategy of using KCI to transfer the freezing theorem to the family $(X_n)_n$. Details are obvious and thus left to the reader.\qed

\subsection{Proof of Theorem \ref{mainGFF}.}  
%%%%%%%%%%%%%%%%%%%%%%%%%%%
In what follows, $(X_t)_t$ is the family defined by \eqref{defxt} and  and $$M^\gamma_t(dx)=e^{\gamma X_t-\frac{\gamma^2}{2}t}\,dx,\quad \quad M'(dx)=\lim_{t\to\infty}(2\E[X_t(x)^2]-X_t(x))e^{2X_t(x)-2\E[X_t(x)^2]}\,dx.$$ For $t_0>0$, we will also consider $$M'_{t_0,\infty}(dx)=\lim_{t\to\infty}(2\E[(X_t-X_{t_0})(x)^2]-X_t(x)+X_{t_0}(x))e^{2(X_t-X_{t_0})(x)-2\E[(X_t-X_{t_0})(x)^2]}\,dx.$$

For each $t_0>0$, we consider the MFF like field 
$$X^{MFF}_{t_0,t}(x)=\sqrt{\pi}\int_{\r^2\times [e^{-2t},e^{-2t_0}[}p(\frac{s}{2},x,y)W(dy,ds)$$ with covariance kernel
$$G_{t_0,t}(x,y)= \int_{e^{-2t}}^{e^{-2t_0}}p(s,x,y)\,ds.$$
We further introduce the corresponding measures for $\gamma>2$ 
 $$M^{\gamma,MFF}_{t_0,t}(dx)=e^{\gamma X_{t_0,t}^{MFF}(x)-\frac{\gamma^2}{2}\E[X^{MFF}_{t_0,t}(x)^2]}\,dx$$
 and the derivative multiplicative chaos
 $$M^{',MFF}_{t_0,\infty}(dx)=\lim_{t\to\infty}(2t-2t_0 - X_{t_0,t}^{MFF}(x))e^{2X^{MFF}_{t_0,t}(x)-2\E[(X^{MFF}_{t_0,t}(x))^2]}\,dx.$$
The strategy that we followed for the MFF works for this process as well and Theorem \ref{mainMFF} works for some constant
 \begin{equation}
\label{eqGFF1}\underset{t\to\infty}{\lim} \E\Big[ \exp(- t^{\frac{3\gamma}{4}}e^{(t-t_0)(\frac{\gamma}{\sqrt{2}}-\sqrt{2})^2}M_{t_0,t}^{\gamma,MFF}(f))\Big]=\E\Big[\exp(-C _{t_0}(\gamma) \int_{\r^2}f(x)^\frac{2}{\gamma}M^{',MFF}_{t_0,\infty}(dx))\Big],
\end{equation}
It is clear that the constant  $C _{t_0}(\gamma)$ is likely to depend depend on $t_0$. Indeed,  observe that the covariance kernel of $X^{MMF}_{t_0,t}$ is the same as $X^{MMF}_{t}$ up to a multiplicative change of spatial coordinates so that this variation in the covariance structure should affect $C _{t_0}(\gamma)$. Actually we can even explicitly calculate this dependence
\begin{Lemma}\label{constant} 
Let us set $C (\gamma)=C _{0}(\gamma)$. The constant  $C _{t_0}(\gamma)$ satisfies:
\begin{equation}
\label{moveto}C _{t_0}(\gamma)e^{-2t_0+\frac{4}{\gamma}t_0}=C (\gamma),\quad \forall t_0\geq 0.
\end{equation}
\end{Lemma}

\noindent {\it Proof.}  It suffices to apply  \eqref{eqGFF1} at two different scales $t_0$ and $t_0+s$. Then in the relation corresponding to $t_0+s$, we replace the function $f$ by $e^{s(\frac{\gamma}{\sqrt{2}}-\sqrt{2})^2}e^{\gamma X_{t_0,t_0+s}^{MFF}(x)-\frac{\gamma^2}{2}\E[X_{t_0,t_0+s}^{MFF}(x)^2]}f(x)$, which remains a compactly supported continuous function. It is random but independent of the measure $M_{t_0,t}^{\gamma,MFF}(dx)$. By identification of both limits, we get the relation $C _{t_0+s}(\gamma)e^{-2s+\frac{4}{\gamma}s}=C _{t_0}(\gamma)$. \qed

\vspace{1mm}
Equipped with this relation, we can now try to apply the freezing theorem to a process that we call {\it switch process}.  Equipped with this relation, we can now try to apply the freezing theorem to a process that we call switch process. Basically the switch process is a Gaussian interpolation between the MFF and the GFF. We will plug this switch process in \eqref{eqGFF1} in order transfer by interpolation the property \eqref{eqGFF1} to the GFF. For $t_0\leq t$,   the  switch process is defined by
$$S_{t_0,t}(x)=X_{t_0}(x)+X^{MFF}_{t_0,t}(x) $$ and we also consider the associated measure
$$M^{\gamma,switch}_{t_0,t}(dx)=e^{\gamma S_{t_0,t}(x)-\frac{\gamma^2}{2}t}\,dx.$$

To evaluate to which extent the switch process is a good interpolation between the MFF and the GFF, we need to evaluate how  the covariance kernel of the switch process evolves with $t_0$. To this purpose, we set 
$$\forall x,y\in D,\quad G_{D,t_0,t}(x,y)=G_{D,t}(x,y)-G_{D,t_0}(x,y).$$

Consider a domain $D'\subset D$ such that ${\rm dist}(D',D^c)>0$. We have 
\begin{equation}\label{unifgreen}
\lim_{t_0\leq t\to\infty} \sup_{x,y\in D'}|G_{D,t_0,t}(x,y)-G_{t_0,t}(x,y)|=0.
\end{equation}
This comes from the following lemma, the proof of which is postponed to the end of this subsection
\begin{Lemma}\label{cvunif}
For all subset $D'$ of $D$ such that ${\rm dist}(D',D)>0$, the following convergence holds uniformly on $D'\times D'$:
$$ \lim_{t\to 0}p_D(t,\cdot,\cdot)=p(t,\cdot,\cdot)$$ where $p(t,x,y)$ stands for the transition densities of the whole planar Brownian motion (i.e. not killed on the boundary of $D$).
\end{Lemma}
 
Let us now begin with the interpolation procedure.  By independence of $X_{t_0}$ and $X_{t_0,t}^{MFF}$, we can apply  \eqref{eqGFF1} to the function $$f_{\eqref{eqGFF1}}(x)=f(x)e^{t_0(\frac{\gamma}{\sqrt{2}}-\sqrt{2})^2} e^{\gamma X_{t_0}(x)-\frac{\gamma^2}{2} t_0 }$$ and get after a straightforward calculation involving \eqref{moveto}
\begin{align}
\label{eqGFF2}
\underset{t\to\infty}{\lim} &\E\Big[ \exp(- t^{\frac{3\gamma}{4}}e^{t(\frac{\gamma}{\sqrt{2}}-\sqrt{2})^2}\int_{\r^2}f(x) M^{\gamma,switch}_{t_0,t}(dx))\Big]\nonumber\\
&=\E\Big[\exp(-C (\gamma) \int_{\r^2}f(x)^\frac{2}{\gamma}e^{2X_{t_0}(x)-2t_0}M^{',MFF}_{t_0,\infty}(dx))\Big].
\end{align}

Let $\epsilon>0$ be fixed. From \eqref{unifgreen}, we can choose $T$ such that for all $T\leq t_0\leq t$, 
\begin{equation}
\label{eqGFF3} \sup_{x,y\in D'}|G_{D,t_0,t}(x,y)-G_{t_0,t}(x,y)|\leq \epsilon.
\end{equation}
Let us set $g_{t_0,t}(x)=e^{\frac{\gamma^2}{2}(\E[(X_{t}(x)-X_{t_0}(x))^2]-(t-t_0))} $. From \eqref{eqGFF3}, we have $e^{-\frac{\gamma^2}{2}\epsilon}\leq g_{t_0,t}(x)\leq e^{ \frac{\gamma^2}{2}\epsilon}$ for  all $T\leq t_0\leq t$. We will use this relation in the forthcoming lines. By Kahane's convexity inequalities and \eqref{eqGFF3}, we have for all $T\leq t_0\leq t$
\begin{align*}
\E\Big[ \exp(- t^{\frac{3\gamma}{4}}&e^{t(\frac{\gamma}{\sqrt{2}}-\sqrt{2})^2}\int_{\r^2}f(x)M^{\gamma}_{t}(dx))\Big]\\
 &\leq \E\Big[ \exp(- t^{\frac{3\gamma}{4}}e^{t(\frac{\gamma}{\sqrt{2}}-\sqrt{2})^2}e^{\epsilon^{1/2}Z-\epsilon/2}\int_{\r^2} f(x)g_{t_0,t}(x)M^{\gamma,switch}_{t_0,t}(dx))\Big]\\
 &\leq \E\Big[ \exp(- t^{\frac{3\gamma}{4}}e^{t(\frac{\gamma}{\sqrt{2}}-\sqrt{2})^2}e^{\epsilon^{1/2}Z-\epsilon/2}e^{-\frac{\gamma^2}{2}\epsilon}\int_{\r^2} f(x) M^{\gamma,switch}_{t_0,t}(dx))\Big]
\end{align*}
for some standard Gaussian random variable $Z$ independent of everything. We just explain some subtlety: observe that the definition of $M^\gamma_t$ does not involve a renormalization by the variance  $\E[X_t(x)^2]$ but $t$ instead. To apply KCI, one needs to compare measure involving a renormalization by the variance. So the function $g_{t_0,t}(x)$ appearing in the first inequality just results from the switching of variance required to apply KCI.

By taking the $\limsup$ as $t\to\infty$ in the above relation and by using \eqref{eqGFF2}, we deduce:
\begin{align*}
\limsup_{t\to\infty}\E&\Big[ \exp(- t^{\frac{3\gamma}{4}}e^{t(\frac{\gamma}{\sqrt{2}}-\sqrt{2})^2}\int_{\r^2}f(x)M^{\gamma}_{t}(dx))\Big]\\
&\leq \limsup_{t\to\infty}\E\Big[ \exp(- t^{\frac{3\gamma}{4}}e^{t(\frac{\gamma}{\sqrt{2}}-\sqrt{2})^2}e^{\epsilon^{1/2}Z-\epsilon/2}e^{-\frac{\gamma^2}{2}\epsilon}\int_{\r^2} f(x)M^{\gamma,switch}_{t_0,t}(dx))\Big]\\
&=\E\Big[\exp(-C (\gamma)e^{2\epsilon^{1/2}Z/\gamma-\epsilon/\gamma - \gamma \epsilon}\int_{\r^2} f(x)^\frac{2}{\gamma}e^{2X_{t_0}(x)-2 t_0}M^{',MFF}_{t_0,\infty}(dx))\Big].
\end{align*}
Now we want to apply once again KCI to the derivative martingale to replace the $M^{',MFF}_{t_0,\infty}$ part by $M^{',t_0}_\infty$. Recall that this is possible because we know that the Seneta-Heyde norming \cite{Rnew12} holds for both these measures. The control of covariance kernels is provided by \eqref{unifgreen} (notice that the uniform control w.r.t. $t$ is necessary to apply KCI for $t=\infty$) We get
\begin{align}
\limsup_{t\to\infty}\E&\Big[ \exp(- t^{\frac{3\gamma}{4}}e^{t(\frac{\gamma}{\sqrt{2}}-\sqrt{2})^2}\int_{\r^2}f(x)M^{\gamma}_{t}(dx))\Big]\nonumber\\
&\leq \E\Big[\exp(-C (\gamma)e^{2\epsilon^{1/2}Z/\gamma-\epsilon/\gamma - \gamma \epsilon+\epsilon^{1/2}Z'-\epsilon/2} 
\int_{\r^2} f(x)^\frac{2}{\gamma}e^{2X_{t_0}(x)-2t_0}M'_{t_0,\infty}(dx))\Big].\label{eqGFF4}
\end{align}
for some other standard Gaussian random variable $Z'$ independent of everything. By using the relation $e^{2X_{t_0}(x)-2\E[X_{t_0}]}M'_{t_0,\infty}(dx)=M'(dx)$, we see that \eqref{eqGFF4}
can be reformulated as 
\begin{align}
\limsup_{t\to\infty}\E&\Big[ \exp(- t^{\frac{3\gamma}{4}}e^{t(\frac{\gamma}{\sqrt{2}}-\sqrt{2})^2}\int_{\r^2}f(x)M^{\gamma}_{t}(dx))\Big]\nonumber\\
&\leq \E\Big[\exp(-C (\gamma)\,e^{2\epsilon^{1/2}Z/\gamma+\epsilon^{1/2}Z'-\epsilon(1\gamma + \gamma  +1/2)} 
\int_{\r^2} f(x)^\frac{2}{\gamma}e^{2\E[X_{t_0}]-2t_0}M' (dx))\Big].\label{eqGFF42}
\end{align}
By using the uniform convergence on $D'$ of $(\E[X_t(x)^2]-  t)_t$ as $t\to\infty$  towards the conformal radius  $\ln C(x,D)$ (see \cite[Lemma 6.1]{lacoin}), we deduce
\begin{align*}
\limsup_{t\to\infty}\E&\Big[ \exp(- t^{\frac{3\gamma}{4}}e^{t(\frac{\gamma}{\sqrt{2}}-\sqrt{2})^2}\int_{\r^2}f(x)M^{\gamma}_{t}(dx))\Big]\nonumber\\&\leq \E\Big[\exp(-C (\gamma)e^{2\epsilon^{1/2}Z/\gamma+\epsilon^{1/2}Z'-\epsilon(1\gamma + \gamma  +1/2)} \int_{\r^2}f(x)^\frac{2}{\gamma} C(x,D)^2M'(dx))\Big]. 
\end{align*}
 Since $\epsilon$ can be chosen arbitrarily small, the upper bound for the limit in Theorem \ref{mainGFF} when the cut off family $(X_t)_t$ has covariance $G_{D,t}$ is proved. The lower bound follows from a similar argument. Then we can use the same arguments as in the case of the MFF to extend the convergence to cut-off families uniformly close to $(G_{D,t})_t$.\qed
 
\vspace{2mm}
\noindent {\it Proof of Lemma \ref{cvunif}.} Recall the standard formula \cite[section 3.3]{Morters}
\begin{equation*}
\Delta(s,x,y):=p(s,x,y)-p_D(s,x,y)= \E^{x} \Big[ 1_{\{T^x_D  \leq s\}}  \frac{1}{2\pi(s-T^{x}_D)}e^{-\frac{|B^{x}_{T_D^{x}}-y |^2}{2(s-T^{x}_D)}}   \Big]
\end{equation*}
where $B^x_{t}$ is a standard Brownian motion starting from $x$ and $T^x_D=\text{inf} \{t \geq 0, \; B^x_t\not \in D \}$. If we denote   $\delta={\rm dist}(D',D^c)$, we deduce\qed
$$\Delta(s,x,y)\leq  \E^{x} \Big[ 1_{\{T^x_D  \leq s\}}  \frac{1}{2\pi(s-T^{x}_D)}e^{-\frac{\delta^2}{2(s-T^{x}_D)}}   \Big].$$
Now observe that the mapping $u\mapsto u e^{-u}$ is decreasing for $u\geq 1$. Therefore for $s\geq \delta^2/2$, we have
$$\Delta(s,x,y)\leq    \frac{1}{ \pi \delta^2 s}e^{-\frac{\delta^2}{2s}}   ,$$ which obviously completes the proof of the lemma.\qed

%%%%%%%%%%%%%%%%%%%%%%%%%%%%%%%%%%%%%%%%%%%%%%%%%%%%%%%%%%%%%%%%%%%%%%%%%%%%%%%%%%%%%%%%%%%%%%%%%%%%%%%%%%%%%%%%%%%%%%%%%%%%%%%%%%%%%%%%%%%%%%%%%%%%%%%%%%%%%%%%%%%%%%%%%%%%%%%%%%%%%%%%%%%%%%%%%%%%%%%%%%%%%%%%%%%%%%%%%%%%%%%%%%%%%%%%%%%%%%%%%%%%%%%%%%%%%%%%%%%%%%%%%
\appendix%%%%%%%%%%%%%%%%%%%%%%%%%%%%%%%%%%%%%%%%%%%%%%%%%%%%%%%%%%%%%%%%%%%%%%%%%%%%%%%%%%%%%%%%%%%%%%%%%%%%%%%%%%%%%%%%%%%%%%%%%%%%%%%%%%%%%%%%%%%%%%%%%%%%%%%%%%%%%%%%%%%%%%%%%%%%%%%%%%%%%%%%%%%%%%%%%%%%%%%%%%%%%%%%%%%%%%%%%%%%%%%%%%%%%%%%%%%%%%%%%%%%%%%%%%%%%%%%%%%%%%%%%
\section{Toolbox of technical results}\label{app:toolbox}
%%%%%%%%%%%%%%%%%%%%%%%%%%% 
In this subsection, we gather some results in \cite{Rnew7,Rnew12,Mad13} in order to have a paper  self contained, at least as much as possible.

We first recall a Lemma that can be found in \cite{Rnew7} 
\begin{Lemma} 
\label{rRhodesdec}
For any fixed $u\neq x$, the process $(Y_t(u))_{t\geq 0}$ can be decomposed as:
\begin{eqnarray*}
Y_t(u)=P_t^{x}(u)+Z_t^{x}(u) - \zeta_t^x(u),\qquad \forall t>0,
\end{eqnarray*}
where, for $t>0$

- $\zeta_t^{x}(u):=\sqrt{2\d} t-\sqrt{2\d} \int_0^{t}k(e^{s}(x-u))ds$, 

- $P_t^{x}(u):=\int_0^{t}k(e^{s}(x-u))dY_s(x)$ is measurable with respect to the $\sigma$-algebra generated by $(Y_t(x))_{t\geq 0}$,

-$(Z_t^{x}(u))_{t\geq 0} $ is a centered Gaussian process independent of $(Y_t(x))_{ t\geq 0}$  with covariance kernel:
\begin{equation}
 \label{covgeneralZ}\E\left( Z_t^x(u)  Z_{t'}^x(v) \right):= \int_0^{t \wedge t'} \left[ k(e^s(u-v))-k(e^s(x-u))k(e^{s}(x-v))\right]ds  ,\qquad \forall\, x,u,v\in \r^\d.
\end{equation}
\end{Lemma}

The following  lemma can be found in \cite{Mad13} (we refer to the Lemmas 3.1, 3.2 and 3.3) 
\begin{Lemma}
\label{lto0}
For any $\theta \in \r_{+}^* $ and $\epsilon >0$
\begin{align}
\label{eq0lto}
&\lim_{t,R\to\infty,\frac{e^t+1}{R+1}\in\N^*} \P\left(  |\gamma^{-1}\ln \theta|M_t^{\sqrt{2\d}}(e^{-t} [0,1]^\d ) \geq \epsilon \theta^{-\frac{\sqrt{2\d}}{\gamma}}  \right) + \P\left(   M_t'(e^{-l}{\rm BZ}_{R,t})\geq \epsilon \theta^{-\frac{\sqrt{2\d}}{\gamma}}  \right) \leq \epsilon,\\
\label{eqlto2}
&\lim_{t\to\infty}  \P\left(  w_{Y_t(\cdot)}^{1/3}(\frac{1}{t}e^{-t})\geq e^{\frac{t}{3}} \right) =\lim_{t\to\infty}\P\left(\underset{x,y\in [0,e^t]^\d,\, |x-y|\leq \frac{1}{t}}{\sup} \frac{\left|Y_t(\frac{x}{e^t})-Y_t(\frac{y}{e^t}) \right|}{|x-y|^{1/3}}\geq 1\right)=0,\\
\label{eq1lto}
&\lim_{t\to\infty}  \P\left(\forall x\in [0,1]^\d,\, -10\sqrt{2\d} \,t\leq Y_t(x)\leq -\kappa_\d\ln t\right)=1.
\end{align}
\end{Lemma}

In this section we will use the following two lemmas from \cite{Mad13}:
\begin{Lemma}
\label{maxfx}
We can find a constant $c_3>0$ such that for any $t'> 2$ and $R\geq 1$ such that $\frac{e^{t'}+1}{R+1}\in\N^* $
\begin{equation}
\label{eqmaxfx}
\P\Big( \sup_{ x\in [0,R]^\d} \sup_{s\in [\ln t',\infty)}Y_s(x)\geq \chi(x)\Big)\leq  c_3\int_{[0,R]^\d}((\ln t')^\frac{3}{8}+  \chi(x)^{\frac{3}{4}})e^{-\sqrt{2\d}\chi(x)}dx
\end{equation} 
for any  $\chi(\cdot)\in \mathcal{C}_R(t',\,10,\, +\infty)$.
\end{Lemma}

\begin{Lemma}
\label{LemL2}
We can find two constants $c_4,\, c_5>0 $ such that for any $t'\geq 2$, there exists $T(t')>0$ such that for any $L>0,\, R\geq 1$,  $\chi(\cdot)\in \mathcal{C}_R(t', 10, +\infty)$, $t\geq t'$ and  $  a\leq \frac{t}{2} $, 
\begin{eqnarray}
\nonumber
\P\Big(\exists x\in [0,R]^\d,\, \sup_{s\in [\ln t', \frac{t}{2}]}Y_s(x)\leq \chi(x),\, \sup_{s\in [\frac{t}{2},t-a]}Y_s(x) \leq a_t+\chi(x)+L-1,\, \sup_{[t-a,t]}Y_s(x)\in  a_t+\chi(x)+L+ [0,1]\Big)
\\
\label{eqLemL2} \leq c_4(1+a) e^{-c_5 L}\int_{[0,R]^\d} (\sqrt{\ln t'}+\chi(x)) e^{-\sqrt{2\d}\chi(x)}dx .
\end{eqnarray}
\end{Lemma}
 
\noindent {\bf Remark:} Lemma \ref{LemL2} is not explicitly stated in \cite{Mad13} but stems easily from (4.26) in \cite{Mad13}.

Here we reproduce \cite[Proposition 4.1]{Mad13}  
\begin{Lemma}\label{tightmad}
There exist two constants $c_1,c_2>0$ such that for any $t'\geq 2$, there exists $T>0$ such that for any $R\in [1,\ln t']$ and $t\geq T$ 
\begin{equation}\label{eq:psg1}
\P\Big(\exists x_0\in[0,R]^\d,   Y_t(x_0)\geq a_t + \chi(x_0)  \Big)\geq c_2\int_{[0,R]^\d}   \chi(x)  e^{-\sqrt{2\d}\chi(x)} dx
\end{equation}
for any function $\chi \in \mathcal{C}_R(t', \kappa_\d\ln t',\ln t)$.
\end{Lemma}

\bibliographystyle{plain}
 %%%%%%%%%%%%%%%%%%%%%%%%%%%%%%%%%%%%%%%%%%%%%%%%%%

\end{document}